\documentclass[12pt]{article}
\usepackage{amsthm,amsfonts,amssymb,amsmath}
\usepackage{enumerate}
\usepackage[colorlinks=true]{hyperref}

\newcommand{\BP}{{\mathbb {P}}}
\newcommand{\BQ}{{\mathbb {Q}}}
\newcommand{\BR}{{\mathbb {R}}}

\newcommand{\CCC}{{\mathcal {C}}}
\newcommand{\CD}{{\mathcal {D}}}
\newcommand{\CE}{{\mathcal {E}}}

\newcommand{\CJ}{{\mathcal {J}}}

\newcommand{\CL}{{\mathcal {L}}}
\newcommand{\CM}{{\mathcal {M}}}
\newcommand{\CN}{{\mathcal {N}}}
\newcommand{\CO}{{\mathcal {O}}}

\newcommand{\CU}{{\mathcal {U}}}
\newcommand{\CV}{{\mathcal {V}}}

\newcommand{\CX}{{\mathcal {X}}}

\newcommand{\OL}{{\overline{L}}}
\newcommand{\OM}{{\overline{M}}}
\newcommand{\ON}{{\overline{N}}}

\newcommand{\OQ}{{\overline{Q}}}

\newcommand{\OB}{{\overline{B}}}
\newcommand{\OK}{{\overline{K}}}

\renewcommand{\OE}{{\overline{E}}}

\newcommand{\OTheta}{{\overline{\Theta}}}

\newcommand{\Ar}{{\mathrm{Ar}}}

\newcommand{\Fal}{{\mathrm{Fal}}}

\newcommand{\FZ}{{\mathrm{FZ}}}

\newcommand{\an}{{\mathrm{an}}}

\newcommand{\Div}{{\mathrm{Div}}}

\renewcommand{\div}{{\mathrm{div}}}

\newcommand{\charr}{{\mathrm{char}}}

\newcommand{\Gal}{{\mathrm{Gal}}}

\newcommand{\GSp}{{\mathrm{GSp}}}

\newcommand{\Hom}{{\mathrm{Hom}}}

\renewcommand{\Im}{{\mathrm{Im}}}

\newcommand{\intb}{{\mathrm{int}}}

\newcommand{\Lie}{{\mathrm{Lie}}}

\newcommand{\ord}{{\mathrm{ord}}}

\newcommand{\rank}{{\mathrm{rank}}}

\newcommand{\Pic}{\mathrm{Pic}}

\newcommand{\Picc}{{\mathcal{P}\mathrm{ic}}}

\newcommand{\rmod}{{\mathrm{mod}}}

\DeclareMathOperator{\Spec}{Spec}

\newcommand{\vol}{{\mathrm{vol}}}

\newcommand{\wt}{\widetilde}
\newcommand{\wh}{\widehat}

\newcommand{\pair}[1]{\langle {#1} \rangle}

\newcommand{\ds}{\displaystyle}

\newcommand{\ol }{\overline}

\newcommand{\lra}{\longrightarrow}

\newcommand{\kkk}{Let $k$ be either $\ZZ$ or a field. }

\newcommand{\nef}{\mathrm{nef}}
\renewcommand{\vert}{\mathrm{vert}}
\newcommand{\eff}{\mathrm{eff}}

\newcommand{\id}{\mathrm{id}}

\newcommand{\CLL}{\overline{\mathcal L}}
\newcommand{\CMM}{\overline{\mathcal M}}
\newcommand{\CHH}{\overline{\mathcal H}}
\newcommand{\CNN}{\overline{\mathcal N}}

\newcommand{\CDD}{\overline{\mathcal D}}
\newcommand{\CEE}{\overline{\mathcal E}}

\newcommand{\CC}{\mathbb{C}}
\newcommand{\RR}{\mathbb{R}}
\newcommand{\ZZ}{\mathbb{Z}}
\newcommand{\QQ}{\mathbb{Q}}
\newcommand{\FF}{\mathbb{F}}

\newtheorem{thm}{Theorem}[section]

\newtheorem{lem}[thm]{Lemma}
\newtheorem{prop}[thm]{Proposition}

\theoremstyle{definition}
\newtheorem{definition}[thm]{Definition}

\theoremstyle{remark}
\newtheorem{remark}[thm]{Remark}

\begin{document}

\title{Arithmetic bigness and a uniform Bogomolov-type result}
\author{Xinyi Yuan}
\maketitle

\tableofcontents

\section{Introduction}

In this paper, we study some positivity properties in Arakelov geometry, and apply them to obtain a uniform Diophantine property on curves. 
More precisely, the paper achieves the following three goals.
\begin{enumerate}[(1)]
\item Extend Zhang's theory of admissible metrics from projective curves to  families of projective curves (cf. \S \ref{sec 2}). 
\item Prove the bigness of the admissible canonical bundle of the universal family over the moduli space of curves introduced in (1)  (cf. \S \ref{sec 3}). This is the main theorem of this paper.
\item Prove a uniform Bogomolov-type theorem for curves over global fields as a consequence of the bigness result in (2)  (cf. \S \ref{sec 4}).
\end{enumerate}
Our whole treatment is based on the recent theory of adelic line bundles of Yuan--Zhang \cite{YZ2}, a limit version of the intersection theory in algebraic geometry and the arithmetic intersection theory of Arakelov \cite{Ara} and Gillet--Soul\'e \cite{GS}.
Our proof of the bigness in (2) relies on the results of Zhang \cite{Zha1,Zha3}, Cinkir \cite{Cin1} and de Jong \cite{dJo3} on lower bounds of the self-intersection numbers of the admissible canonical bundles of curves over global fields.

Our uniform Bogomolov-type theorem generalizes and strengthens the uniform Bogomolov-type theorem of Dimitrov--Gao--Habegger \cite{DGH1} and K\"uhne \cite{Kuh}. 
This gives an alternative proof of their result on the uniform Mordell--Lang problem proposed by Mazur \cite{Maz}. 
Our approach does not use the o-minimality theory, and works over function fields of arbitrary characteristics, but it relies on more algebraic geometry and Arakelov geometry.

\subsection{Uniform Bogomolov-type result}

Let $C$ be a smooth projective curve of genus $g>1$ over $\overline\QQ$. 
In our convention, all curves are assumed to be geometrically connected over the base field.
The Bogomolov conjecture proved by Ullmo \cite{Ull} asserts that for any divisor $\alpha$ on $C$ of degree 1, there is a constant $c>0$ such that
$$
\#\{x\in C(\ol \QQ): \wh h(x-\alpha)\leq c\}<\infty.
$$ 
Here $\wh h:J(\overline\QQ)\to \RR$ denotes the N\'eron--Tate height function on the Jacobian variety $J$ of $C$. 
The proof of the loc. cit. is based on the celebrated equidistribution theorem of Szpiro--Ullmo--Zhang \cite{SUZ}. 

Recently, a uniform version of this theorem was proved by 
Dimitrov--Gao--Habegger \cite{DGH1} and K\"uhne \cite{Kuh}.
In fact, the new gap principle in \cite[Thm. 4.1]{Gao3}, as a combination of \cite[Prop. 7.1]{DGH1} and \cite[Thm. 3]{Kuh}, asserts that there are constants $c_1,c_2>0$ depending only on $g>1$ such that 
for any smooth projective curve $C$ over $\overline \QQ$ of genus $g$, and for any $y\in C(\overline \QQ)$,
$$
\#\{x\in C(\ol \QQ): \wh h(x-y)\leq c_1\, \max\{h_\Fal(C),1\}\} \leq c_2.
$$ 
Here $h_\Fal(C)=h_\Fal(J)$ denotes the stable Faltings height of the Jacobian variety $J$ of $C$. 
Note that the case $c_1=0$ is a uniform version of the Manin--Mumford conjecture for curves.
It is worth noting that DeMarco--Krieger--Ye \cite{DKY} previously proved the new gap principle in the case that $g=2$, $y$ is a Weierstrass point, and $C$ has a morphism of degree two to an elliptic curve.

The new gap principle has a significant consequence to the uniform Mordell--Lang problem proposed by Mazur \cite[p. 234]{Maz}. 
Recall that the Mordell conjecture was proved by Faltings \cite{Fal1}, and a different proof was given by Vojta \cite{Voj}. 
Vojta's proof was simplified and extended by Faltings \cite{Fal3} to prove the Mordell--Lang conjecture for subvarieties of abelian varieties, and was further simplified by Bombieri \cite{Bom} in the original case of curves. 
The proofs of \cite{Voj, Fal3, Bom} actually gave an upper bound of the number of points of large heights, which was further refined by de Diego \cite{dDi} and R\'emond \cite{Rem}.
Combining the upper bound with the new gap principle, we obtain the uniform bound on the number of rational points in \cite[Thm. 4]{Kuh}, which asserts that there is a constant $c>0$ depending only on $g>1$ such that 
for any smooth projective curve $C$ of genus $g$ over an algebraically closed field $F$ of characteristic 0, for any $y\in C(F)$, and for any subgroup $\Gamma \subset J(F)$ of finite rank,
$$
\#\{x\in C(F): x-y\in\Gamma\} \leq c^{1+\rank\, \Gamma}.
$$ 

Our first result is the following uniform version of the Bogomolov conjecture, which strengthens and generalizes the new gap principle of \cite{DGH1,Kuh}.

\begin{thm} [Theorem \ref{small points2}] \label{small points22}
Let $g>1$ be an integer. 
Then there are constants $c_1,c_2>0$ depending only on $g$ satisfying the following properties. 
Let $K$ be either a number field or a function field of one variable over a field $k$. 
Then for any geometrically integral, smooth and projective curve $C$ of genus $g$ over $K$, and for any 
line bundle $\alpha\in \Pic(C_{\ol K})$ of degree 1, with the extra assumption that $(C_{\ol K},\alpha)$ is non-isotrivial over $\overline k$ in the case that $K$ is a function field of one variable over a field $k$, 
one has
\small
$$
\#\left\{x\in C(\ol K): \wh h(x-\alpha)\leq c_1\big( \max\{h_\Fal(C),1\}+
\wh h((2g-2)\alpha-\omega_{C/K})\big) \right\} \leq c_2.
$$ 
\normalsize
\end{thm}

Here a \emph{number field} means a finite extension of $\QQ$, and a \emph{function field of one variable over $k$} means a finitely generated field extension of $k$ of transcendence degree 1.  
In the function field case, we say that $(C_{\ol K},\alpha)$ is \emph{non-isotrivial} over $\overline k$ if it is not isomorphic to the base change from $\overline k$ to $\overline K$ of any pair $(C_0,\alpha_0)$ consisting of a smooth projective curve $C_0$ over $\overline k$ and a line bundle $\alpha_0\in \Pic(C_{0})$ of degree 1. 
We refer to \S \ref{sec normalization} for the normalization of the heights involved in the theorem.

Our theorem is stronger and more general than the new gap principle of \cite{DGH1,Kuh} by the following aspects:
\begin{enumerate}[(1)]
\item it has an extra non-negative term $\wh h((2g-2)\alpha-\omega_{C/K})$ in the formula;
\item it allows $\alpha$ to be in $\Pic^1(C_{\overline K})$ instead of just in $C(\overline K)$;
\item it is valid for global fields of all characteristics instead of just number fields, and the constants $c_1$ and $c_2$ are uniform for all these fields. 
\end{enumerate}
Our proof is very different from that of \cite{DGH1,Kuh}, and the key ingredient is a bigness result of adelic line bundles on the universal curve.
We will come back to that in the next subsection.

In the function field case, assuming that $C_{\overline K}$ (instead of $(C_{\overline K},\alpha)$) is non-isotrivial, Looper--Silverman--Wilms \cite[Thm. 1.2]{LSW} proves a similar bound with explicit constants independently. By Theorem \ref{fiberwise11} below, 
their bound can be converted to our bound without the extra term $\wh h((2g-2)\alpha-\omega_{C/K})$. We refer to the explanation after Theorem \ref{fiberwise11} for more details on their result.

It is also worth noting that an explicit uniform bound on the Mordell--Lang problem over function fields of positive characteristic $p>0$ was previously obtained by Buium--Voloch \cite{BV} using jet schemes. 



\subsection{Potential bigness}

Our exposition is based on the theory of adelic line bundles on quasi-projective varieties of Yuan--Zhang \cite{YZ2}, which generalizes the projective case of Zhang \cite{Zha2}.
In \cite{YZ2}, the base ring $k$ is set to be either $\ZZ$ or an arbitrary field, and the main definitions give a notion of adelic line bundles on $X/k$ for any quasi-projective integral scheme $X$ over $k$.
Roughly speaking, an adelic line bundle on $X$ is a reasonable limit of a sequence of hermitian line bundles (or usual line bundles if $k$ is a field) on projective models of $X$ over $k$.
As a convention, we only require hermitian line bundles to have continuous metrics (instead of smooth metrics). 
The adelic line bundles have nice functorial properties, analytification properties, 
intersection theory, positivity properties, and volume theory.

\kkk
Let $S$ be a quasi-projective and flat normal integral scheme over $k$.
Let $\pi:X\to S$ be a smooth relative curve over $S$ of genus $g>1$, i.e.,
a smooth projective morphism whose fibers are curves of genus $g$. 
Denote by $J\to S$ the relative Jacobian scheme over $S$. 

Our N\'eron--Tate height is based on a canonically defined line bundle $\Theta$ on $J$ satisfying the following properties:
\begin{enumerate}[(1)]
\item
$\Theta$ is symmetric and rigidified along the identity section of $J\to S$;
\item 
the restriction of $\Theta$ to geometric fibers of $J\to S$ 
are algebraically equivalent to twice of theta divisors. 
\end{enumerate}
As a consequence, $\Theta$ is relatively ample. 
The construction of $\Theta$ is known to experts and reviewed in Definition \ref{qqq}. 

By \cite[Thm. 6.1.1]{YZ2}, there is a nef adelic line bundle $\OTheta$ on $J/k$ extending $\Theta$ such that $[2]^*\OTheta=4\OTheta$ in $\wh\Pic(J/k)$. 
Here $\wh\Pic(J/k)$ is the group of adelic line bundles on $J/k$ in the sense of \cite{YZ2}. 
In particular, an adelic line bundle is \emph{strongly nef} if it is the limit of nef hermitian (or usual) line bundles on projective models of $J$ over $k$ under the boundary topology, and a further mild limit process on strongly nef adelic line bundles gives the notion of \emph{nef} adelic line bundles.  
This generalizes the construction of Zhang \cite{Zha2} from the case $S=\Spec \ZZ$ to general base $S$. 

Let $\alpha$ be a line bundle on $X$ of degree 1 on fibers of $X\to S$. 
Then we have an immersion
$$
i_\alpha: X\lra J, \quad x\longmapsto x-\alpha.
$$
We obtain an adelic line bundle $i_\alpha^*\OTheta$ on $X$ by pull-back. 

On the other hand, consider the morphism 
$$\tau:J\times_SX \lra J\times_SJ, \quad (y,x)\longmapsto (y,y+(2g-2)x-\omega_{X/S}).$$
This morphism agrees with the $J$-morphism
$$
i_{\omega-Q}: X_J\lra J_J,\quad x\longmapsto (2g-2)x-(\omega_{X_J/J}-Q).
$$
Here we denote $X_J=J\times_S X$ and $J_J=J\times_S J$, viewed as $J$-schemes via the first projections $q_1:X_J\to J$ and $p_1:J_J\to J$. 
Here $J_J$ is canonically isomorphic to the Jacobian scheme of $X_J$ over $J$, and 
$Q$ is a universal line bundle on $J\times_S X$ in a suitable sense.
Denote by $\OTheta_J=p_2^*\OTheta$ the adelic line bundle on $J_J$ by the pull-back via the second projection $p_2:J\times_SJ\to J$.
Then we obtain an adelic line bundle $\tau^*(\OTheta_J)$ on $X_J$ by pull-back. 

Finally, we have the following bigness result, where many relevant terms will be explained after the statement. 

\begin{thm} [Theorem \ref{bigness7}] \label{bigness77}
Let $k$ be either $\ZZ$ or a field. Let $S$ be a quasi-projective and flat normal integral scheme over $k$, and let $\pi:X\to S$ be a smooth relative curve over $S$ of genus $g>1$. 
Assume that the family $\pi:X\to S$ has maximal variation.
Then the following hold:
\begin{enumerate}[(1)]
\item  The adelic line bundle $\pi_*\pair{i_\alpha^*\OTheta,i_\alpha^*\OTheta}$ is nef and big on $S$. 
\item  The adelic line bundle $q_{1*}\pair{\tau^*(\OTheta_J),\tau^*(\OTheta_J)}$ is nef and big on $J$. 
\end{enumerate}
\end{thm}

Here we say that the family $\pi:X\to S$ has \emph{maximal variation} if the moduli morphism $S\to M_{g,k}$ is generically finite, where $M_{g,k}$ denotes the coarse moduli scheme of smooth curves of genus $g$ over $k$.

Recall that a nef adelic line bundle is \emph{big} if its top self-intersection number
(as defined in \cite[Prop. 4.1.1]{YZ2}) is strictly positive.
The notation
$$\pi_*\pair{\cdot,\cdot}:\wh\Pic(X)\times \wh\Pic(X)\lra \wh\Pic(S)$$ 
denotes the Deligne pairing introduced in \cite[Thm. 4.1.3]{YZ2}.

Let us explain how the bigness result implies our uniform Bogomolov-type theorem in Theorem \ref{small points22}.
Consider the line bundle $\OL=\tau^*(\OTheta_J)$ on $X_J$. 
Over the $m$-fold fiber product
$$
(X_J)^m_{/J}= (X_J)\times_J\cdots \times_J (X_J),
$$ 
we have a nef adelic line bundle
$$
\OL_m=m_\boxtimes\OL
=p_1^*\OL+\cdots +p_m^*\OL.
$$
Denote $d=\dim J$. 
If $m\geq d$, expand the top self-intersection number 
$$
\OL_m^{d+m}
=(p_1^* \OL + p_2^* \OL + \cdots + p_m^* \OL)^{d+m}.
$$
The expansion includes the term
$$
(p_1^* \OL)^{2}  \cdots (p_d^* \OL)^{2}
\cdot (p_{d+1}^* \OL)  \cdots (p_m^* \OL)
=a^{m-d} (q_{1*}\pair{\tau^*(\OTheta_J),\tau^*(\OTheta_J)})^d>0.
$$
Here $a$ denotes the degree of $\OL$ on the generic fiber of $q_1:X_J\to J$,
and the top self-intersection number 
$(q_{1*}\pair{\tau^*(\OTheta_J),\tau^*(\OTheta_J)})^d$ on $J$ is strictly positive by Theorem \ref{bigness77}(2).

Therefore, $\OL_m$ is big on $(X_J)^m_{/J}$ for all integers $m\geq d$.
We say that $\OL$ is \emph{potentially big} on $X_J/J/k$ in this situation.
We remark that our notion of \emph{potential bigness} has some similarity with the notion of \emph{correlation}
of Caporaso--Harris--Mazur \cite{CHM}.

With the bigness result, our next key step to prove Theorem \ref{small points22} is to apply the height inequality (and some variants) in \cite[Thm. 5.3.5]{YZ2}.
To illustrate the idea, we only consider the case that $K$ is a number field (and thus $k=\ZZ$). 
Take $S$ to be a fine moduli space of smooth curves of genus $g$ over $\ZZ[1/N]$ with a suitable full level-$N$ structure, and take $X\to S$ to be the universal curve.
Fix an integer $m\geq d$.
As a consequence of the height inequality, for any adelic line bundle $\OM$ on $J/\ZZ$, 
there is a non-empty Zariski open subscheme $U$ of $(X_J)^m_{/J}$ such that 
$$
h_{\OL_m}(x)\geq \epsilon\, h_{\OM}(\pi_m(x)) ,\quad \forall x\in U(\ol \QQ).
$$ 
Here $\pi_m:(X_J)^m_{/J}\to J$ denotes the structure morphism.

For convenience, we sketch an idea to prove the height inequality taking advantage of the condition that $\OL_m$ is nef and big.
We can first reduce to the case that $\OM$ is nef by \cite[Lem. 5.1.6(1)]{YZ2}.
Then the key is \cite[Thm. 5.2.2(2)]{YZ2}, the adelic version of the bigness theorems of Siu \cite{Siu} and Yuan \cite{Yua}.
It implies that for any nef adelic line bundle $\OM$ on $J$, viewed as an adelic line bundle on 
$(X_J)^m_{/J}$ via pull-back,
 we have
$$
\wh\vol(\OL_m-\epsilon\OM)\geq \OL_m^{d+m}-(d+m)\OL_m^{d+m-1}\cdot \epsilon\OM>0
$$ 
for some rational number $\epsilon>0$.
Then some multiple of $\OL_m-\epsilon\OM$ has a nonzero effective section. 
Away from the zero locus $Z$ of this effective section, we have 
$$
h_{\OL_m-\epsilon\OM}(x)\geq  0,\quad \forall x\in ((X_J)^m_{/J}\setminus Z)(\ol \QQ).
$$

Once we have the height inequality, choose the adelic line bundle $\OM$ on $J$ by
$$
\OM= \OTheta + \ol\lambda_{S} +\CO(c).
$$
Here $\OTheta$ is as above, $\ol\lambda_{S}$ is the (hermitian) Hodge (line) bundle on $S$ associated to the abelian scheme $J\to S$, and  
$\CO(c)$ with $c>0$ is a hermitian line bundle on $\Spec\ZZ$ of arithmetic degree $c$. Here $\ol\lambda_{S}$ and $\CO(c)$ are viewed as adelic line bundles on $J/\ZZ$ by pull-back.
For any point $y\in J(\ol K)$ with image $s\in S(\ol K)$, we have
$$
h_{\OM}(y)
= 2\, \wh h(y)
+ h_{\rm Fal}(X_s) +c.
$$
This essentially dominates the term 
$$\max\{h_\Fal(C),1\}+
\wh h((2g-2)\alpha-\omega_{C/K})$$ 
in Theorem 
\ref{small points22}.

Finally, with an induction and some extra arguments, this eventually implies Theorem \ref{small points22}. This gives the idea of the proof of the theorem in the number field case.

The function field case is proved similarly by applying the above argument to moduli spaces over $k$ (instead of over $\ZZ$) and paying
special attention to definitions of height functions and non-isotriviality.
The uniformity of $(c_1,c_2)$ on $K$ is obtained by applying the above arguments for all $K$ to a single moduli space over $\ZZ$. 

There are two other consequences of the bigness of $\OL_m$ for $m\geq \dim J$. The first consequence is that in the arithmetic case $k=\ZZ$,
the morphism $(X_J)^m_{/J} \times_\ZZ\QQ\to (J_J)^m_{/J}\times_\ZZ\QQ$ for all integers $m\geq \dim J$ satisfies the relative Bogomolov conjecture proposed by \cite[Conj. 1.2]{DGH2}. 
This can be easily seen from the above height inequality. 
The second consequence of the bigness is that if $k$ is a field of any characteristic, 
the morphism $(X_J)^m_{/J} \to (J_J)^m_{/J}$ is non-degenerate for all integers $m\geq \dim J$, in terms of the definition of non-degeneracy in \cite[\S6.2.2]{YZ2}. 

The morphism $X^m_{/S} \to J^m_{/S}$ also satisfies similar properties for $m\geq \dim S$.
Moreover, by a similar method, we deduce that if $k$ is a field, and $\pi:X\to S$ has maximal variation, then the Faltings--Zhang morphism
$$
i_{\FZ,m}: X^{m+1}_{/S}\lra J^m_{/S} ,\quad
(x_0,\cdots, x_m)\longmapsto (x_1-x_0,\cdots, x_m-x_0)
$$ 
is non-degenerate for all $m\geq \dim S+1$. 
These non-degeneracy results generalize \cite[Thm. 1.2(i), Thm 1.2']{Gao2} to base fields $k$ of all characteristics. We refer to \S\ref{sec non-degeneracy} for more details.

The proof of the non-degeneracy by \cite{Gao2} is based on the mixed Ax--Schanuel theorem for the universal family of abelian varieties by \cite{Gao1}, and the proof of the latter uses the o-minimality theory.
Recently, Bl\'azquez-Sanz--Casale--Freitag--Nagloo \cite{BCFN} announced a different approach of the mixed Ax--Schanuel theorem without using the o-minimality theory. 
Both approaches work only in characteristic 0, while our result is valid for families of curves in all characteristics. 
The proof of the new gap principle in \cite{DGH1,Kuh} depends crucially on the non-degeneracy of the Faltings--Zhang map proved in \cite{Gao2}, while our proof of Theorem \ref{small points22} does not need the non-degeneracy, but deduces the non-degeneracy as a consequence of the potential bigness.

Another crucial tool of the proof of \cite{Kuh} is an equidistribution theorem, while  our crucial tool is the arithmetic bigness. 
For a more precise comparison, recall that the original Bogomolov conjecture was proved by Ullmo \cite{Ull} in terms of the equidistribution theorem of \cite{SUZ}, and a second proof in terms of bounding the self-intersection number of the admissible canonical bundle was obtained along the line of 
Zhang \cite{Zha1,Zha3}, Cinkir \cite{Cin1} and de Jong \cite{dJo3}. 
Then the treatment of \cite{Kuh} is a family version of that of \cite{SUZ,Ull}, and our treatment is a family version of that of \cite{Zha1, Zha3, Cin1, dJo3}.

\subsection{Admissible canonical bundle}

Let $K$ be a number field or a function field of one variable. 
Let $C$ be a smooth projective curve of genus $g>1$ over $K$.
Let $\omega_{C/K}$ be the relative dualizing sheaf, and $\Delta \subset C^2$ be the diagonal divisor.
By the construction of Zhang \cite{Zha1}, 
there are a canonical adelic line bundle $\overline\omega_{C/K,a}$ on $C$ extending $\omega_{C/K}$, and a canonical adelic line bundle $\overline\CO(\Delta)_a$ on $C^2$ extending $\CO(\Delta)$. 
The metrics of $\overline\omega_{C/K,a}$ and $\overline\CO(\Delta)_a$ at an archimedean place are the Arakelov metrics introduced by Arakelov \cite{Ara}, and the metrics of $\overline\omega_{C/K,a}$ and $\overline\CO(\Delta)_a$ at a non-archimedean place are the admissible metrics introduced by Zhang \cite{Zha1}.

The adelic line bundles $\overline\omega_{C/K,a}$ and $\overline\CO(\Delta)_a$
satisfy many nice properties.
For example, the canonical isomorphism $\CO(\Delta)|_{\Delta}\to \omega_{C/K}^\vee$
is an isometry. 
The curvatures of $\overline\CO(\Delta)_a$ on closed fibers of the projections $p_i:C^2\to C$ are proportional to the curvature of $\overline\omega_{C/K,a}$ on $C$ at all places. 

Note that \cite{Zha1} was written before \cite{Zha2}, but the treatment works without much difficulty in the terminology of \cite{Zha2} as we review in \S\ref{sec appendix}.
One goal of this paper is to introduce a family version of the admissible canonical bundle in the terminology of \cite{YZ2}. 

\kkk
Let $S$ be a quasi-projective normal scheme over $k$.
Let $\pi:X\to S$ be a smooth relative curve over $S$ of genus $g>1$.
Let $\omega_{X/S}$ be the relative dualizing sheaf, and $\Delta \subset X^2_{/S}$ be the diagonal divisor.
In Theorem \ref{admissible2}, we will introduce a canonical adelic line bundle 
$\overline\omega_{X/S,a}$ on $X$ extending $\omega_{X/S}$, and a canonical adelic line bundle $\overline\CO(\Delta)_a$ on $X^2_{/S}$ extending $\CO(\Delta)$. 
The extensions are uniquely determined by the properties that at any point 
$s\in S$, whose residue field is a number field if $k=\ZZ$ and is a function field of one variable over $k$ if $k$ is a field, the pull-backs 
$\overline\omega_{X/S,a}|_{X_s}$ 
and 
$\overline\CO(\Delta)_a|_{X_s^2}$
are canonically isomorphic to the admissible adelic line bundles 
$\overline\omega_{X_s/s,a}$ on $X_s$
and 
$\overline\CO(\Delta_s)_a$ on $X_s^2$
of \cite{Zha1}. 

While the construction of Zhang \cite{Zha1} uses graph theory, our construction
uses the canonical metrics of line bundles on the relative Jacobian scheme. 
In fact, we can define $\overline\omega_{X/S,a}$ by the formula
$$\overline\omega_{X/S,a}
=\frac{1}{4g(g-1)}i_\omega^*\OTheta+ \frac{1}{64g^2(g-1)^4}\pi^*\pi_*\pair{i_\omega^*\OTheta,i_\omega^*\OTheta}$$
in $\wh\Pic(X/k)_\QQ$.
Here $\OTheta$ is the adelic line bundle on the relative Jacobian $J$ as above, and $i_\omega$ is the canonical morphism
$$i_\omega:X\lra J,\quad
x\longmapsto (2g-2)x-\omega_{X/S}.$$
By the formula, we immediately see that $\overline\omega_{X/S,a}$ is nef on $X$.
We refer to Theorem \ref{isomorphism3} for more formulae of this type.

\subsection{Bigness of admissible canonical bundle}

The following is our main theorem of this paper.

\begin{thm}[Theorem \ref{bigness5}] \label{bigness55}
Let $k$ be either $\ZZ$ or a field. 
Let $S$ be a quasi-projective and flat normal integral scheme over $k$. 
Let $\pi:X\to S$ be a smooth relative curve over $S$ of genus $g>1$ with maximal variation.
Then the adelic line bundle $\overline\omega_{X/S,a}$ is nef and big on $X$.
\end{thm}

Once $\overline\omega_{X/S,a}$ is nef and big on $X$, the Deligne pairing 
$\pi_*\pair{\ol\omega_{X/S,a}\, ,\ol\omega_{X/S,a}}$ is nef and big on $S$ by a general fact of adelic line bundles. 
However, we will first prove the bigness of $\pi_*\pair{\ol\omega_{X/S,a}\, ,\ol\omega_{X/S,a}}$ and then use it to prove the theorem by some minor arguments.

As for Theorem \ref{bigness77}(2), it is easily obtained by the identity 
$$
q_{1*}\pair{\tau^*(\OTheta_J), \tau^*(\OTheta_J)}= 
16(g-1)^3\OTheta+16g(g-1)^3\pi_J^*\pi_*\pair{\ol\omega_{X/S,a}\, ,\ol\omega_{X/S,a}}
$$
in $\wh\Pic(J/k)$.
A similar formula gives Theorem \ref{bigness77}(1).

Now we sketch our proof of Theorem \ref{bigness55}.
As mentioned above, the key is to prove that 
$\pi_*\pair{\overline\omega_{X/S,a},\overline\omega_{X/S,a}}$ is nef and big on $S$. 
The process is to align the relevant results of \cite{Zha1, Zha3, Cin1, dJo3} into a family.

Let us first sketch the geometric case that $k$ is a field.
Replacing $S$ by a finite extension if necessary, we can assume that $\pi:X\to S$ has a {stable compactification} $\ol\pi:\ol X\to \ol S$, i.e. a projective variety 
$\ol S$ over $k$ with an open immersion $S\to \ol S$, a stable relative curve 
$\ol\pi:\ol X\to \ol S$ of genus $g$, and an open immersion $X\to \ol X$ compatible with the previous morphisms. 
Our proof takes the following steps.

\medskip\noindent\emph{Step 1}.
There is an effective adelic divisor $\ol E_S$ on $S$ such that
$$\pi_*\pair{\overline\omega_{X/S,a},\overline\omega_{X/S,a}}
= \overline\pi_*\pair{\omega_{\overline X/\overline S},\omega_{\overline X/\overline S}}
- \CO(\ol E_S)
$$
in $\wh\Pic(S)_\QQ$.
Here $\overline\pi_*\pair{\omega_{\overline X/\overline S},\omega_{\overline X/\overline S}}$ is viewed as an element of $\wh\Pic(S)_\QQ$ by the natural map $\Pic(\ol S)_\QQ\to \wh\Pic(S)_\QQ$.
The divisor $\ol E_S$ is in fact defined by the identity.
Then it has underlying divisor 
$0\in \Div(S)$, and thus is totally determined by its Green's function $g_{\ol E_S}$ on the Berkovich analytic space $S^\an$.
We have an ``explicit'' description of $g_{\ol E_S}$ in that its value at any discrete valuation $v$ of $k(S)/k$ is given by the $\epsilon$-invariant of the curve $X_{H_v}$ over the valuation field $H_v$ of $v$ defined by Zhang \cite{Zha1} in terms of graph theory. This determines $g_{\ol E_S}$ by continuity, and we say that $\ol E_S$ is the globalization of the $\epsilon$-invariant.

\medskip\noindent\emph{Step 2}.
The Noether formula gives
$$\overline\pi_*\pair{\omega_{\overline X/\overline S},\omega_{\overline X/\overline S}}
=12\lambda_{\overline S}- \CO(\Delta_{\overline S})
$$
in $\Pic(\ol S)_\QQ$.
Here $\lambda_{\overline S}=\det \ol\pi_*\omega_{\ol X/\ol S}$ is the Hodge line bundle of $\ol X$ over $\ol S$, and $\Delta_{\overline S}$ is the divisor of $\ol S$ with support equal to
$\ol S\setminus S$ measuring the singularity of $\ol X$ over $\ol S$.
Both $\lambda_{\overline S}$ and $\Delta_{\overline S}$ are the well-known tautological divisors in the theory of moduli spaces of curves.

\medskip\noindent\emph{Step 3}.
The difference $(2g-2)\Delta_{\ol S}-\ol E_S$ is an effective adelic divisor in $\wh\Div(S)$.
As both $\Delta_{\ol S}$ and $\ol E_S$  have underlying divisor 
$0\in \Div(S)$, it suffices to check $(2g-2)g_{\Delta_{\ol S}}\geq g_{\ol E_S}$ on the Berkovich analytic space $S^\an$.
By continuity, we only need to check it at any discrete valuation $v$ of $k(S)/k$, or equivalently compare the $\epsilon$-invariant of the curve $X_{H_v}$ defined by Zhang \cite{Zha1} 
and the classical $\delta$-invariant of $X_{H_v}$ counting the number of nodes of reduction. 
The comparison is done by graph theory. 

Combining these steps, we have
$$\pi_*\pair{\overline\omega_{X/S,a},\overline\omega_{X/S,a}}
= 12\lambda_{\overline S}- (2g-1)\CO(\Delta_{\overline S})+\eff.
$$
Here ``$\eff$'' means an effective adelic divisor, and later we will write
``$\nef$'' for a nef adelic line bundle. 

\medskip\noindent\emph{Step 4}.
There is an effective adelic divisor $\ol \Phi_S$ on $S$ such that
$$
\pi_*\pair{\overline\omega_{X/S,a},\overline\omega_{X/S,a}}
= \frac{2}{3g-1}\CO(\ol \Phi_S)+\nef.
$$
Similarly to $\ol E_S$, the arithmetic divisor $\ol \Phi_S$ has underlying divisor 
$0\in \Div(S)$, and its Green's function $g_{\ol \Phi_S}$ at any discrete valuation $v$ of 
$k(S)/k$ is given by the $\varphi$-invariant of the curve $X_{H_v}$ defined by Zhang \cite{Zha3} in terms of graph theory. 
Then $\ol \Phi_S$ is the globalization of Zhang's $\varphi$-invariant.
The above equality is a family version of de Jong \cite[Thm. 8.1]{dJo3}. 

\medskip\noindent\emph{Step 5}.
The difference $\overline\Phi_S-\frac{1}{39}\Delta_{\ol S}$ is an effective adelic divisor in $\wh\Div(S)$.
Similar to Step 3, it suffices to compare the values of the Green's function at all discrete valuations $v$ of $k(S)/k$. 
Then it follows from the result of Cinkir \cite[Thm. 2.11]{Cin1} in graph theory.

As a consequence, we have
$$
\pi_*\pair{\overline\omega_{X/S,a},\overline\omega_{X/S,a}}
= \frac{2}{39(3g-1)}\CO(\Delta_{\ol S})+\nef+\eff.
$$

\medskip\noindent\emph{Step 6}.
Take a linear combination of the equalities respectively at the end of Step 3 and Step 5 to cancel the term $\CO(\Delta_{\ol S})$. 
We obtain
$$
\left(1+ \frac{39(3g-1)}{2} (2g-1) \right)
\pi_*\pair{\overline\omega_{X/S,a},\overline\omega_{X/S,a}}
= 12\lambda_{\overline S}+\nef+\eff.
$$
Then the bigness of 
$\pi_*\pair{\overline\omega_{X/S,a},\overline\omega_{X/S,a}}$ follows from the classical result that the Hodge (line) bundle 
$\lambda_{\overline S}$
is nef and big on $\overline S$. 

Now we prove  the bigness in the arithmetic case $k=\ZZ$. 
The above proof is not valid in this case as Step 2 becomes rather subtle in the arithmetic case. 
However, we already know the bigness of the image of $\pi_*\pair{\overline\omega_{X/S,a},\overline\omega_{X/S,a}}$
in $\wh\Pic(S_\QQ/\QQ)$ by the above proof.
Then it suffices to prove that 
$$\pi_*\pair{\overline\omega_{X/S,a},\overline\omega_{X/S,a}}= \CO(c)+\eff+\nef$$
for some rational number $c>0$. 
Here $\CO(c)$ is the trivial line bundle on $\ZZ$ with metric given by $\|1\|=e^{-c}$, and viewed as an adelic line bundle on $S$ by pull-back. 

Step 4 still works in the arithmetic case, so it suffices to prove that there exists a constant $c_0>0$ depending only on $g$ such that
$\ol\Phi_S-c_0$ is an effective adelic divisor on $S$.
It suffices to prove that there exists $c_0>0$ depending only on $g$ such that
 $\varphi(X_s)\geq c_0$ for any $s\in S(\CC)$.
Let $M_g$ be a fine moduli space of curves of genus $g$ over $\CC$ with a suitable level structure, and $\overline M_g$ be a suitable compactification of $M_g$. 
By Zhang \cite{Zha3}, we know that $\varphi>0$ on $M_g$. 
To prove that $\inf \varphi>0$ on $M_g$, we prove that $\varphi$ tends to infinity along the boundary $\overline M_g\setminus M_g$.
This is a consequence of Step 5 for $k=\CC$, as the image of $\ol\Phi_S$ in $\Div(M_g/\CC)$ governs its growth behavior along $\overline M_g\setminus M_g$.

\subsection{Uniform fiberwise bigness} \label{sec intro pointwise}

By Zhang \cite{Zha1}, the Bogomolov conjecture for a single smooth projective curve $C$ over a number field $K$ (or a function field of one variable) is equivalent to $\ol\omega_{C/K,a}^2>0$. The following theorem is a uniform version of this positivity and also a fiberwise version of Theorem \ref{bigness55}.  
We thank Ziyang Gao for raising the question of proving such a bound.

\begin{thm} [Theorem \ref{fiberwise2}, Theorem \ref{fiberwise1}] \label{fiberwise11}
Let $g>1$ be an integer. 
Then there are constant $c_3>0$ and $c_4>0$ depending only on $g$ satisfying the following properties. 
Let $K$ be either a number field or a function field of one variable over a field $k$. 
Then for any geometrically connected, smooth and projective curve $C$ over $K$ of genus $g$, with the extra assumption that $C_\OK$ is non-isotrivial over $\overline k$ in the function field case, 
one has
$$
c_3 \cdot \max\{h_\Fal(C),1\} \leq  [\ol\omega_{C/K,a}^2] \leq  c_4 \cdot \max\{h_\Fal(C),1\}.
$$ 
\end{thm}

The stable Faltings height $h_\Fal(C)$ is normalized as in \S \ref{sec normalization}.
The term
$\ol\omega_{C/K,a}$ is Zhang's admissible canonical bundle of $C$ over $K$, and the normalized intersection number
 $[\ol\omega_{C/K,a}^2]$ is as follows.
If $K$ is a number field, 
$$
[\ol\omega_{C,a}^2]= \frac{1}{[K:\QQ]} \ol\omega_{C/K,a}^2,
$$
where $\ol\omega_{C/K,a}^2$ is the arithmetic intersection number.
If $K$ is a function field of one variable over $k$,
$$
[\ol\omega_{C,a}^2]= \ol\omega_{C/K,a}^2,
$$
where the intersection number $\ol\omega_{C/K,a}^2$ on $C$ is 
normalized by multiplicity functions given by degrees over $k$.

The theorem implies that $\max\{h_\Fal(C),1\}$ and $[\ol\omega_{C/K,a}^2]$ are equivalent invariants in both the uniform Bogomolov-type problem and the conjectural arithmetic Bogomolov--Miyaoka--Yau inequality.  

Most pieces of the theorem are known to experts previously, while the new piece here is only the first inequality in the number field case.
Our new ingredient is still to use adelic divisors to bound asymptotic behavior of relevant archimedean invariants near the boundary of the moduli space of curves. 
We will see that the theorem is a consequence of Theorem \ref{bigness55}, but to get more explicit constants $c_3,c_4$, we will mainly work on single curves by refining 
 many ingredients and strategies of the proof of 
Theorem \ref{bigness55}. Roughly speaking, our treatment still follows the line of \cite{Zha1, Zha3, Cin1, dJo3}, but we also use ideas and results of \cite{Wil2, Wil3} to get more optimal constants.  
We refer to Theorem \ref{fiberwise2} and Theorem \ref{fiberwise1} for explicit constants,  proofs, and more historical accounts. 

As mentioned above, Looper--Silverman--Wilms \cite{LSW} proves a uniform Bogomolov-type theorem over function fields, which is closely related to Theorem \ref{small points22}.
Namely, in the setting of Theorem \ref{small points22} and Theorem \ref{fiberwise11}, 
in the function field case, assuming that $C_{\ol K}$ is non-isotrivial over $\bar k$, 
\cite[Thm. 1.2]{LSW} takes the form
$$
\#\left\{x\in C(\ol K): \wh h(x-\alpha)\leq c_1'\, \ol\omega_{C,a}^2 \right\} \leq c_2'.
$$ 
The constants $(c_1',c_2')$ depend only on $g$ and are amazingly explicit. 
With Theorem \ref{fiberwise11}, we can replace the term $\ol\omega_{C,a}^2$ by 
$\max\{h_\Fal(C),1\}$ and get a bound closer to that in Theorem \ref{small points22}. 
The proof in the loc. cit. is also based on Zhang's admissible pairings, but its arguments are on individual curves instead of on moduli spaces of curves.

\subsection{Notations and terminology}\label{sec notation}

By a \emph{variety} over a field $k$, we mean an integral scheme, separated of finite type over $k$. 

By a \emph{curve} over a field $k$, we mean a geometrically connected projective scheme $C$ of dimension 1 over a field $k$.
The curve is called \emph{semistable} if $C_{\bar k}$ is reduced and all singular points of $C_{\bar k}$ are ordinary double points.  
The curve is called \emph{stable} if it is semistable with arithmetic genus $g>1$, and any rational irreducible component of $C_{\bar k}$ intersects other irreducible components at three or more points. 

By a \emph{relative curve over a scheme $S$}, we mean a projective and flat morphism $\pi:X\to S$ of purely relative dimension 1 with geometrically connected fibers. 
The curve is called \emph{smooth (resp. stable, semistable)}, if every fiber of  
$\pi:X\to S$ is smooth (resp. stable, semistable).
It is said to be \emph{of genus $g$}, if every fiber of $\pi:X\to S$ has arithmetic genus $g$.

Let $S$ be a quasi-projective scheme over a base ring $R$, and let $\pi:X\to S$ be a smooth relative curve of genus $g>1$. 
We say that $\pi:X\to S$ has \emph{maximal variation} (or \emph{maximal variation of moduli}) if the moduli morphism $S\to M_{g,R}$ associated to $\pi$ is generically finite. Here $M_{g,R}$ denotes the coarse moduli scheme of smooth curves of genus $g$ over $R$.

By a \emph{line bundle} on a scheme, we mean an invertible sheaf on the scheme. 
We often write or mention tensor products of line bundles additively, so $aL-bM$ means
$L^{\otimes a}\otimes M^{\otimes (-b)}$
for line bundles $L,M$ and integers $a,b$.

Let $X\to S$ be a morphism of schemes. For $m>0$, we denote by
$$
X_{/S}^m=X\times_SX\times_S\cdots \times_S X
$$
the $m$-fold fiber product of $X$ over $S$. We take the convention $X_{/S}^0=S$.
If $L$ is a line bundle on $X$, then \emph{the box tensor}
$$
m_{\boxtimes}L=p_1^*L+p_2^*L+\cdots+p_m^*L
$$
is a line bundle on $X_{/S}^m$. Here $p_i:X_{/S}^m\to X$ is the projection to the $i$-th component.

To treat adelic line bundles of \cite{YZ2}, we follow the uniform terminology of \cite[\S1.5]{YZ2}.
In particular, most of the time, our base ring $k$ is $\ZZ$ or a field.

\subsubsection*{Acknowledgments}
A major motivation of this paper is to understand the recent works of Dimitrov--Gao--Habegger and K\"uhne on the uniform Mordell--Lang conjecture. This paper is certainly inspired by many of their ideas, though the eventual approach is very different.

The author is indebted to Shou-Wu Zhang for many crucial discussions, and the paper grows out of the joint work with him on the theory of adelic line bundles on quasi-projective varieties.

The author would like to thank Alexandru Buium, Jan-Hendrik Evertse, Ziyang Gao, Rafael von K\"anel, Yanshuai Qin, Jos\'e Felipe Voloch, Robert Wilms, and Wei Zhang for many 
helpful communications.
The author would like to thank Ko Aoki, Ruoyi Guo, Qirui Li, Shiquan Li, Cong Wen, Yinchong Song, Jiaming Zhang, and Xiaoxiang Zhou for pointing out many mistakes and typos in earlier versions of the paper.

The author is also grateful to the anonymous referees for pointing out many mistakes and making many great suggestions for this paper.

\section{Admissible adelic line bundles} \label{sec 2}

The goal of this section is to generalize the theory of admissible line bundles of Zhang \cite{Zha1} from curves to relative curves using the terminology of Yuan--Zhang \cite{YZ2}. 
Note that the original theory of \cite{Zha1} is reviewed in our appendix \S\ref{sec appendix}.
In \S\ref{sec adelic}, we review the theory of adelic line bundles of \cite{YZ2}. 
In \S\ref{sec admissible}, we introduce our theory of admissible line bundles on relative curves.
In \S\ref{sec pullback}, we consider some formulae on intersections and pull-back of admissible adelic line bundles.

\subsection{Review on adelic line bundles} \label{sec adelic}

The goal of this subsection is to sketch the notion of adelic line bundles on quasi-projective varieties of Yuan--Zhang \cite{YZ2}, which generalizes the more classical adelic line bundles of projective varieties over number fields of Zhang \cite{Zha2}. 
In the end, we prove two basic results about adelic line bundles which will be used in this paper.

\subsubsection*{Adelic divisors}

\kkk
Let $\CU$ be a flat and quasi-projective integral scheme over $k$.
Let us first recall the definition of adelic divisors on $\CU/k$.  

Let $\CX$ be a \emph{projective model} of $\CU$ over $k$, i.e., a projective integral scheme over $k$ with an open immersion $\CU\to\CX$ over $k$.
In the spirit of \cite[\S2.2]{YZ2}, take the fiber product
$$
\wh\Div(\CX,\CU)=\wh\Div(\CX)_\QQ \times_{\Div(\CU)_\QQ} \Div(\CU),
$$
whose elements are arithmetic divisors of mixed coefficients. 
In the arithmetic case (that $k=\ZZ$), $\wh\Div(\CX)$ is the group of arithmetic divisors on $\CX$, where the Green's functions are assumed to be continuous (away from the singularities). 
In the geometric case (when $k$ is a field), $\wh\Div(\CX)$ means the usual 
$\Div(\CX)$.

By abuse of terminology, in the following, if $k$ is a field, then ``arithmetic divisor'' (resp. ``hermitian line bundle'') means
``divisor'' (resp. ``line bundle'').

Define the group of \emph{model adelic divisors} by
$$\wh\Div (\CU/k)_\rmod=\lim_{\substack{\lra\\ \CX}}\wh\Div(\CX,\CU).$$
Here the limit is over the system of projective models $\CX$ of $\CU$ over $k$.
An element of $\wh\Div (\CU/k)_{\rmod}$ is \emph{effective} if it is the image of an effective element of some $\wh\Div(\CX,\CU)$, where an element of $\wh\Div(\CX,\CU)$ is \emph{effective} if its images in $\wh\Div(\CX)_\QQ$ and $\Div(\CU)$ are both effective. 

Fix a \emph{boundary divisor} $(\CX_0,\CEE_0)$ of $\CU$, i.e,  a projective model $\CX_0$ of $\CU$ and a strictly effective arithmetic divisor $\CEE_0$ on $\CX_0$ such that the support of the finite part $\CE_0$ is exactly $\CX_0\setminus \CU$.
We have a \emph{boundary norm} 
$$\|\cdot\|_{\CEE_0}:\wh\Div (\CU/k)_\rmod
\lra [0,\infty]$$
by 
$$
\|\CDD\|_{\CEE_0}:=\inf\{\epsilon\in \BQ_{>0}: \ 
 -\epsilon \CEE_0 \leq
\CDD \leq  \epsilon \CEE_0\}.
$$
Here the inequalities are defined in terms of effectivity. 
It further induces a \emph{boundary topology} on $\wh\Div (\CU/k)_{\rmod}$, which does not depend on the choice of $(\CX_0,\overline\CE_0)$.

Let $\wh \Div  (\CU/k)$ be the \emph{completion} of $\wh \Div  (\CU/k)_{\rmod}$ with respect to the boundary topology. 
An element of $\wh \Div(\CU/k)$ is called an \emph{adelic divisor} on $\CU/k$.

By definition, an adelic divisor is represented by a Cauchy sequence in $\wh \Div  (\CU/k)_\rmod$, i.e., a sequence $\{\CDD_i\}_{i\geq 1}$ in $\wh \Div  (\CU/k)_\rmod$ satisfying the property that there is a sequence $\{\epsilon_i\}_{i\geq 1}$ of positive rational numbers converging to $0$ such that 
$$
 -\epsilon_i \CEE_0 \leq
\CDD_{i'}-\CDD_{i} \leq  \epsilon_i \CEE_0,\quad\ i'\geq i\geq 1.
$$

There is a canonical map
$$
\wh\Div(\CU/k) \lra \Div(\CU),\quad
\{\CDD_i\}_{i\geq 1}\longmapsto\CD_1|_{\CU}.
$$ 
We usually write $\CDD=\{\CDD_i\}_{i\geq 1}$ and $\CD=\CD_1|_{\CU}$, and call $\CD$ the \emph{underlying divisor} of $\CDD$.

\subsubsection*{Adelic line bundles}

\kkk
Let $\CU$ be a flat and quasi-projective integral scheme over $k$. 

Let $\CX$ be a projective model of $\CU$ over $k$.
In the spirit of \cite[\S2.2]{YZ2}, 
let $\wh\Picc(\CX)$ be the category of hermitian line bundles on $\CX$, 
and $\wh\Picc(\CX)_\QQ$ be the category of hermitian $\QQ$-line bundles on $\CX$. 
In the arithmetic case (when $k=\ZZ$), $\wh\Picc(\CX)$ is the category of hermitian line bundles with continuous metrics on $\CX$. 
In the geometric case (when $k$ is a field), $\wh\Picc(\CX)$ means the usual 
$\Picc(\CX)$.

As a convention, categories of various line bundles are defined to be groupoids; i.e., the morphisms in them are defined to be isomorphisms (or isometries) of the line bundles. 
To illustrate the category of various $\QQ$-line bundles, take $\wh\Picc(\CX)_\QQ$ for example. 
An object of $\wh\Picc(\CX)_\QQ$ is a pair $(a,\CL)$ (or just written as $a\CL$)
with $a\in \QQ$ and $\CL\in\wh\Picc(\CX)$, and a morphism of two such objects is defined to be an element of
$$\Hom(a\CL,a'\CL')=\varinjlim_m \Hom(am\CL, a'm\CL'),$$
where $m$ runs through positive integers such that $am$ and $a'm$ are both integers, so that $am\CL$ and $a'm\CL'$ are viewed as objects of $\wh\Picc(\CX)$, 
and ``$\Hom$'' on the right-hand side are viewed in $\wh\Picc(\CX)$.

Let $(\CX_0,\CEE_0)$ be a boundary divisor as above.
Define the \emph{category $\wh\Picc (\CU/k)$ of adelic line bundles} on $\CU/k$ as follows.
An object of $\wh\Picc (\CU/k)$ is a pair
$(\CL, (\CX_i,\overline \CL_i, \ell_{i})_{i\geq 1})$ where:
\begin{enumerate}[(1)]
\item $\CL$ is an object of $\Picc(\CU)$, i.e., a line bundle on $\CU$;

\item  $\CX_i$ is a projective model of $\CU$ over $k$;

\item  $\overline \CL_i$ is an object of $\wh\Picc(\CX_i)_\QQ$, i.e. a hermitian $\QQ$-line bundle on $\CX_i$;

\item $\ell_i:\CL\to \CL_i|_{\CU}$ is an isomorphism in $\Picc(\CU)_\QQ$, where $\CL_i$ is the underlying $\QQ$-line bundle of $\CLL_i$ on $\CX_i$.
\end{enumerate}
The sequence is required to satisfy the \emph{Cauchy condition} that
the sequence $\{\wh \div(\ell_i \ell_1^{-1})\}_{i\geq 1}$ is a 
Cauchy sequence in $\wh\Div(\CU/k)_\rmod$ under the boundary topology.
Here $\ell_i \ell_1^{-1}: \CL_1|_{\CU} \to \CL_i|_{\CU}$ is viewed as a rational section of the underlying line bundle of
$\CLL_i-\CLL_1$, so that $\wh \div(\ell_i \ell_1^{-1})$ is a well-defined element of $\wh\Div(\CU/k)_\rmod$. 

A morphism from an object $(\CL, (\CX_i,\overline \CL_i, \ell_{i})_{i\geq 1})$ to another 
$(\CL',(\CX_i',\overline \CL_i', \ell_{i}')_{i\geq 1})$ is an isomorphism  $\iota:\CL\to \CL'$ of the integral line bundles on $\CU$ such that the sequence 
$\{ \wh\div(\ell_i'\iota \ell_i^{-1}) \}_{i\geq1}$
of $\wh \Div (\CU/k)_\rmod$ converges to 0 in $\wh \Div (\CU/k)$
under the boundary topology.
Note that the model arithmetic divisor $\wh\div(\ell_i'\iota \ell_i^{-1})$ is equal to
$\wh\div(\ell_{i}'\ell_1'^{-1})-\wh\div(\ell_{i}\ell_1^{-1})+\wh\div(\ell_1'\iota \ell_1^{-1})
$.

An object of $\wh\Picc (\CU/k)$ is called an \emph{adelic line bundle} on $\CU$.
Define $\wh\Pic (\CU/k)$ to be the \emph{group} of isomorphism classes of objects of $\wh\Picc (\CU/k)$. 

There is a canonical forgetful functor
$$
\wh\Picc(\CU/k) \lra \Picc(\CU),\quad
(\CL, (\CX_i,\overline \CL_i, \ell_{i})_{i\geq 1})\longmapsto \CL
$$ 
We usually write $\CLL=(\CL, (\CX_i,\overline \CL_i, \ell_{i})_{i\geq 1})$ and call $\CL$ the \emph{underlying line bundle} of $\CLL$.

There is a canonical surjection
$$\wh\Div(\CU/k)\lra \wh\Pic(\CU/k), \quad
\CDD\longmapsto \CO(\CDD).$$
The kernel is the image of the group of principal arithmetic divisors on projective models of $\CU$.

An adelic  divisor is called \emph{effective} if it is equal to a limit of effective arithmetic divisors (of mixed coefficients). 
An adelic line bundle is called \emph{effective} if it is the image of an effective adelic divisor. 
An adelic line bundle (resp. adelic divisor) is called \emph{strongly nef} if it is isomorphic (resp. equal) to a limit of nef hermitian line bundles (resp. adelic divisors) under the boundary topology. 
An adelic line bundle  (resp. adelic divisor) $\CLL$ on $\CU$ is \emph{nef} if there exists a strongly nef adelic line bundle (resp. adelic divisor) $\CMM$ on $\CU$ such that $a\CLL+\CMM$ is strongly nef for all positive integers $a$.
An adelic line bundle (resp. adelic divisor) is \emph{integrable} if it is 
isomorphic (resp. equal) to the difference of two strongly nef ones.

Denote by $\wh\Picc(\CU/k)_\nef$ (resp. $\wh\Picc(\CU/k)_\intb$) the subcategory of 
$\wh\Picc(\CU/k)$ consisting of nef (resp. integrable) adelic line bundles.
Denote by $\wh\Pic(\CU/k)_\nef$ (resp. $\wh\Pic(\CU/k)_\intb$) the subsets of nef (resp. integrable) elements in $\wh\Pic(\CU/k)$. 

If the base ring $k$ is clear, we usually omit the dependence on $k$ of the groups or categories of adelic objects. For example,  
$\wh\Div(\CU/k)$, $\wh\Picc(\CU/k)$, $\wh\Pic(\CU/k)$ are written as 
$\wh\Div(\CU)$, $\wh\Picc(\CU)$, $\wh\Pic(\CU)$.

\subsubsection*{Analytification}

\kkk
Let $\CU$ be a flat and quasi-projective integral scheme over $k$. 
There is a natural \emph{Berkovich analytic space} $\CU^\an=(\CU/k)^\an$.
In fact, if $\CU$ has an open affine cover $\{\Spec A_i\}_i$, then $\CU^\an=\cup_i \CM(A_i)$, where $\CM(A_i)$ is the set of multiplicative semi-norms on $A_i$, assumed to be trivial over $k$ if $k$ is a field.
The set $\CU^\an$ is endowed with the coarsest topology such that every $\CM(A_i)$ is open in $\CU^\an$ and that $|f|:\CM(A_i)\to \RR$ is continuous for all $f\in A_i$.
Then $\CU^\an$ is Hausdorff, path-connected and locally compact.

An \emph{arithmetic divisor on $\CU^\an$} is a pair $\CDD=(\CD, g_\CD)$ consisting of a divisor $\CD$ on $\CU$ and a Green's function $g_\CD$ of $\CD$ on $\CU^\an$, i.e., 
 a continuous function $g_\CD:\CU^\an\setminus |\CD|^\an \to \RR$ with logarithmic singularity along $\CD$ in the sense that, for any rational function $f$ on a Zariski open subset $\CV$ of $\CU$ satisfying $\div(f)=\CD|_\CV$, the function
$g_\CD+\log |f|$ can be extended to a continuous function on $\CV^\an$.
Note that the definition uses the term $g_\CD+\log |f|$ instead of $g_\CD+2\log |f|$, so our convention is different from that used by Gillet--Soul\'e.

A \emph{metrized line bundle on $\CU^\an$} is a pair $\CLL=(\CL,\|\cdot\|)$ consisting of a line bundle $\CL$ on $\CU$ and a continuous metric $\|\cdot\|$ of $\CL$ on $\CU^\an$; i.e., an assignment of a $H_v$-metric $\|\cdot\|_v$ to the fiber $\CL\otimes H_v$ over every $v\in \CU^\an$, assumed to vary continuously as $v$ varies.
Here $H_v$ denotes the completed residue field of $v$. 

Denote by $\wh \Div (\CU^\an)$ the group of arithmetic divisors on $\CU^\an$.
Denote by $\wh \Picc (\CU^\an)$ (resp. $\wh \Pic (\CU^\an)$) the category (resp. group of isometry classes) of metrized line bundles on $\CU^\an$.

By \cite[Prop. 3.3.1, Prop. 3.4.1]{YZ2}, there are injective analytification maps
$$
\wh\Div(\CU)\lra \wh\Div(\CU^\an),
$$
$$
\wh\Pic(\CU)\lra \wh\Pic(\CU^\an),
$$
and a fully faithful analytification functor
$$
\wh\Picc(\CU)\lra \wh\Picc(\CU^\an).
$$
These maps will bring lots of convenience in our treatment later.

\subsubsection*{Intersection theory}

Here we recall the intersection theory in \cite[\S4.1]{YZ2}. 

\kkk
Let $\CU$ be a flat and quasi-projective integral scheme over $k$ of absolute dimension $d$. 
There is an intersection pairing
$$\wh\Pic (\CU/k)_\intb^{d}\lra \BR, \qquad  (\CHH_1, \cdots, \CHH_{d})\longmapsto\CHH_1 \cdots \CHH_{d},$$
defined as limits of the arithmetic (or geometric) intersection pairings of the projective case.

\kkk
Let $f:\CU\to \CV$ be a projective flat morphism of relative dimension $n$ of 
flat and quasi-projective integral schemes over $k$. 
Assume that $\CV$ is normal.
There is a relative intersection pairing
$$\wh \Picc (\CU/k)_\intb^{n+1}\lra \wh\Picc (\CV/k)_\intb, \qquad 
(\CLL_1, \cdots, \CLL_{n+1})\longmapsto f_*\pair{\CLL_1, \cdots, \CLL_{n+1}}.
$$
This is defined as the limit of the Deligne pairing.
The pairing of nef adelic line bundles is still nef. 
For simplicity, we may abbreviate $f_*\pair{\CLL_1, \cdots, \CLL_{n+1}}$
as $\pair{\CLL_1, \cdots, \CLL_{n+1}}$ if no confusion can occur.

\subsubsection*{Volume and bigness}

Here we recall the notions of volume and bigness of adelic line bundles in \cite[\S5.1-5.2]{YZ2}. 

\kkk
Let $\CU$ be a flat and quasi-projective integral scheme over $k$.
Let $\CLL$ be an adelic line bundle on $\CU$ with underlying line bundle $\CL$ on $\CU$.  
Define 
$$\wh H^0(\CU, \CLL):=\{s\in H^0(\CU, \CL): \|s(x)\|\leq 1,\forall x\in \CU^\an\}.$$
Here the metric $\|s(x)\|$ on the Berkovich space $\CU^\an$ is defined via the analytification 
functor
$\wh\Picc(\CU)\to \wh\Picc(\CU^\an)$. 
Elements of $\wh H^0(\CU, \CLL)$ are called \emph{effective sections} of $\CLL$ on $\CU$. 

If $k=\ZZ$, then $\wh H^0(\CU, \CLL)$ is a finite set, and we denote 
$$\wh h^0(\CU, \CLL):=\log\#\wh H^0(\CU, \CLL);$$ 
if $k$ is a field, then $\wh H^0(\CU, \CLL)$ is a finite-dimensional vector space over $k$, and we denote 
$$\wh h^0(\CU, \CLL):=\dim_k \wh H^0(\CU, \CLL).$$

Define the \emph{volume}
$$
\wh\vol(\CU,\CLL)
:=\lim_{m\to \infty} \frac{d!}{m^{d}}\wh h^0(\CU, m\CLL).
$$
Here $d$ is the absolute dimension of $\CU$. 
By \cite[Thm. 5.2.1]{YZ2}, the limit $\wh\vol(\CU,\CLL)$ always exists. 
The adelic line bundle $\CLL$ is said to be \emph{big} on $\CU$ if $\wh\vol(\CU,\CLL)>0$.

If $\CLL$ is \emph{nef}, the adelic Hilbert--Samuel formula in \cite[Thm. 5.2.2]{YZ2} asserts that
$$
\wh\vol(\CU,\CLL)
=\CLL^d.
$$
In this case, $\CLL$ is big if and only if $\CLL^d>0$.

\subsubsection*{Varieties over global fields}

By further direct limits, the above definitions and notations are extended to 
flat and essentially quasi-projective integral schemes over $k$ in \cite{YZ2}. 
The notion of ``essentially quasi-projective scheme'' is introduced in \cite[\S2.3]{YZ2}, which can be realized as an intersection of (possibly infinitely many) open subschemes of a quasi-projective scheme.
We do not need this generality in this paper, but we only need the following case of quasi-projective varieties over a number field or a function field of one variable. 

\kkk
If $k$ is a field, let $K$ be the function field of a projective regular curve $B$ over $k$. 
If $k=\ZZ$, let $K$ be a number field and denote $B=\Spec O_K$.
 
Let $X$ be a quasi-projective variety over $K$, which is viewed as a scheme over $k$. 
By a \emph{quasi-projective model} $\CU$ of $X$ over $k$, we mean a flat and quasi-projective integral scheme $\CU$ over $B$, together with an open immersion $X\to \CU_K$
over $K$.
Define 
$$
\wh\Div(X/k)=  \varinjlim_{\CU} \wh\Div(\CU/k),$$
$$
\wh\Pic(X/k)=  \varinjlim_{\CU} \wh\Pic(\CU/k),
$$ 
$$
\wh\Picc(X/k)=  \varinjlim_{\CU} \wh\Picc(\CU/k).$$
Here the limits are over all quasi-projective models $\CU$ of $X$ over $k$.
For the last direct limit, an object of $\wh\Picc(X/k)$ is a pair $(\CLL, \CU)$, where $\CU$ is a quasi-projective model of $X$ over $k$ and $\CLL$ is an object of $\wh\Picc(\CU/k)$. 
A morphism $(\CLL, \CU)\to (\CLL', \CU')$ between two objects of $\wh\Picc(X/k)$ is an isomorphism $\iota: \CL|_X\to \CL'|_X$ in $\Picc(X)$ satisfying the property that for some quasi-projective model $\CV$ of $X$ over $k$ endowed with open immersions $\psi:\CV\to \CU$ and $\psi':\CV\to \CU'$ extending the identity morphism $X\to X$, the isomorphism $\iota: \CL|_X\to \CL'|_X$ can be extended to an isomorphism $\CL|_\CV\to \CL'|_\CV$ in $\Picc (\CV)$ and induces
an isomorphism $\CLL|_\CV\to \CLL'|_\CV$ in $\wh \Picc (\CV/k)$. 
Here we take the convention $\CL|_X=(\CL|_\CU)|_X$, and
if $\CLL=(\CL, (\CX_i,\overline \CL_i, \ell_{i})_{i\geq 1})$ in $\wh\Picc (\CU/k)$, then $\CLL|_\CV=(\CL|_\CV, (\CX_i,\overline \CL_i, \ell_{i}|_\CV)_{i\geq 1})$ in 
$\wh\Picc (\CV/k)$.

Again, we often abbreviate the terms as 
$$\wh\Div(X), \quad
\wh\Picc(X), \quad
\wh\Pic(X).$$ 

If $X$ is projective over $K$, then $\wh\Picc(X)$ is essentially the category of adelic line bundles on $X$ introduced in Zhang \cite{Zha2}. 

The notions of effectivity, nefness and integrability of adelic line bundles (or adelic divisors) are also transferred to quasi-projective varieties over $K$ by taking direct limits. 
The intersection pairings are also valid in the current situation by taking direct limits.

The Berkovich analytic space $X^\an=(X/k)^\an$ is defined similar to $\CU^\an=(\CU/k)^\an$, which does not use the property that $\CU$ is quasi-projective over $k$.
The notions of metrized line bundles and arithmetic divisors on $X^\an$ are defined similarly. 
Then the direct limit process induces
injective analytification maps
$$
\wh\Div(X)\lra \wh\Div(X^\an),
$$
$$
\wh\Pic(X)\lra \wh\Pic(X^\an),
$$
and a fully faithful analytification functor
$$
\wh\Picc(X)\lra \wh\Picc(X^\an).
$$

Define the set of \emph{effective sections} of $\overline L$ on $X$ by
$$\wh H^0(X, \overline L):=\{s\in H^0(X, L): \|s(x)\|\leq 1,\forall x\in X^\an\}.$$
If $k=\ZZ$, denote 
$$\wh h^0(X, \OL):=\log\#\wh H^0(X, \OL);$$ 
if $k$ is a field, denote 
$$\wh h^0(X, \OL):=\dim_k \wh H^0(X, \OL).$$
Then $\wh h^0(X, \OL)$ is a well-defined (finite) number.
Define the \emph{volume}
$$
\wh\vol(X,\OL)
:=\lim_{m\to \infty} \frac{d!}{m^{d}}\wh h^0(X, m\OL).
$$
Here $d=\dim X+1$. 
Then the limit $\wh\vol(X,\OL)$ always exists. 
The adelic line bundle $\OL$ is said to be \emph{big} on $X$ if $\wh\vol(X,\OL)>0$.

If $\OL$ is \emph{nef}, 
the adelic Hilbert--Samuel formula asserts that
$$
\wh\vol(X,\OL)
=\OL^d.
$$
In this case,  $\OL$ is big if and only if $\OL^d>0$.


\subsubsection*{Some basic results}

We will need the following basic properties of adelic line bundles, which are not in \cite{YZ2}.

\begin{lem} \label{basic3}
\kkk
Let $X$ and $S$ be quasi-projective and flat integral schemes over $k$. 
Let $\pi:X\to S$ be a projective and flat morphism over $k$ with geometrically connected fibers. Then the following are true:

\begin{enumerate}[(1)]
\item The canonical map $\pi^*:\wh\Pic(S)\to \wh\Pic(X)$ is injective.
\item Assume that $S$ is normal or $\pi:X\to S$ admits a section. 
Let $\OL$ be an integrable adelic line bundle on $X$ such that the underlying line bundle $L=\pi^*M$ for some line bundle $M$ on $S$. 
Assume that for each $v\in S^\an$, under an isomorphism 
$L|_{X_{H_v}}\simeq \CO_{X_{H_v}}$ induced by $L=\pi^*M$,
the metric of $L|_{X_{H_v}}$ on $X_v^\an$ induced by $\OL$ corresponds to a constant multiple of the trivial metric of $\CO_{X_{H_v}}$ on $X_v^\an$.
Then there is an integrable adelic line bundle $\OM$ on $S$ with underlying line bundle $M$ such that $\OL$ is isomorphic to $\pi^*\OM$. 
\end{enumerate}

\end{lem}

\begin{proof}
By \cite[Prop. 3.4.1]{YZ2}, there is an injective analytification map
$$
\wh\Pic(X)\lra \wh\Pic(X^\an).
$$
By \cite[\S8.1, Prop. 4]{BLR}, $\pi^*:\Pic(S)\to \Pic(X)$ is injective.
Then (1) follows from the fact that the maps
$C(S^\an)\to C(X^\an)$ between the space of continuous functions is injective.

For (2), if $\pi:X\to S$ admits a section $x:S\to X$, then we can set $\OM=x^*\OL$, and check that it satisfies the requirement. 

Assume that $S$ is normal in the following. It suffices to prove the existence of an integrable adelic $\QQ$-line bundle $\OM$ on $S$ extending $M$ such that $\OL$ is isomorphic to $\pi^*\OM$ as adelic $\QQ$-line bundles. In fact, once the adelic $\QQ$-line bundle $\OM$ exists, it is automatically an adelic line bundle, since $M$ is a line bundle
(instead of only a $\QQ$-line bundle).

Let $\ON$ be an integrable adelic line bundle on $X$. 
Denote by $d$ the degree of the underlying line bundle $N$
of $\ON$ on the generic fiber of $X\to S$, and assume that $d>0$.
Set 
$$
\OM=d^{-1} \pair{\OL, \ON, \cdots, \ON} \in \wh\Picc(S)_\QQ.
$$ 
The motivation of this definition comes from \cite[Lem. 4.6.1(2)]{YZ2}. 

To check that $\pi^*\OM$ is isomorphic to $\OL$,   
we first note that their underlying line bundles are canonically isomorphic by 
\cite[Prop. 5.2.1.a]{MG}. 
By the analytification map, it suffices to check that the metrics of $L=\pi^*M$ 
over $X^\an$ induced by $\pi^*\OM$ and $\OL$ are equal.
For any $v\in S^\an$, the metric of $\pi^*M$ is equal to that of $L$ on the fiber $X_v^\an$ by the integration formula in \cite[Thm. 4.6.2]{YZ2} and  \cite[\S4.2.2]{YZ2}.
This finishes the proof.
\end{proof}

\begin{lem} \label{basic8}
\kkk
Let $f:X\to Y$ be a projective and flat morphism of relative dimension $n$ over $k$.
Here $X, Y$ are quasi-projective and flat integral schemes over $k$. 
Assume that $Y$ is normal of dimension $d\geq1$. 
Let $\OL_1,\cdots, \OL_{n+1}$ be integrable adelic line bundles on $X$.
\begin{enumerate}[(1)]
\item
If $\OL_1,\cdots, \OL_{n+1}$ are nef and big, then 
$f_*\pair{\OL_1,\cdots, \OL_{n+1}}$ is nef and big on $Y$.
\item
If $\OL_1$ is effective and $\OL_2,\cdots, \OL_{n+1}$ are nef, then  
$f_*\pair{\OL_1,\cdots, \OL_{n+1}}$ is pseudo-effective on $Y$ in the sense that 
$$f_*\pair{\OL_1,\cdots, \OL_{n+1}}\cdot \OM_1\cdots \OM_{d-1}\geq 0$$
for any nef adelic line bundles $\OM_1,\cdots, \OM_{d-1}$ on $Y$. 
\end{enumerate}
\end{lem}

\begin{proof}

Part (2) is a consequence of the identity 
$$f_*\pair{\OL_1,\cdots, \OL_{n+1}}\cdot \OM_1\cdots \OM_{d-1}
=\OL_1\cdots \OL_{n+1}\cdot f^*\OM_1\cdots f^*\OM_{d-1}.
$$
If $k=\ZZ$, the identity is a consequence of \cite[Lem. 4.6.1(1)]{YZ2} applying to the composition $X\to Y\to S$ with $S=\Spec\ZZ$. 
If $k$ is a field, by blowing-up $Y$, we can obtain a fibration $Y\to S$ for $S=\BP_k^1$. Then the identity still follows. 

For part (1), note that $f_*\pair{\OL_1,\cdots, \OL_{n+1}}$ is nef by 
\cite[Thm. 4.1.3]{YZ2}.
Let $\OM$ be a nef and big adelic line bundle on $Y$, which can be obtained as a model adelic line bundle induced by an ample hermitian line bundle on a projective model of $Y$ over $k$.
By \cite[Thm. 5.2.2(2)]{YZ2},
$$
\wh\vol(\OL_1-\epsilon f^*\OM)\geq \OL_1^{d+n}-(d+n)\OL_1^{d+n-1}\cdot \epsilon f^*\OM>0
$$ 
for some rational number $\epsilon>0$.
Then some positive multiple of $\ON=\OL_1-\epsilon f^*\OM$ is effective.
It follows that 
$$
 f_*\pair{\OL_1,\cdots, \OL_{n+1}}
= f_*\pair{\ON,\OL_2,\cdots, \OL_{n+1}}
+\epsilon\, f_*\pair{ f^*\OM,\OL_2,\cdots, \OL_{n+1}}
$$
By (2), $\OM_1= f_*\pair{\ON,\OL_2,\cdots, \OL_{n+1}}$ is pseudo-effective.
By \cite[Lem. 4.6.1(2)]{YZ2}, $\OM_2=\epsilon f_*\pair{ f^*\OM,\OL_2,\cdots, \OL_{n+1}}$ is a positive multiple of $\OM$, and thus it is big and nef. 

For $i=0, \cdots, d-1$, we have 
$$\OM_2^{i} \cdot (\OM_1+\OM_2)^{d-i}-\OM_2^{i+1} \cdot (\OM_1+\OM_2)^{d-i-1}
= \OM_2^{i} \cdot (\OM_1+\OM_2)^{d-i-1}\cdot \OM_1\geq0, $$
 since $\OM_2$ and $\OM_1+\OM_2$ are nef, and $\OM_1$ is pseudo-effective.
It follows that 
$$
(\OM_1+\OM_2)^d
\geq \OM_2\cdot(\OM_1+\OM_2)^{d-1}
\geq \OM_2^2(\OM_1+\OM_2)^{d-2}
\geq \cdots
\geq \OM_2^d
>0.$$
This finishes the proof.
\end{proof}

\subsection{Admissible adelic line bundles} \label{sec admissible}

We refer to \S\ref{sec appendix} for a review of admissible metrics of Arakelov \cite{Ara} and Zhang \cite{Zha1}.
The goal of this subsection is to extend the theory to quasi-projective families of curves. Our main result here is the following family version of Theorem \ref{admissible1}. 

\begin{thm} \label{admissible2}
Let $k$ be either $\ZZ$ or a field. 
Let $S$ be a quasi-projective and flat normal integral scheme over $k$. 
Let $\pi:X\to S$ be a smooth relative curve of genus $g>0$.
Denote by $\Delta:X\to X\times_S X$ the diagonal morphism. Then the following are true:
\begin{enumerate}[(a)]
\item
There is an adelic line bundle $\overline\omega_{X/S,a}$ in $\wh\Picc(X/k)$ with underlying line bundle $\omega_{X/S}$, 
such that for any $v\in S^\an$, the metric of $\omega_{X_{H_v}/H_v}$ on $X_{H_v}^\an$ induced by 
$\overline\omega_{X/S,a}$
is equal to the canonical admissible metric $\|\cdot\|_a$. 
Moreover, $\overline\omega_{X/S,a}$ is nef and unique up to isomorphism. 

\item
There is an adelic line bundle $\overline\CO(\Delta)_a$ in $\wh\Picc(X\times_S X/k)$ with underlying line bundle $\CO(\Delta)$, 
such that for any $v\in S^\an$, the metric of $\CO(\Delta_{H_v})$ on $(X_{H_v}^2)^\an$ induced by $\overline\CO(\Delta)_a$
is equal to the canonical admissible metric $\|\cdot\|_{\Delta,a}$. 
Moreover, $\overline\CO(\Delta)_a$ is integrable and unique up to isomorphism. 
\end{enumerate}
\medskip\noindent 
Moreover, the adelic line bundles satisfy the following extra properties:
\begin{enumerate}[(1)]
\item The canonical isomorphism
$$
\omega_{X/S} \lra \Delta^*\CO(-\Delta) 
$$
induces an isomorphism 
$$
\ol\omega_{X/S,a} \lra \Delta^*\ol\CO(-\Delta)_a.
$$
\item The canonical isomorphisms
$$
p_{1*}\pair{\CO(\Delta), p_2^* \omega_{X/S}} \lra \Delta^*p_2^* \omega_{X/S}\lra \omega_{X/S}
$$
induce an isomorphism
$$
p_{1*}\pair{\ol\CO(\Delta)_a, p_2^* \ol\omega_{X/S,a}} \lra \ol\omega_{X/S,a}.
$$
Here $p_1,p_2: X\times_S X\to X$ denote the two projections.
\end{enumerate}
\end{thm}

A few terms in the above statement require explanation. 
By abuse of notations, we write the diagonal divisor $\Delta(X)$ as $\Delta$. 
It is actually a Cartier divisor in $X\times_S X$ (cf. \cite[\S8.2, Lem. 6]{BLR}), so $\CO(\Delta)$ is a line bundle on $X\times_S X$. 

Recall that \S\ref{sec adelic} contains a brief review of the Berkovich space $S^\an$.
For any $v\in S^\an$, $H_v$ denotes the completed residue field of
$v$. We refer to \cite[\S3.1]{YZ2} for more details. 

The canonical admissible metrics are introduced for all $v\in S^\an$ 
in \S\ref{sec all val}. 
They are the Arakelov metrics in the archimedean case, the Zhang metrics in the non-archimedean case, and the canonical metric in the trivially valued case.

Finally, the uniqueness in (a) and (b) is a consequence of \cite[Prop. 3.4.1]{YZ2}, which gives a canonical fully faithful analytification functor
$$
\wh\Picc(X) \lra \wh\Picc(X^\an). 
$$

The proof of the existence is similar to that of Theorem \ref{admissible1}.
We first prove weak versions of (a) and (b), and then modify them to satisfy (1) and (2).
These weak versions are the admissible extensions generalizing the notions in \S\ref{sec appendix admissible}.

\subsubsection{Basics on Jacobian schemes} \label{sec Jacobian}

Here we present some basic constructions of abelian schemes and Jacobian schemes. For more details, we refer to \cite[Chap. 6]{MFK} and \cite[Chap. 9]{BLR}.
These constructions will be used in the proof of Theorem \ref{admissible1} and elsewhere of this paper.

Let $S$ be a noetherian scheme.
Let $\pi:X\to S$ be a smooth relative curve of genus $g>0$.
Denote by $J=\Pic^0_{X/S}$ the Jacobian scheme over $X$ over $S$. 
Let $\alpha$ be a line bundle on $X$ of degree $d$ on the fibers of $X\to S$. 
We have a finite $S$-morphism 
$$
i_\alpha: X\lra J, \quad x\longmapsto dx-\alpha.
$$
The right-hand side is understood in terms of the functor of points. 
We will also need the morphism
$$
(i_\alpha,i_\alpha): X\times_SX\lra J\times_SJ.
$$

Now we consider line bundles on $J$.
Without taking a base change of $S$, we have to deal with the case that there is no line bundle on $X$ of degree 1 on fibers over $S$.
Then it is not possible to define a theta divisor as in the case of algebraically closed fields.

However, there is still a canonical principal polarization $J\to J^\vee$ over $S$ in this generality, as proved in \cite[\S6.1, Prop. 6.9]{MFK}. 
Then we have a \emph{Poincar\'e line bundle} $P$ on $J\times_S J$, obtained by the {Poincar\'e line bundle} of $J\times_S J^\vee$ via the polarization. 
As a convention, the Poincar\'e line bundle $P$ is rigidified along zero sections in the sense that there are fixed isomorphisms $(e,\id)^*P\to \CO_J$ and $(\id,e)^*P\to \CO_J$. 
Here $e:S\to J$ is the identity section, and $(e,\id)$ and $(\id,e)$ are the morphisms $J\to J\times_SJ$ by the identity section. 
This determines the class of $P$ in $\Pic(J\times_S J)$ uniquely. 

\begin{definition} \label{qqq}
Define a line bundle on $J$ by
$$
\Theta=\Delta_J^*(P^\vee),
$$
where $\Delta_J:J\to J\times_S J$ is the diagonal morphism. 
\end{definition}

This is the construction outlined right before \cite[\S6.2, Prop. 6.10]{MFK}. 
By construction, $\Theta$ is symmetric and rigidified along the identity section $e:S\to J$, since the similar property holds for $P$.
By Theorem \ref{isomorphism}(4), $\Theta$ is algebraically equivalent to twice of a theta divisor on each geometric fiber of $J\to S$. 
Therefore, $\Theta$ is relatively ample over $S$.

\subsubsection{Admissible extensions}

Recall that in \S\ref{sec appendix admissible} (together with \S\ref{sec ara}, \S\ref{sec all val}), we have a notion of admissible metrics of line bundles on abelian varieties, smooth curves, square of smooth curves over a complete valuation field. 
Now we extend this to a family version as follows. 

Resume the notation of Theorem \ref{admissible2}. 
Namely, let $k$ be either $\ZZ$ or a field. 
Let $S$ be a quasi-projective and flat normal integral scheme over $k$. 
Let $\pi:X\to S$ be a smooth relative curve of genus $g>0$.
We introduce the following definition. 

\begin{definition} \label{def admissible family}
Let $W$ be $X$, $X\times_SX$, or an abelian scheme over $S$. 
Let $\OL$ be an integrable adelic line bundle on $W/k$ with underlying line bundle $L$.
We say that $\OL$ is \emph{admissible} if for any $v\in S^\an$, the metric of $L|W_{H_v}$ on $W_{H_v}^\an$ induced by $\OL$
is admissible. 
In this case, we say that $\OL$ is an \emph{admissible adelic extension} of $L$. 
\end{definition}

\begin{prop} \label{prop admissible family}
Let $W$ be $X$, $X\times_SX$, or an abelian scheme over $S$. 
Then any line bundle $L$ on $W$ admits an admissible adelic extension in 
$\wh\Pic(W/k)_\intb$, unique up to translation by $\ker(\wh\Pic(S/k)_\intb\to \Pic(S))$. 
\end{prop}

\begin{proof}
The uniqueness follows from Lemma \ref{basic3}(2).
The existence is similar to the exposition in \S\ref{sec appendix admissible}.
For example, the case of abelian schemes is given by the dynamically invariant extension constructed in \cite[Thm. 6.1.1]{YZ2}.

To treat the other cases, let us first introduce some notations.
For any projective and flat $S$-schemes $Y$, denote by $\Pic^0(Y)$ the subgroup of $\Pic(Y)$ of line bundles on $Y$ which are algebraically equivalent to 0 on fibers of $Y\to S$. 
Denote by $\Pic_{Y/S}(S)$ the relative Picard group, and denote by $\Pic_{Y/S}^0(S)$ the subgroup of line bundles which are algebraically equivalent to 0 on fibers of $Y\to S$.

For any projective and flat $S$-schemes $Y_1,Y_2$, denote by $\Pic^{00}(Y_1\times_SY_2)$ the subgroup of $\Pic(Y_1\times_SY_2)$
of line bundles which are algebraically equivalent to 0 on fibers of both projections
$p_i:Y_1\times_SY_2\to Y_i$.

For the case $W=X$, let $\alpha$ be a line bundle on $X$ of degree $d> 0$ on the fibers of $X\to S$. It defines a finite morphism $i_\alpha: X\to J.$
By pull-back of an admissible adelic extension of the line bundle $\Theta$ from $J$ to $X$, 
it suffices to prove that the map $i_\alpha^*: \Pic^0(J) \to \Pic^0(X)$ has a torsion cokernel, or equivalently that the map 
$i_\alpha^*: \Pic^0(J)/\Pic(S) \to \Pic^0(X)/\Pic(S)$ has a torsion cokernel.
As a consequence of the Leray spectral sequence (cf. \cite[\S8.1, Prop. 4]{BLR}), there is a canonical isomorphism  
$\Pic^0(J)/\Pic(S)\to \Pic^0_{J/S}(S)$
and a canonical injection $\Pic^0(X)/\Pic(S)\to \Pic^0_{X/S}(S)$. 
Then its suffices to prove that the map 
$i_\alpha^*: \Pic^0_{J/S}(S)\to \Pic^0_{X/S}(S)$ has a torsion cokernel.

Note that we have a canonical isomorphism  
$\Pic^0_{J/S}\to \Pic^0_{X/S}$, which can be established via flat descent. Under the isomorphism, 
$i_\alpha^*: \Pic^0_{J/S}\to \Pic^0_{X/S}$ is isomorphic to the multiplication 
$[d]:\Pic^0_{X/S}\to \Pic^0_{X/S}$, which can be verified fiberwise.
This proves the case $W=X$.  

For the case $W=X\times_SX$, let $\alpha$ be as above, and consider the morphism
$(i_\alpha,i_\alpha): X\times_SX\to J\times_SJ.$ 
By pull-back of admissible extensions of $p_i^*\Theta$ from $J\times_SJ$ to $X\times_SX$, it suffices to prove that the canonical map 
$$i_\alpha^*: \Pic^{00}(J\times_SJ) \lra \Pic^{00}(X\times_SX)$$ 
has a torsion cokernel.

Apply the above result to the relative curve $p_1:X\times_SX\to X$, whose relative Jacobian is just $q_1:X\times_SJ\to X$. We see that the canonical map
$$
\Pic^{00}(X\times_SJ)\lra \Pic^{00}(X\times_SX)
$$
has a torsion cokernel. Similarly, apply the above result to the relative curve $q_2:X\times_SJ\to J$, whose relative Jacobian is just $p_2:J\times_SJ\to J$. We see that the canonical map
$$
\Pic^{00}(J\times_SJ)\lra \Pic^{00}(X\times_SJ)
$$
has a torsion cokernel. 
The result is proved by combining these two results. 
\end{proof}

\subsubsection{Construction of canonical admissible extensions}

Now we prove Theorem \ref{admissible2}. Let $X/S/k$ be as in the theorem.

By Proposition \ref{prop admissible family}, there is an admissible adelic extension $\overline\omega_{X/S,a}'$ of $\omega_{X/S}$ in $\wh\Picc(X/k)$
and an admissible adelic extension $\overline\CO(\Delta)_a'$ of $\CO(\Delta)$ in $\wh\Picc(X\times_S X/k)$. 
We will modify $(\overline\omega_{X/S,a}',\overline\CO(\Delta)_a')$ to meet Theorem \ref{admissible2}. 
Namely, we will prove that there are integrable adelic line bundles $\ol\beta_1$ and $\ol\beta_2$ in 
$\wh\Picc(S)$ with underlying line bundles $\beta_1$ and $\beta_2$ in $\Picc(S)$,  together with isomorphisms $\beta_1 \to \CO_S$ and $\beta_2 \to \CO_S$, such that the adelic line bundles
$$
\overline\omega_{X/S,a}: =\overline\omega_{X/S,a}'+ \pi^* \ol\beta_1,\qquad
\overline\CO(\Delta)_a:= \overline\CO(\Delta)_a' + (\pi,\pi)^*\ol\beta_2
$$
satisfy (1) and (2) of Theorem \ref{admissible2}.
Here $(\pi,\pi):X\times_S X\to S$ is the structure morphism. 

We first look at (2), which requires an isomorphism
$$
p_{1*}\pair{\ol\CO(\Delta)_a, p_2^* \ol\omega_{X/S,a}} \lra \ol\omega_{X/S,a}
$$
extending the canonical isomorphism of the corresponding underlying line bundles.
By \cite[Lem. 4.6.1(2)]{YZ2}, there is a canonical isomorphism
$$
p_{1*}\pair{\ol\CO(\Delta)_a, (\pi,\pi)^* \ol\beta_1} \lra \pi^* \ol\beta_1.
$$
Thus (2) is equivalent to an isomorphism
$$
p_{1*}\pair{\ol\CO(\Delta)_a, p_2^* \ol\omega_{X/S,a}'} \lra \ol\omega_{X/S,a}'.
$$
By \cite[Lem. 4.6.1(2)]{YZ2} again, it suffices to have a canonical isomorphism
$$
p_{1*}\pair{\ol\CO(\Delta)_a', \ p_2^* \ol\omega_{X/S,a}'}+(2g-2)\ol\beta_2 \lra \ol\omega_{X/S,a}'.
$$
Define $\ol\beta_2\in \wh\Picc(S)$ with trivial underlying line bundle satisfying the above isomorphism. 

This establishes (2), which determines $\ol \beta_2$ and thus 
$\overline\CO(\Delta)_a$. Similarly, we use (1) to determine $\ol \beta_1$ and thus
$\overline\omega_{X/S,a}$.
This proves Theorem \ref{admissible2}.

\subsection{Intersection and pull-back} \label{sec pullback}

The goal of this subsection is to consider some properties of admissible 
adelic line bundles.
We first present the adjunction formula and the arithmetic Hodge index theorem, and then prove some explicit formulae about pull-back of admissible adelic line bundles 
from the Jacobian schemes. 

For convenience, we introduce the notations for this subsection. 
Let $k$ be either $\ZZ$ or a field. 
Let $S$ be a quasi-projective and flat normal integral scheme over $k$. 
Let $\pi:X\to S$ be a smooth relative curve of genus $g>0$.
Denote by $\Delta:X\to X\times_S X$ the diagonal morphism.
Let $\ol\omega_{X/S,a}$ and $\ol\CO(\Delta)_a$ be the canonical admissible extensions in Theorem \ref{admissible2}. 

For any section $x:S\to X$ of $\pi:X\to S$, 
note that $x(S)$ is a Cartier divisor on $X$ (cf. \cite[\S8.2, Lem. 6]{BLR}). 
Denote by $\CO(x)$ the line bundle on $X$ associated to the Cartier divisor $x(S)$. 
Denote by $\ol\CO(x)_a$ the pull-back of $\ol\CO(\Delta)_a$ via the morphism $(x,\id):X\to X\times_SX$, which is an admissible adelic extension of $\CO(x)$ on $X/k$.
We call $\ol\CO(x)_a$ the \emph{canonical admissible adelic extension} of $\CO(x)$ on $X/k$.

By Definition \ref{qqq}, we have a symmetric, relatively ample and rigidified line bundle $\Theta$ on $J$.
By \cite[Thm. 6.1.1(3)]{YZ2}, there is a nef adelic line bundle $\OTheta$ on $J/k$ extending $\Theta$ such that $[2]^*\OTheta=4\OTheta$ in $\wh\Pic(J/k)$.

\subsubsection{Adjunction formula and Hodge index theorem}

We first have the following adjunction formula, which generalizes \cite[Thm. 4.2]{Zha1}.

\begin{prop} \label{adjunction}
Let $x:S\to X$ be a section of $\pi:X\to S$.
\begin{enumerate}[(1)]
\item There is a canonical isomorphism
$$\pi_*\pair{\ol\CO(x)_a, \OL}
\lra x^*\OL$$
for any admissible adelic line bundle $\OL$ on $X/k$. 
\item There is a canonical isomorphism
$$\pi_*\pair{\ol\CO(x)_a, \ol\CO(x)_a}\lra -\pi_*\pair{\ol\CO(x)_a, \ol\omega_{X/S,a}}.$$

\end{enumerate}
\end{prop}
\begin{proof}

We first prove (1). 
There is a canonical isomorphism $\pi_*\pair{\CO(x), L} \to x^*L$. 
By the analytification functor $\wh\Picc(S)\to \wh\Picc(S^\an)$, it suffices to prove that the norm of this isomorphism is $1$ at any $v\in S^\an$. 
By the integration formulae in \cite[Thm. 4.6.2]{YZ2} and \cite[\S4.2.2]{YZ2}, the logarithm of this norm is given by 
$$
\int_{X_v^\an} \log\|1\|_{v} c_1(\OL)_v.
$$
Here $1$ is the canonical section of $\CO(x)$. 
This vanishes by Theorem \ref{admissible1}.

For (2), recall that Theorem \ref{admissible2}(1) gives an isomorphism 
$\ol\omega_{X/S,a} \to \Delta^*\ol\CO(-\Delta)_a.$
Take the pull-back of this isomorphism via $x:S\to X$, we obtain an isomorphism 
$x^*\ol\omega_{X/S,a} \to x^*\ol\CO(-x)_a.$
Combining with (1), this gives (2).
\end{proof}

Now we have the following consequence of the arithmetic Hodge index theorem of  \cite{YZ2}.

\begin{thm} \label{hodge index}

Let $M$ be a line bundle on $X$ with degree 0 on fibers of $X\to S$.
Let $\OM$ be an admissible adelic extension of $M$ on $X$.
Then the following hold.
\begin{enumerate}[(1)]
\item Denote by $\wh\Pic(X/k)_{\vert}$ the kernel of the forgetful map 
$\wh\Pic(X/k)_{\intb}\to \Pic(X)$. Then for any $\overline V\in \wh\Pic(X/k)_{\vert}$,
$$\pi_*\pair{\overline M,\overline V}= 0$$ 
in $\wh\Pic(S/k)$.
\item
Denote by $\iota: S\to J$ the section of $J\to S$ corresponding to the line bundle $M$ on $X$. Then
$$\pi_*\pair{\overline M, \overline M}= -\iota^*\OTheta$$
in $\wh\Pic(S/k)_{\QQ}$.
\end{enumerate}
\end{thm}
\begin{proof}
This is a variant of the arithmetic Hodge index theorem of \cite[Thm. 6.5.4]{YZ2}. 
It suffices to prove (1). 
The proof is similar to that of Proposition \ref{adjunction}(1).
In fact, fix an isomorphism $V\simeq \CO_X$. 
The isomorphism induces a canonical isomorphism $\pi_*\pair{\OM,\ol V} \to \CO_S$. 
By the analytification functor $\wh\Picc(S)\to \wh\Picc(S^\an)$, it suffices to prove that the norm of this isomorphism is $1$ at any $v\in S^\an$. 
By definition, the logarithm of this norm is given by 
$$
\int_{X_v^\an} \log\|1\|_v c_1(\OM)_v.
$$
Here $1$ is viewed as a section of $V$ via $V\simeq \CO_X$.  
By Theorem \ref{admissible1}, $c_1(\OM)_v=0$ since $\deg(M_v)=0$. 
This finishes the proof.
\end{proof}

\subsubsection{Pull-back formula: general}

Resume the above notations. Namely, let $k$ be either $\ZZ$ or a field, let $S$ be a quasi-projective and flat normal integral scheme over $k$, and let $\pi:X\to S$ be a smooth relative curve of genus $g>0$.

We first present the following general pull-back formula, which generalizes its primitive version in Theorem \ref{isomorphism}.

\begin{thm} \label{isomorphism5}
Let $\alpha$ be a line bundle on $X$ of degree $d$ on the fibers of $\pi:X\to S$.
Let $\ol\alpha$ be an admissible adelic extension of  $\alpha$ on $X$. 
\begin{enumerate}[(1)]
\item 
There is an identity
$$
i_\alpha^*\OTheta= d^2\ol\omega_{X/S,a}+2d\ol\alpha-\pi^*\pi_*\pair{\ol\alpha,\ol\alpha}
$$
in $\wh\Pic(X)_\QQ$.

\item
There are identities
$$
\pi_*\pair{i_\alpha^*\OTheta,i_\alpha^*\OTheta}
= d^4\pi_*\pair{\ol\omega_{X/S,a}\, \ol\omega_{X/S,a}}
+4d^3\pi_*\pair{\ol\omega_{X/S,a}\, \ol\alpha}
-(4g-4)d^2\pi_*\pair{\ol\alpha,\ol\alpha}
$$
and
$$
(g-1)\pi_*\pair{i_\alpha^*\OTheta,i_\alpha^*\OTheta}
= gd^4\pi_*\pair{\ol\omega_{X/S,a}\, \ol\omega_{X/S,a}}
+d^2\iota_\alpha^*\OTheta
$$
in $\wh\Pic(S)_\QQ$.
Here $\iota_\alpha:S\to J$ is the section corresponding to the line bundle $(2g-2)\alpha-d\omega_{X/S}$ on $X$. 

\item Assume  $\alpha=\CO(x)$ for a section $x:S\to X$ of $\pi:X\to S$. 
Then 
$$
i_\alpha^*\OTheta= \ol\omega_{X/S,a}+2\ol\CO(x)_a+\pi^*x^*\ol\omega_{X/S,a}
$$
in $\wh\Pic(X)_\QQ$, and
$$
\pi_*\pair{i_\alpha^*\OTheta,i_\alpha^*\OTheta}
= \pi_*\pair{\ol\omega_{X/S,a}\, \ol\omega_{X/S,a}}
+4g\,x^*\ol\omega_{X/S,a}
$$
in $\wh\Pic(S)_\QQ$.
\end{enumerate}
\end{thm}
\begin{proof}

It is easy to see that (3) is a consequence of (1) and (2) for $\ol\alpha=\ol\CO(x)_a$ by the adjunction formula in Proposition \ref{adjunction}. 

To prove (2) by (1), by \cite[Lem. 4.6.1(2)]{YZ2}, 
$$
\pair{i_\alpha^*\OTheta,i_\alpha^*\OTheta}
= \pair{d^2\ol\omega_{X/S,a}+2d\ol\alpha,\, d^2\ol\omega_{X/S,a}+2d\ol\alpha}
-4gd^2\pair{\ol\alpha,\ol\alpha}.
$$
This gives
$$
\pair{i_\alpha^*\OTheta,i_\alpha^*\OTheta}
= d^4\pair{\ol\omega_{X/S,a}\, \ol\omega_{X/S,a}}
+4d^3\pair{\ol\omega_{X/S,a}\, \ol\alpha}
-(4g-4)d^2\pair{\ol\alpha,\ol\alpha}.
$$
and
$$
(g-1)\pair{i_\alpha^*\OTheta,i_\alpha^*\OTheta}
=  gd^4 \pair{\ol\omega_{X/S,a}\, \ol\omega_{X/S,a}}
-d^2 \pair{(2g-2)\ol\alpha-d\ol\omega_{X/S,a}, \, (2g-2)\ol\alpha-d\ol\omega_{X/S,a}}.
$$
Then it suffices to prove 
$$
\pair{(2g-2)\ol\alpha-d\ol\omega_{X/S,a}, \, (2g-2)\ol\alpha-d\ol\omega_{X/S,a}}=-\iota_\alpha^*\OTheta
$$
in $\wh\Pic(S)_\QQ$.
This follows from the arithmetic Hodge index theorem in Theorem \ref{hodge index}.

It remains to prove (1). 
Denote
$$
\ON=i_\alpha^*\OTheta- d^2\ol\omega_{X/S,a}-2d\ol\alpha+\pi^*\pi_*\pair{\ol\alpha,\ol\alpha}.
$$
We need to prove that $\ON=0$ in $\wh\Pic(X)_\QQ$. 

First, we prove that for any section $x:S\to X$ of $\pi:X\to S$, the pull-back $x^*\ON=0$ in 
$\wh\Pic(S)_\QQ$. 
In fact, $i_\alpha:X\to J$ maps $x$ to the section $\iota:S\to J$ corresponding to the line bundle $d\CO(x)-\alpha$ on $X$.
By the arithmetic Hodge index theorem in Theorem \ref{hodge index}, 
$$x^*i_\alpha^*\OTheta=\iota^*\OTheta
=-\pi_*\pair{d\ol\CO(x)_a-\ol\alpha,\, d\ol\CO(x)_a-\ol\alpha}.$$ 
By the adjunction formula in Proposition \ref{adjunction},
$$
x^*\ol\alpha=\pi_*\pair{\ol\CO(x)_a,\, \ol\alpha}
$$
and
$$
x^*\ol\omega_{X/S,a}=\pi_*\pair{\ol\CO(x)_a,\, \ol\omega_{X/S,a}}=-\pi_*\pair{\ol\CO(x)_a,\, \ol\CO(x)_a}.
$$ 
Then $x^*\ON=0$ by bi-linearity of the Deligne pairing. 

Second, we extend the result to the statement that, for any flat and quasi-projective normal integral scheme $S'$ over $S$, and for any $S$-morphism $y:S'\to X$, the pull-back $y^*\ON=0$ in 
$\wh\Pic(S')_\QQ$. 
In fact, denote by $\pi':X'\to S'$ the base change of $\pi:X\to S$ by $S'\to S$. 
Then $y:S'\to X$ induces a section $y':S'\to X'$ of $\pi':X'\to S'$.
Denote by $\ON'$ the pull-back of $\ON$ to $\wh\Pic(X')_\QQ$. 
It suffices to prove $y'^*\ON'=0$ in $\wh\Pic(S')_\QQ$.
This is reduced to the above case. 

Third, in above statement, take $S'=X$ and take $y:S'\to X$ to be the identity section. Then 
$y^*\ON=0$ in $\wh\Pic(S')_\QQ$ just means $\ON=0$ in 
$\wh\Pic(X)_\QQ$.  
This finishes the proof. 
\end{proof}

\subsubsection{Pull-back formula: explicit examples}

As above, let $k$ be either $\ZZ$ or a field, let $S$ be a quasi-projective and flat normal integral scheme over $k$, and let $\pi:X\to S$ be a smooth relative curve of genus $g$.
The following special examples of the theorem will be very useful for our later considerations.

\begin{thm} \label{isomorphism3}
Assume that $g>1$. 
\begin{enumerate}[(1)]
\item Let $\omega=\omega_{X/S}$ be the canonical bundle of $X$ over $S$. Then 
$$i_\omega^*\OTheta= 4g(g-1)\overline\omega_{X/S,a}-\pi^*\pi_*\pair{\overline\omega_{X/S,a},\overline\omega_{X/S,a}},$$
in  $\wh\Pic(X)_\QQ$, and 
$$\pi_*\pair{i_\omega^*\OTheta,i_\omega^*\OTheta}= 16g(g-1)^3 \pi_*\pair{\overline\omega_{X/S,a},\overline\omega_{X/S,a}}$$
in $\wh\Pic(S)_\QQ$.
Moreover, $\overline\omega_{X/S,a}$ is nef on $X$, and 
$\pi_*\pair{\overline\omega_{X/S,a},\overline\omega_{X/S,a}}$ is nef on $S$. 

\item Let $p_i:X\times_S X\to X$ be the projection to the $i$-th factor for $i=1,2$, and let
$j$ be the morphism 
$$j:X\times_SX \lra J, \quad (x,y)\longmapsto y-x.$$
Then
$$
j^*\OTheta=2\ol\CO(\Delta)_a+p_1^*\overline\omega_{X/S,a}+p_2^*\overline\omega_{X/S,a}
$$
in $\wh\Pic(X\times_SX)_\QQ$,
and 
$$
p_{1*}\pair{j^*\OTheta,j^*\OTheta}=
4g\, \overline\omega_{X/S,a}+\pi^*\pi_*\pair{\overline\omega_{X/S,a},\overline\omega_{X/S,a}}
$$
in $\wh\Pic(X)_\QQ$.


\item 
Let $\tau$ be the morphism 
$$\tau:J\times_SX \lra J\times_SJ, \quad (y,x)\longmapsto (y,y+(2g-2)x-\omega_{X/S}).$$
Then 
$$
q_{1*}\pair{\tau^*(\OTheta_J), \tau^*(\OTheta_J)}= 
16(g-1)^3\OTheta+16g(g-1)^3\pi_J^*\pi_*\pair{\ol\omega_{X/S,a}\, ,\ol\omega_{X/S,a}}
$$
in $\wh\Pic(J)_\QQ$.
Here $\pi_J:J\to S$ denotes the structure morphism. 
\end{enumerate}
\end{thm}
\begin{proof}
Note that (1) is a direct consequence of Theorem \ref{isomorphism5}(1)(2).
The nefness of $\overline\omega_{X/S,a}$ comes from the nefness of $\OTheta$ by the formula
$$\overline\omega_{X/S,a}
=\frac{1}{4g(g-1)}i_\omega^*\OTheta+ \frac{1}{64g^2(g-1)^4}\pi^*\pi_*\pair{i_\omega^*\OTheta,i_\omega^*\OTheta}.$$

For (2), the key is to interpret $j$ in the form of Theorem \ref{isomorphism5}(3).
Denote $X_X=X\times_S X$ and $J_X=X\times_S J$, viewed as $X$-schemes via the first projections $p_1:X_X\to X$ and $p_1:J_X\to X$. Here we use $p_1$ twice by abuse of notations. 
Then $p_1:J_X\to X$ is canonically isomorphic to the Jacobian scheme of $p_1:X_X\to X$.
View $\Delta:X\to X\times_SX=X_X$ as a section of $p_1:X_X\to X$, which defines an $X$-morphism 
$$
i_\Delta=i_{\CO(\Delta)}: X_X\lra  J_X,\quad x\longmapsto x-\Delta.
$$
Denote by $\OTheta_X$ the pull-back of the adelic line bundle $\OTheta$ via the projection $J_X\to J$. 

Apply Theorem \ref{isomorphism5}(3) to the $X$-morphism $i_\Delta: X_X\to  J_X$. 
We obtain 
$$
i_\Delta^*(\OTheta_X)=2\ol\CO(\Delta)_a+p_1^*\overline\omega_{X/S,a}+p_2^*\overline\omega_{X/S,a}
$$
in $\wh\Pic(X_X)$,
and
$$
p_{1*}\pair{i_\Delta^*(\OTheta_X),i_\Delta^*(\OTheta_X)}=
4g\, \overline\omega_{X/S,a}+\pi^*\pi_*\pair{\overline\omega_{X/S,a},\overline\omega_{X/S,a}}
$$
in $\wh\Pic(X)$.

On the other hand, we can write $i_\Delta$ as
$$
i_\Delta: X\times_SX \lra  X\times_SJ,\quad (x,y)\longmapsto (x,y-x).
$$
Thus 
$j:X\times_SX\to J$ is equal to the composition of $i_\Delta$ with 
the projection $p_2:X\times_SJ\to J$. 
This gives $j^*(\OTheta)=i_\Delta^*(\OTheta_X)$.
Then (2) is proved.

To prove (3), we first consider the case that $\pi:X\to S$ has a section $x:S\to X$. 
This gives an immersion 
$$(\id,i_x): J\times_SX \lra J\times_SJ.$$
Then we have a universal line bundle $Q=(\id,i_x)^*P$ on $J\times_SX$, where $P$ is the Poincar\'e line bundle on $J\times_S J$. 
As in \cite[\S6.5.2]{YZ2} or \cite[Thm. 6.1.2]{YZ2}, there is a canonical extension $\OQ$ of $Q$ in $\wh\Pic(J\times_SX)_\intb$ such that $[2]_X^*\OQ=2\OQ$. 
Here $[2]_X:J\times_SX\to J\times_SX$ is the base change of 
$[2]:J\to J$.

Write $X_J=J\times_SX$ and $J_J=J\times_SJ$, viewed as $J$-schemes by projections to the first factors.
Then $J_J$ is canonically isomorphic to the Jacobian scheme of $X_J$ over $J$.
We claim that the morphism
$$\tau:J\times_SX \lra J\times_SJ, \quad (y,x)\longmapsto (y,y+(2g-2)x-\omega_{X/S})$$
is equal to 
$$
i_{\omega-Q}: X_J\lra J_J,\quad x\longmapsto (2g-2)x-(\omega_{X/S}-Q).
$$

To prove the claim, we can assume that $S$ is the spectrum of an algebraically closed field. Then for any point $\beta\in J(S)$, the fiber of $i_{\omega-Q}: X_J\to J_J$ above $\beta$ is isomorphic to the morphism 
$$
i_{\omega-Q_\beta}: X\lra J, \quad x\longmapsto (2g-2)x-(\omega_{X_J/J}-Q_\beta).
$$
By the universal property, $Q_\beta=Q|_{X\times_S\beta}=\beta$ as a line bundle on $X$. 
This proves the claim.

With the claim, apply Theorem \ref{isomorphism5}(2) to $i_{\omega-Q}:X_J\to J_J$.
It gives
$$
q_{1*}\pair{\tau^*(\OTheta_J),\, \tau^*(\OTheta_J)}
= 16g(g-1)^3\pi_J^*\pi_*\pair{\ol\omega_{X/S,a},\, \ol\omega_{X/S,a}}
+4(g-1)\iota^*\OTheta.
$$
Here $\iota:J\to J$ is the morphism corresponding to the line bundle 
$$(2g-2)(\omega_{X/S}-Q)-(2g-2)\omega_{X/S}=-(2g-2)Q$$ 
on $J\times_SX$, and thus it is the multiplication morphism $[-(2g-2)]:J\to J$.
This gives $\iota^*\OTheta=4(g-1)^2\OTheta$.  
This proves the result assuming that $X\to S$ has a section. 

As in the proof of Theorem \ref{isomorphism5}, it is easy to extend this result to the general case. 
In fact, take the base change of everything by $X\to S$, which converts to the case with sections, and then we can recover the original case by Lemma \ref{basic3}(1).
\end{proof}

\section{Bigness of the admissible canonical bundle} \label{sec 3}

The goal of this section is to study the positivity of the admissible canonical line bundle 
$\overline\omega_{X/S,a}$ in Theorem \ref{admissible2}.
The following is the main result of this section. 

\begin{thm}  \label{bigness5}
Let $k$ be either $\ZZ$ or a field. 
Let $S$ be a quasi-projective and flat normal integral scheme over $k$. 
Let $\pi:X\to S$ be a smooth relative curve over $S$ of genus $g>1$ with maximal variation.
Then the admissible canonical bundle $\overline\omega_{X/S,a}$ is nef and big on $X$. 
\end{thm}

Recall that the relative curve $\pi:X\to S$ has \emph{maximal variation} if the moduli morphism $S\to M_{g,k}$ is generically finite, where $M_{g,k}$ denotes the coarse moduli scheme of smooth curves of genus $g$ over $k$.
Note that we require $g>1$ in this section, while we only require $g>0$ in the previous section.
An equivalent form of the theorem is as follows.

\begin{thm} \label{bigness1}
Let $k$ be either $\ZZ$ or a field. 
Let $S$ be a quasi-projective and flat normal integral scheme over $k$. 
Let $\pi:X\to S$ be a smooth relative curve over $S$ of genus $g>1$ with maximal variation.
Then the Deligne pairing $\pi_*\pair{\overline\omega_{X/S,a},\overline\omega_{X/S,a}}$ is nef and big on $S$.
\end{thm}

As a convention, if $\dim S=0$ (and thus $k$ is a field), all (adelic) line bundles on $S/k$ are nef and big. 
By Lemma \ref{basic8}(1), the Deligne pairing of nef and big adelic line bundles is still nef and big. Then Theorem \ref{bigness1} is an easy consequence of Theorem \ref{bigness5}. 
However, in our treatment, we will first prove Theorem \ref{bigness1} and then use it to deduce Theorem \ref{bigness5} by the following argument. 

\begin{proof}[Proof of Theorem \ref{bigness5} by Theorem \ref{bigness1}]
By Theorem \ref{isomorphism3}(1), $\overline\omega_{X/S,a}$ is nef on $X$, and
$$4g(g-1)\overline\omega_{X/S,a}=i_\omega^*\OTheta+\pi^*\pi_*\pair{\overline\omega_{X/S,a},\overline\omega_{X/S,a}}.$$
It suffices to estimate the self-intersection of the right-hand side. 
Denote by $d=\dim S$. 
Then the self-intersection of the right-hand side contains a term
$$
i_\omega^*\OTheta\cdot \big(\pi^*\pi_*\pair{\overline\omega_{X/S,a},\overline\omega_{X/S,a}}\big)^d
=a \big(\pi_*\pair{\overline\omega_{X/S,a},\overline\omega_{X/S,a}}\big)^d
>0.
$$ 
Here $a=4g(g-1)(2g-2)$ is the degree of $i_\omega^*\OTheta$ on the generic fiber of $X\to S$. This finishes the proof. 
\end{proof}

The idea of the proof of Theorem \ref{bigness1} is explained in the introduction.
The first three subsections of this section are preparatory results for the proof, and then the last two subsections give the proofs of the geometric case and the arithmetic case respectively.

\subsection{Basics on the moduli space of curves}

In this subsection, we review some standard notions on moduli spaces of curves and stable curves.

\subsubsection{Basics on the moduli space} \label{sec basic moduli}

We refer to \cite[\S1]{DM} (or our \S\ref{sec notation}) for quick definitions of stable curves.

\kkk
Fix an integer $g>1$.
Denote by $\CM_{g}$ the moduli stack of smooth curves of genus $g$ over $k$, and by $\overline\CM_{g}$ the moduli stack of stable curves of genus $g$ over $k$. Denote by $\pi_g:\CCC_g\to \CM_{g}$ and $\pi_g:\ol\CCC_g\to \CMM_{g}$ respectively the universal curves of the stacks.

By \cite[Thm. 5.2]{DM}, $\CM_{g}$ is smooth over $k$, and $\overline\CM_{g}$ is proper and smooth over $k$. 
Denote by $M_g$ and $\overline M_{g}$ the coarse moduli schemes of $\CM_{g}$ and $\overline\CM_{g}$ respectively. 
We refer to \cite{DM} for more details on the moduli stacks.

By the standard theory of stacks, line bundles and Cartier divisors on the stacks are defined in terms of descent. We refer to \cite[\S1]{AC} and \cite[\S4]{CH} for explanation on line bundles on $\CM_{g}$ and $\ol\CM_{g}$. 
Or one can always take a fine level structure to avoid stacks.
Over $\ol\CM_{g}$, there are a Hodge bundle $\lambda$ and tautological divisors $\Delta_0,\cdots, \Delta_{[g/2]}$ defined as follows. 

The Hodge bundle on $\CMM_{g}$ is defined by
$$
\lambda_{\CMM_g}:=\det \pi_{g*} (\omega_{\ol\CCC_g/ \CMM_{g}})
=\wedge^g \pi_{g*} (\omega_{\ol\CCC_g/ \CMM_{g}}),
$$
where $\omega_{\ol\CCC_g/ \CMM_{g}}$ is the relative dualizing sheaf. 
It is a line bundle on $\CMM_{g}$. 
Note that in the literature, the term ``Hodge bundle'' may also mean the vector bundle  
$\pi_{g*} (\omega_{\ol\CCC_g/ \CMM_{g}})$, but in this paper it always means the determinant of $\pi_{g*} (\omega_{\ol\CCC_g/ \CMM_{g}})$.

By \cite[Thm. 5.2]{DM},
$$\Delta:=\overline\CM_{g}\setminus \CM_{g}$$
is a divisor of normal crossing in $\overline\CM_{g}$.
In terms of ``prime divisors'', we can write
$$\Delta=\Delta_0\cup \Delta_1\cup \cdots \cup \Delta_{[g/2]},$$
where $\Delta_i$ for $i>0$ (resp. $i=0$) is irreducible and parametrizes non-smooth stable curves $C$ such that the partial normalization of $C$ at one of its nodes consists of two connected components of arithmetic genera $i$ and $g-i$ (resp. is connected).

The following universal Noether formula is obtained by Mumford \cite[p. 102]{Mum} as an easy application of Grothendieck's Riemann--Roch theorem.

\begin{thm}[Noether formula] \label{noether}
In $\Pic(\overline\CM_{g})$,  
$$
12\lambda_{\CMM_g}= \pi_{g*}\pair{\omega_{\ol\CCC_g/ \CMM_{g}},\omega_{\ol\CCC_g/ \CMM_{g}}}+ \CO_{\CMM_g}(\Delta).
$$
\end{thm}

\subsubsection{Pull-back to arbitrary families} \label{sec bundle on stable}

Let $\overline S$ be an integral noetherian  scheme.
Let $\pi:\overline X\to \overline S$ be a stable relative curve of genus $g$. Assume that the generic fiber of $\overline X$ is smooth. 
By the moduli property, there is a morphism $\iota:\overline S\to \CMM_g$ such that $\overline X\to \overline S$ is isomorphic to the base change of the universal curve $\ol\CCC_g\to \CMM_g$.
We have the Hodge bundle
$$
\lambda_{\overline S}=\iota^*\lambda \simeq \det \pi_* \omega_{\overline X/\overline S}=\wedge^g \pi_* \omega_{\overline X/\overline S}.
$$
We also have the boundary divisors
$$
\Delta_{\overline S}=\iota^*\Delta, \quad \Delta_{\overline S,i}=\iota^*\Delta_i,\quad i=0,1,\cdots, [g/2].
$$
They are effective Cartier divisors on $\overline S$.
We have 
$$
\Delta_{\overline S}=\Delta_{\overline S,0}+\Delta_{\overline S,1}+\cdots+\Delta_{\overline S,[g/2]}.
$$

Via pull-back, the Noether formula in Theorem \ref{noether} gives
$$
12\lambda_{\overline S}= \pi_{*}\pair{\omega_{\overline X/\overline S},\omega_{\overline X/\overline S}}+ \CO(\Delta_{\overline S}).
$$
in $\Pic(\overline S)$.
An isomorphism 
$$
12\lambda_{\overline S}\lra  \pi_{*}\pair{\omega_{\overline X/\overline S},\omega_{\overline X/\overline S}}+ \CO(\Delta_{\overline S})
$$
of line bundles over $S$ is called \emph{semi-canonical} if it is obtained as the base change by the moduli morphism $S\to\CMM_g$ of an isomorphism 
$$
12\lambda_{\CMM_g}\lra \pi_{g*}\pair{\omega_{\ol\CCC_g/ \CMM_{g}},\omega_{\ol\CCC_g/ \CMM_{g}}}+ \CO_{\CMM_g}(\Delta).
$$
Here $\CMM_g$ is the moduli stack of stable curves of genus $g$ over $\ZZ$ in the setting of Theorem \ref{noether}. 
As $\overline\CM_{g}$ is proper and smooth over $\ZZ$ with geometrically connected fibers, we have $\Gamma(\overline\CM_{g}, \CO_{\overline\CM_{g}})=\ZZ$ and 
$\Gamma(\overline\CM_{g}, \CO_{\overline\CM_{g}}^\times)=\ZZ^\times=\{\pm 1\}$. 
Then the isomorphism for $12\lambda_{\CMM_g}$ is unique up to multiplication by $\{\pm 1\}$. As a consequence, a semi-canonical isomorphism for $12\lambda_{\overline S}$ is also unique up to multiplication by $\{\pm 1\}$.

If $\overline S$ is normal, the divisor $\Delta_i$ on $\overline S$ can be defined explicitly as follows. 
It suffices to define the multiplicity $\ord_v(\Delta_i)$ of $\Delta_i$ along a codimension one point $v$ of $\overline S$. 
The problem is further reduced to the curve $\overline X_{\CO_{\overline S,v}}$ over $\CO_{\overline S,v}$.
This is essentially contained in the example of \cite[p. 464]{CH}. 
In fact, by a finite unramified base change of $\CO_{\overline S,v}$, we can assume that all nodes of $\overline X_v$ are rational over the residue field $k(v)$. 
For any node $x$ of $\overline X_v$, denote its type by $i(x)\in \{0,1,\cdots, [g/2]\}$. 
The multiplicity $m(x)$ of $x$ is defined to be a positive integer such that the local equation of $\overline X$ at $x$ is of type $t_1t_2=\varpi^{m(x)}$, where $\varpi\in \CO_{\overline S,v}$ is a generator of the maximal ideal.
With these data, we have the formula
$$
\ord_v(\Delta_i) = \sum_{x\in \overline X_v,\, i(x)=i} m(x).
$$

\subsubsection{Equivalent definitions} \label{sec hodge bundle}

There is another equivalent definition of the Hodge bundle, which is essential in the definition of the Faltings height. 

Let $S$ be an integral noetherian  scheme, and let $\pi:X\to S$ be a stable relative curve of genus $g>1$.
By \cite[\S9.4, Thm. 1]{BLR}, the Picard functor  $\Pic_{X/S}$ is representable by an algebraic space over $S$, and more importantly the relative identity component $ J=\Pic_{X/S}^0$ is represented by a (separated and smooth) semi-abelian group scheme over $S$.
The \emph{Hodge bundle} of $ J$ over $ S$ is defined as 
$$
\underline\omega_{S}:=\det( e^*  \Omega_{ J/ S}^1),
$$
where $e: S\to  J$ is the identity section. 
Then we have the following result, which is well-known to the experts. 

\begin{lem}\label{hodge bundle equivalence}
There is a canonical isomorphism 
$$
i:\lambda_{S}\lra \underline\omega_{ S}.
$$
Moreover, if $S$ is a quasi-projective variety over $\CC$ and $\pi:X\to S$ is smooth, then the canonical isomorphism induces an equality 
$\|\cdot\|_{\rm det}=i^*\|\cdot\|_{\rm Fal}$
of hermitian metric of $\lambda_{S}$ on $S$, where the two metrics are defined as follows. 
\begin{itemize}
\item[(1)]  The Faltings metric $\|\cdot\|_{\rm Fal}$ of $\underline\omega_{S}$ on $S$ is defined such that for any point $s\in S(\CC)$ and any section
$$\alpha\in  \underline\omega_{ S}(s)= e_s^*\Omega_{J_s/s}^g\simeq  \Gamma(J_s, \Omega_{J_s/s}^g),$$
the metric gives
$$
\|\alpha\|_{\rm Fal}^2=\frac{i^{g^2}}{2^g} \int_{J_s} \alpha\wedge\bar \alpha.
$$

\item[(2)]
The determinant metric 
$\|\cdot\|_{\rm det}=\det \|\cdot\|_{\rm nat}$ induced by  the process
$\lambda_S=\det(\pi_*\omega_{X/S})$, where the natural metric $\|\cdot\|_{\rm nat}$ on 
$\pi_*\omega_{X/S}$ is defined  
such that for any point $s\in S(\CC)$ and any section
$$\beta\in  (\pi_*\omega_{X/S})(s)= \Gamma(X_s, \omega_{X_s/s}),$$
the metric gives
$$
\|\beta\|_{\rm nat}^2=\frac{i}{2} \int_{X_s} \beta\wedge\bar \beta.
$$
\end{itemize}
\end{lem}
\begin{proof}
The equality of the metrics follows from Szpiro \cite[Lem. 3.2.1]{Szp}.
For the canonical isomorphism, it suffices to establish a canonical isomorphism
$$
\pi_*\omega_{ X/ S} \lra e^*  \Omega_{ J/ S}^1.
$$
By duality, it suffices to establish a canonical isomorphism 
$$
\Lie( J/ S)\lra R^1\pi_*\CO_{ X}.
$$
This is a consequence of the deformation theory (cf. \cite[\S8.4, Thm. 1]{BLR}).
\end{proof}

Let $S$ be a flat and quasi-projective integral scheme over $\QQ$ or over $\ZZ$. 
Let $\pi:X\to S$ be a smooth relative curve of genus $g>1$.
Then the Hodge line bundle $\lambda_S$ and the determinant metric 
$\|\cdot\|_{\rm det}$ forms an adelic line bundle $\ol \lambda_S$ over $S/\ZZ$.
This follows from \cite[\S5.5]{YZ2} via the isometry in Lemma \ref{hodge bundle equivalence}.
We denote this adelic line bundle by $\ol\lambda_S$.

\subsubsection{Stable compactification} \label{sec stable compactification}

\kkk
Let $S$ be a quasi-projective integral scheme over $k$. 
Let $\pi:X\to S$ be a smooth relative curve of genus $g>1$.

By a \emph{stable compactification} of $\pi:X\to S$, we mean a projective integral scheme $\ol S$ over $k$ with an open immersion $S\to \ol S$, a stable relative curve $\ol\pi:\ol X\to \ol S$ of genus $g$, and an open immersion $X\to \ol X$ compatible with the previous morphisms. 

For any smooth relative curve $\pi:X\to S$ as above, there is an integral scheme $S'$ with a finite and surjective morphism $S'\to S$ such that the base change $\pi':X'\to S'$ of $\pi:X\to S$ by $S'\to S$ has a stable compactification. Moreover, $S'$ can be taken to be normal by a further normalization process.
If $\dim S=1$, this fact is a consequence of the stable reduction theorem (or the properness of the stack $\CMM_{g}$). 
If $\dim S>1$, we sketch a proof as follows.

We claim that there is a projective normal integral scheme $\CMM_g'$ over $\ZZ$, together with a (representable) finite and surjective morphism $\CMM_{g}'\to \CMM_{g}$. 
In fact, by the abstract result about Deligne--Mumford stacks in \cite[Thm. 16.6]{LMB}, there is a (representable) finite and surjective morphism $\CMM_{g}'\to \CMM_{g}$ from a scheme $\CMM_{g}'$. 
Then $\CMM_g'$ is a proper scheme over $\ZZ$, and we need to prove that it is actually projective over $\ZZ$. 
For that, the key is that the coarse moduli scheme 
$\overline M_g$ of $\CMM_g$ is projective over $\ZZ$ (cf. \cite[Thm. 7.2]{CLM}).
Then it suffices to prove that the composition $\CMM_{g}'\to \CMM_{g} \to \overline M_g$ is a finite morphism of schemes. 
Since the composition is a morphism of proper schemes over $\ZZ$, it suffices to 
check that it is quasi-finite, and thus it suffices to check that for any algebraically closed field $F$, the composition $\CMM_{g}'(F)\to \CMM_{g}(F) \to \overline M_g(F)$ as a map of sets has finite fibers. 
This holds since the first arrow has finite fibers and the second arrow is bijective. 
Thus $\CMM_g'$ is projective, and we can assume that it is integral and normal by passing to the normalization of a suitable irreducible component. 

The scheme $\CMM_{g}'$ is equipped with a ``universal'' stable curve $\overline\CCC_{g}'\to\CMM_g'$ defined as the base change of the universal curve $\overline\CCC_{g}\to \CMM_{g}$ by $\CMM_{g}'\to \CMM_{g}$.

Return to the existence of the stable compactification.
We can take $S'\to S$ to be the base change of $\CMM_{g}'\to \CMM_{g}$ by the moduli map $S\to \CMM_g$, and a stable compactification is given by the Zariski closure of $S'$ in $\CMM_g'$.
In this process, we may also need to replace $S'$ by a suitable irreducible component to make it integral.

As we have just seen, the caveat that neither $\CM_g$ nor $\CMM_g$ is (representable by) a scheme can be overcome by the above scheme $\CMM_g'$ with its universal family with the cost of passing to a finite surjective extension.
Another approach is to add a full level-$N$ structure (for $N\geq3$) to the definition of 
$\CM_g$ to get a (fine) moduli scheme $\CM_{g,N}$. 
This approach is more common and has the advantage of being explicit, but it has other issues including that it requires the base scheme to avoid the prime factors of $N$, and that there is no natural level-$N$ structure to add to $\CMM_g$. 
The second issue has a weak solution by van Geemen--Oort \cite{GO}, which constructs a stable compactification of the universal family of $\CM_{g,N}$.

In this paper, both approaches can solve our related issues. 
However, we will usually take the approach using level structures, since it is more common in the literature.

\subsection{Passing to the relative dualizing sheaf} \label{sec def E}

Our first step to treat $\pair{\overline\omega_{X/S,a},\overline\omega_{X/S,a}}$ in Theorem \ref{bigness1} is to convert it to the relative dualizing sheaf of a stable compactification.

Let $k$ be a field. Note that we exclude the case $k=\ZZ$ here.
Let $S$ be a quasi-projective normal integral scheme over $k$.
Let $\pi:X\to S$ a smooth relative curve of genus $g>1$ with a stable compactification 
$\overline\pi:\ol X\to \ol S$ over $k$.

Note that the relative dualizing sheaf $\omega_{\overline X/\overline S}$ is naturally viewed as an adelic line bundle on $X/k$ by the functor $\Picc(\ol X)\to \wh\Picc(X/k)$. 
Consider the difference 
$$
\overline\pi_*\pair{\omega_{\overline X/\overline S},\omega_{\overline X/\overline S}}
-\pi_*\pair{\overline\omega_{X/S,a},\overline\omega_{X/S,a}}
\simeq \CO(\ol E_S)
$$
in $\wh\Pic(S)$. 
Here we explain the right-hand side as follows. 
The underlying line bundle (over $S$) of the left-hand side is canonically isomorphic to the trivial bundle $\CO_S$.
The section $1\in \CO_S$ corresponds to a rational section $t$ of the left-hand side, and the adelic divisor 
$$\ol E_S:=\wh\div(t)$$ 
on $S/k$ has underlying divisor 0 on $S$. 

There is a reasonably explicit description of $\ol E_S$ by 
the graph-theoretic approach of \cite{Zha1}. 
In fact, to describe $\ol E_S$, it suffices to describe its image under the analytification map
$$\wh\Div(S/k)\lra \wh\Div(S^\an)$$
in \cite[Prop. 3.3.1]{YZ2}. 
As the underlying divisor of $\ol E_S$ on $S$ is 0, it suffices to describe its total Green's function 
$$\wt g_{\ol E_S}: S^\an \lra \RR.$$
By \cite[Lem. 3.1.1(1)]{YZ2}, the set of $v\in S^\an$ corresponding to discrete valuations of $k(S)$ is dense in $S^\an$. 
By continuity, it suffices to describe $\wt g_{\ol E_S}(v)$ for any $v\in S^\an$ corresponding to a discrete valuation.
In this case, by \cite[Thm. 4.4]{Zha1}, 
$$
\wt g_{\ol E_S}(v)/\log e_v= \epsilon(\Gamma_v)=\int_{\Gamma_v} g_\mu(x,x)((2g-2)\mu+\delta_{K_{X,v}}).
$$
Here $\Gamma_v$ is the reduction graph of $X$ over the discrete valuation field $H_v$, and $e_v=|\varpi_v|^{-1}$, where $\varpi_v$ is a generator of the maximal ideal of the valuation ring $O_{H_v}$. 
Note that the loc. cit. is about the case $e_v=e$, but the result can be transferred to general $e_v$ by the norm-equivariance property in \cite[\S3.2]{YZ2}. 

By this description, we have the following result. 
Recall that an adelic divisor is effective if it is the limit of effective model adelic divisors.

\begin{lem} \label{effective1}
The difference $(2g-2)\Delta_{\ol S}-\ol E_S$ is an effective adelic divisor in $\wh\Div(S/k)$.
\end{lem}
\begin{proof}
Note that both adelic divisors $\Delta_{\ol S}$ and $\ol E_S$ have trivial underlying divisors in $\Div(S)$.
By \cite[Lem. 5.1.2]{YZ2}, it suffices to prove that the corresponding total Green's function $(2g-2)\wt g_{\Delta_{\ol S}}-\wt g_{\ol E_S}$ is non-negative on $S^\an$.
By continuity and density, it suffices to prove that 
$(2g-2)\wt g_{\Delta_{\ol S}}(v)\geq \wt g_{\ol E_S}(v)$
for discrete valuations $v\in S^\an$. 
Note that the reduction graph $\Gamma_v$ is a metrized graph as recalled in \S\ref{sec graph theory}, 
and $\wt g_{\Delta_{\ol S}}(v)/\log e_v$ is exactly the total length $\ell(\Gamma_v)$ of edges of $\Gamma_v$.
Thus it suffices to check $\epsilon(\Gamma_v)\leq (2g-2)\ell(\Gamma_v)$.

This is basic in graph theory. 
In fact, by \cite[eq. (4.1.4)]{Zha3},
$$
\epsilon(\Gamma_v)=\int_{\Gamma_v\times \Gamma_v} r(x,y) \delta_{K_{X,v}}(x) \mu(y).
$$
Here $r(x,y)$ is the resistance function defined right before \cite[Prop. 3.3]{Zha1} or equivalently in \cite[\S6, Def. 8]{BF}.
The inequality follows from the bound $r(x,y)\leq \ell(\Gamma_v)$, which in turn follows from the bound $r(x,y)\leq \ell(L_{x,y})$, where $L_{x,y}$ is a path in $\Gamma_v$ connecting $x,y$. 

The intuitive bound $r(x,y)\leq \ell(L_{x,y})$ is proved in \cite[Thm. D]{KR} as a consequence of Rayleigh's monotonicity law. 
To fit our terminology, we sketch a proof following the idea of \cite[\S6, Exer. 12]{BF}.
The key is the triangle inequality 
$$
r(x,w)\leq r(x,z)+r(z,w), \quad x,y,z\in \Gamma_v.
$$
In fact, set $y=x$ in \cite[\S6, Thm. 8]{BF}, we have an equality
$$
j_z(x,x)-j_z(w, x) = j_w(x, x)-j_w(z, x).
$$
Here the $j$-function is a generalized Green's function defined in \cite[\S6, Cor. 3]{BF}. 
In particular, $r(x,y)=j_y(x,x)=j_x(y,y)$ as in \cite[\S6, Def. 8]{BF}.
By \cite[\S6, Exer. 9]{BF}, $j_z(w, x)\geq 0$ and $j_w(z, x)\leq j_w(z, z)$, so the equality implies 
$$
j_z(x,x)\geq  j_w(x, x)-j_w(z, z).
$$
This gives the triangle inequality.

With the triangle inequality, we can reduce the proof of $r(x,y)\leq \ell(L_{x,y})$ to the case that $L_{x,y}$ is an edge of the graph $\Gamma_v$ with vertices $x,y$.
Here we can assume that $x,y$ lie in the vertex set of $\Gamma_v$ by taking a subdivision. 
Then the inequality follows from \cite[\S6, Thm. 10]{BF}.
\end{proof}

It is reasonable to say that $\ol E_S$ is the 
\emph{globalization} of the $\epsilon$-invariant of \cite{Zha1}, and 
$\Delta_{\ol S}$ is the \emph{globalization} of the $\ell$-invariant (or the $\delta$-invariant).
In the following, we will see an adelic divisor $\ol \Phi_S$, the \emph{globalization} of the $\varphi$-invariant of Zhang \cite{Zha3}.

\subsection{Lower bound of the pairing}

We also need a lower bound of 
$\pair{\overline\omega_{X/S,a},\overline\omega_{X/S,a}}$, which is given by a globalization of the $\varphi$-invariant of Zhang \cite{Zha3}.

\subsubsection{Zhang's $\varphi$-invariant}

We first recall the $\varphi$-invariant of curves over local fields from \cite{Zha3}. 
Let $C$ be a smooth projective curve of genus $g>1$ over a complete field $K$. 
The invariant $\varphi(C)$ is defined as the local intersection number
$$
\varphi(C)
=-\int_{(C\times C)^\an}g_{\Delta,a} \, c_1(\CO(\ol \Delta)_a)^2.
$$
As in \cite[Thm 1.3.1]{Zha3}, it has explicit expressions as follows.

If $K$ is non-archimedean with a discrete valuation, then 
$$
\varphi(C)=-\frac{1}{4}\delta(C)+\frac{1}{4}\int_{\Gamma(C)} g_{\mu}(x,x) ((10g+2) \mu-\delta_{K_{C}}).
$$
Here the terms in the integration are as in \S\ref{sec graph theory}.
The invariant $\delta(C)$ is the total length $\ell(\Gamma)$ of edges of the reduction graph $\Gamma$ of $C$ over $K$, viewed as a metrized graph as recalled in \S\ref{sec graph theory}.
If $C$ has split semistable reduction, $\delta(C)$ is also the total number of nodes on the minimal regular model of $C$ over $O_K$;
in general, if $K'$ is a finite extension of $K$ such that $C$ has split semistable reduction over $K'$, then $\delta(C)=\delta(C_{K'})/[K':K]$.
We refer to \cite[Lem 3.5.4]{Zha3} for more details on this expression. 

If $K=\RR$, define $\varphi(C)=\varphi(C_\CC)$. 
If $K=\CC$,  the computation of \cite[Prop 2.5.3]{Zha3} gives
$$
\varphi(C)
=\sum_\lambda \sum_{j=1}^g\sum_{k=1}^g \frac{2}{\lambda}\left|  \int_C \varphi_\lambda\, \omega_j\wedge \overline\omega_{k}  \right|^2.
$$
Here $\omega_1,\cdots, \omega_g$ is an orthonormal basis of $\Gamma(C, \Omega_C^1)$ with respect to the inner product
$$
\pair{\alpha, \beta}=\frac{i}{2}\int_C \alpha\wedge \ol \beta,
$$ 
the first summation goes over all positive eigenvalues $\lambda$ of the Laplacian operator 
$$\Delta_{d\mu} f= (\frac{1}{\pi i}\partial\overline\partial f) /d\mu$$
over $C^\infty(C_\CC)$,
and $\varphi_\lambda$ is an eigenvector of $\lambda$ normalized such that 
$\{\varphi_\lambda\}_\lambda$ is orthonormal with respect to the inner product
$$
\pair{f_1,f_2}=\int_C f_1\overline f_2 d\mu.
$$
We refer to \cite[Prop 2.5.3]{Zha3} for more details on this expression.

\subsubsection{Globalization $\ol\Phi_S$ of the $\varphi$-invariant}

\kkk
Let $S$ be a quasi-projective and flat normal integral scheme over $k$.
Let $\pi:X\to S$ be a smooth relative curve of genus $g>1$. 

Consider the morphism $(\pi,\pi):X\times_SX\to S$ and its diagonal divisor $\Delta:X\to X\times_SX$. 
There are canonical isomorphisms
$$
(\pi,\pi)_*\pair{\CO(\Delta),\,\CO(\Delta),\, \CO(\Delta)}
\lra \pi_*\pair{\Delta^*\CO(\Delta),\Delta^*\CO(\Delta)}
\lra \pi_*\pair{\omega_{X/S},\omega_{X/S}}.
$$
This defines a section $s$ of the underlying line bundle of the adelic line bundle 
$$
\pi_*\pair{\ol\omega_{X/S,a},\ol\omega_{X/S,a}}-(\pi,\pi)_*\pair{\ol\CO(\Delta)_a,\,\ol\CO(\Delta)_a,\,\ol\CO(\Delta)_a}.
$$
Via this adelic line bundle, define 
$$\ol\Phi_S=\wh\div(s),$$ 
viewed as an adelic divisor on $S$. 
The underlying divisor of $\ol\Phi_S$ is 0 by definition. 
By the integration formula, 
the total Green's function $\wt g_{\ol\Phi_S}$ on $S^\an$
at any discrete or archimedean valuation $v\in S^\an$
is given by 
$$
\wt g_{\ol\Phi_S}(v)=\int_{X_v^\an} \log \|1\|_{\Delta, a} c_1(\ol\CO(\Delta)_a)^2.
$$
By \cite[Prop. 2.5.3, Lem. 3.5.4]{Zha3}, this is exactly the $\varphi$-invariant of 
$X_{H_v}$ in the standard case $e_v=e$.
It transfers to general $e_v$ by the norm-equivariance property in \cite[\S3.2]{YZ2}.

\subsubsection{The lower bound}

\kkk
Let $S$ be a quasi-projective and flat normal integral scheme over $k$.
Let $\pi:X\to S$ be a smooth relative curve of genus $g>1$.
Denote by $(\pi,\pi):X\times_SX\to S$ the structure morphism as before. 

Let $J=\Pic_{X/S}^0$ be the Jacobian scheme.
Recall that we have a symmetric and relatively ample line bundle $\Theta$ on $J$ by Definition \ref{qqq}. 
Recall that $\OTheta$ is the nef adelic line bundle on $J/k$ extending $\Theta$ such that $[2]^*\OTheta=4\OTheta$ in $\wh\Pic(J)$. 

Let 
$j$ be the morphism 
$$j:X\times_SX \lra J, \quad (x,y)\longmapsto y-x.$$
Recall that in Theorem \ref{isomorphism3}(2), we have obtained a formula
$$
j^*\OTheta=2\ol\CO(\Delta)_a+p_1^*\overline\omega_{X/S,a}+p_2^*\overline\omega_{X/S,a}
$$
in $\wh\Pic(X)_\QQ$.
The following result computes the Deligne pairing $(\pi,\pi)_*\pair{j^*\OTheta,j^*\OTheta,j^*\OTheta}$, which is nef by the nefness of $\OTheta$. 
It is a family version of and inspired by de Jong \cite[Thm. 8.1]{dJo3}, where the latter was in turn based on the main result of Zhang \cite{Zha3}. 

As above, for a discrete valuation $v\in S^\an$, denote $e_v=|\varpi_v|_v^{-1}$, where $\varpi_v$ is a generator of the maximal ideal of the valuation ring $O_{H_v}$;
for an archimedean valuation $v\in S^\an$, denote $e_v=|e|_v$.

\begin{thm} \label{lower bound}
\kkk
There is an identity in $\wh\Pic(S)_\QQ$ given by 
$$
(\pi,\pi)_*\pair{j^*\OTheta,j^*\OTheta,j^*\OTheta}
=(12g-4)\pi_*\pair{\overline\omega_{X/S,a},\overline\omega_{X/S,a}}
- 8\CO(\ol \Phi_S).
$$
\end{thm}

\begin{proof}
The proof follows the framework of \cite[Thm. 8.1]{dJo3}.
Denote $\ol\omega=\overline\omega_{X/S,a}$ and 
$\ol\omega_i=p_i^*\overline\omega_{X/S,a}$ for $i=1,2$. 
Expand
$(\pi,\pi)_*\pair{j^*\OTheta,j^*\OTheta,j^*\OTheta}$ by the formula
$$
j^*\OTheta=2\ol\CO(\Delta)_a+\overline\omega_1+\overline\omega_2.
$$
We have
\begin{eqnarray*}
&&\pair{j^*\OTheta,\,j^*\OTheta,\,j^*\OTheta} \\
&=& 8\pair{\ol\CO(\Delta)_a,\,\ol\CO(\Delta)_a,\,\ol\CO(\Delta)_a}
+ 12\pair{\ol\CO(\Delta)_a,\,\ol\CO(\Delta)_a,\,\overline\omega_1+\overline\omega_2}\\
&&+ 6\pair{\ol\CO(\Delta)_a,\, \overline\omega_1+\overline\omega_2,\,\overline\omega_1+\overline\omega_2}
 + \pair{\overline\omega_1+\overline\omega_2,\, \overline\omega_1+\overline\omega_2,\, \overline\omega_1+\overline\omega_2}.
\end{eqnarray*}
Here we abbreviate the Deligne pairing $(\pi,\pi)_*\pair{\cdot, \cdot, \cdot}$ as $\pair{\cdot, \cdot, \cdot}$ for simplicity, but we will not abbreviate $\pi_*\pair{\cdot, \cdot}$ to avoid confusion. 
For further simplifications, we need the following results:
\begin{enumerate}[(1)]
\item $\pair{\overline\omega_1,\, \overline\omega_1,\, \overline\omega_1}=\pair{\overline\omega_2,\, \overline\omega_2,\, \overline\omega_2}=0.$

\item $\pair{\overline\omega_1,\, \overline\omega_1,\, \overline\omega_2}=\pair{\overline\omega_2,\, \overline\omega_2,\, \overline\omega_1}=(2g-2)\pi_*\pair{\ol\omega,\ol\omega}.$

\item $\pair{\ol\CO(\Delta)_a,\, \overline\omega_i,\, \overline\omega_j}=\pi_*\pair{\ol\omega,\ol\omega}, \quad i,j\in\{1,2\}.$

\item $\pair{\ol\CO(\Delta)_a,\, \ol\CO(\Delta)_a,\, \overline\omega_1}=\pair{\ol\CO(\Delta)_a,\, \ol\CO(\Delta)_a,\, \overline\omega_2}=
-\pi_*\pair{\ol\omega,\ol\omega}.$

\item $\pair{\ol\CO(\Delta)_a,\,\ol\CO(\Delta)_a,\,\ol\CO(\Delta)_a}=
\pi_*\pair{\ol\omega,\ol\omega}-\CO(\ol\Phi_S).
$
\end{enumerate}

It is easy to derive the result from (1)-(5).
In fact, (1)-(4) gives 
$$\pair{j^*\OTheta,\,j^*\OTheta,\,j^*\OTheta} 
= 8\pair{\ol\CO(\Delta)_a,\,\ol\CO(\Delta)_a,\,\ol\CO(\Delta)_a}
+ 6(2g-2) \pi_*\pair{\ol\omega,\ol\omega}.$$
Then apply (5). 

Note that (5) follows from the definition of $\ol\Phi_S$.
The idea to prove the identities in (1)-(4) is to apply \cite[Lem. 4.6.1(2)]{YZ2} to 
$$X\times_SX\stackrel{p_i}{\lra }X\stackrel{\pi}{\lra}S.$$
This gives (for $i=1$)
$$
\pair{\overline\omega_1,\, \overline\alpha,\, \overline\beta}
=\pi_*\pair{\overline\omega, p_{1*}\pair{\overline\alpha,\, \overline\beta}}.
$$
If $\alpha=\omega_1$, then we can apply \cite[Lem. 4.6.1(2)]{YZ2} to compute $p_{1*}\pair{\overline\alpha,\, \overline\beta}$. 
This solves (1), (2) and (3) with $i=j$. 
The remaining case of (3) is a consequence of 
$p_{1*}\pair{\overline\omega,\, \ol\CO(\Delta)_a}=\ol\omega$ from Theorem \ref{admissible2}(2). 

For (4), it suffices to check 
$$
p_{1*}\pair{\ol\CO(\Delta)_a,\, \ol\CO(\Delta)_a}
=-\overline\omega.
$$
There is a canonical isomorphism 
$\Delta^*\ol\CO(\Delta)_a \to -\overline\omega$
by Theorem \ref{admissible2}(1). 
It suffices to have a canonical isomorphism
$$
p_{1*}\pair{\ol\CO(\Delta)_a,\, \ol\CO(\Delta)_a}
\lra \Delta^*\ol\CO(\Delta)_a.
$$
This follows from Proposition \ref{adjunction}(1) applied to the section $\Delta:X\to X\times_SX$ of the smooth relative curve $p_1:X\times_SX\to X$. 
Alternatively, it is also easy to modify the proof of Proposition \ref{adjunction}(1) to the current situation. 
In fact, by the analytification functor $\wh\Picc(X)\to \wh\Picc(X^\an)$, it suffices to prove that the norm of the above morphism is $1$ at any $v\in X^\an$. 
By the integration formulae in \cite[Thm. 4.6.2]{YZ2} and \cite[\S4.2.2]{YZ2},
the logarithm of the norm is equal to
$$
-\int_{X_{H_v}^\an} g_{\Delta,a}(x,\cdot) c_1(x^*\ol\CO(\Delta)_a). 
$$ 
Here $H_v$ is the completed residue field of $v\in X^\an$, $x:\Spec H_v\to X$ is the point corresponding to $v$, and $g_{\Delta,a}$ is the admissible Green's function as in Theorem \ref{admissible1}. 
The integral is 0 by Theorem \ref{admissible1}(1)(2). 
This finishes the proof.
\end{proof}

\subsubsection{Cinkir's bound}

The following result essentially follows from Cinkir \cite{Cin1}, which was conjectured by Zhang \cite{Zha3}.

\begin{thm} \label{effective2}
Let $S$ be a quasi-projective normal integral scheme over a field $k$.
Let $\pi:X\to S$ be a smooth relative curve of genus $g>1$ with a stable compactification $\ol\pi:\ol X\to \ol S$ over $k$.
Then $\ds\overline\Phi_S-\frac{1}{39}\Delta_{\ol S}$ is an effective adelic divisor in $\wh\Div(S)$.
\end{thm}
\begin{proof}
This is similar to Lemma \ref{effective1}. 
It suffices to prove  
$$\wt g_{\ol \Phi_S}(v)-\frac{1}{39}\wt g_{\Delta_{\ol S}}(v)\geq 0$$
for any discrete valuations $v\in S^\an$. 
By \cite[Thm. 2.11]{Cin1},
$$
\varphi(\Gamma_v)\geq \frac{1}{39} \sum_{i=0}^{[g/2]}\delta_i(\Gamma_v)=\frac{1}{39}\ell(\Gamma_v).
$$
This finishes the proof. 

For convenience of readers, we explain the term $\delta_i(\Gamma_v)$ briefly. It is a term for polarized metrized graphs, but in our case it is equal to the total number of nodes of type $i$ on the geometric special fiber of the minimal regular model of $X_{H_v}$ over $O_{H_v}$. Here $H_v$ denotes the completed residue field of $v\in S^\an$ as before. Recall that a node on a semistable curve of arithmetic genus $g$ over an algebraically closed field is of type $i$ for some $i>0$ if the partial normalization of the semistable curve at this node consists of two connected components of arithmetic genera $i$ and $g-i$;
it is of type $0$ if the partial normalization is connected.
This is compatible with the definition of $\Delta_i$ in
\S\ref{sec basic moduli}.
\end{proof}

\subsubsection{A special adelic line bundle}

The $\varphi$-invariant was introduced by Zhang \cite{Zha3} to study the Beilinson--Bloch height of the Gross--Schoen cycle.
Then our adelic divisor $\ol\Phi_S$ can be used to study the variation of the height.
The exposition here will not be used elsewhere in this paper. 

Let $K$ be either a number field, or the function field of one variable over a field $k$.
Let $C$ be a smooth projective curve over $K$ of genus $g>1$. 
Let $\xi$ be a line bundle on $C$ such that $(2g-2)\xi=\omega_{C/K}$ in $\Pic(C)$, which exists by replacing $K$ by a suitable finite extension.
By \cite[Thm. 1.3.1]{Zha3}, 
$$
\pair{\Delta_\xi(C),\Delta_\xi(C)}_{\mathrm{BB}}=\frac{2g+1}{2g-2}
\ol\omega_{C/K,a}^2-\sum_{v}\varphi(C_v)\deg(v).
$$
Here $\Delta_\xi(C)$ is the Gross--Schoen cycle, the modified diagonal cycle on $C^3$, 
 and the left-hand side is its Beilinson--Bloch height. 

Let $S$ be a quasi-projective normal variety over $K$.
Let $\pi:X\to S$ be a smooth relative curve over $S$ of genus $g>1$. 
Define an adelic $\QQ$-line bundle $\OL$ on $S$ by  
$$
\OL=\frac{2g+1}{2g-2} \pi_*\pair{\ol\omega_{X/S,a},\ol\omega_{X/S,a}}
-\CO(\ol\Phi_S). 
$$
Then the height function $h_{\OL}:S(\ol K)\to \RR$ exactly gives
$$
h_{\OL}(y)=\pair{\Delta_\xi(X_y),\Delta_\xi(X_y)}_{\mathrm{BB}},\quad y\in S(\ol K). 
$$
This generalizes the first statement of \cite[Thm. 1.3.5]{Zha3} to the quasi-projective setting.

It is natural to ask wether $\OL$ is nef on $S$. 
If $\pi:X\to S$ has maximal variation, 
we further ask whether $\OL$ is big on $S$.

\subsection{Bigness in the geometric case} \label{sec bigness geometric}

In the geometric case that $k$ is a field, Theorem \ref{bigness1} is a consequence of the following result. 

\begin{thm} \label{bigness2}
Let $k$ be a field. 
Let $S$ be a quasi-projective normal integral scheme over $k$.
Let $\pi:X\to S$ be a smooth relative curve of genus $g>1$ with a stable compactification 
$\ol\pi:\ol X\to \ol S$ over $k$.
In $\wh\Pic(S)_\QQ$, we have 
$$
\pi_*\pair{\overline\omega_{X/S,a},\overline\omega_{X/S,a}}
= \frac{3}{5(2g-1)(3g-1)} \lambda_{\ol S}
+\ol A+\CO(\ol B)
$$ 
for a nef adelic $\QQ$-line bundle $\ol A \in \wh\Pic(S)_{\QQ}$ and an effective adelic divisor $\ol B\in \wh\Div(S)$ with underlying divisor $B=0$ in $\Div(S)$.
\end{thm}
\begin{proof}
This is a combination of all the major results of this section. 
In fact, we first have in $\wh\Pic(S)_\QQ$,
\begin{align*}
\pi_*\pair{\overline\omega_{X/S,a},\overline\omega_{X/S,a}}
= &\ \overline\pi_*\pair{\omega_{\overline X/\overline S},\omega_{\overline X/\overline S}}
- \CO(\ol E_S) \\
= &\ 12\lambda_{\overline S}- \CO(\ol E_S+\Delta_{\overline S}) \\
=&\ 12\lambda_{\overline S}- (2g-1)\CO(\Delta_{\overline S})+\eff.
\end{align*}
Here ``$\eff$'' denotes the adelic line bundle associated to an effective adelic divisor in 
$\wh\Div(S)$ with underlying divisor $0$ in $\Div(S)$.
In fact, the first equality follows from the definition of $\ol E_S$ in \S\ref{sec def E}, 
the second equality follows from the pull-back of the Noether formula in Theorem \ref{noether} described in \S\ref{sec bundle on stable}, 
and the last inequality follows from Lemma \ref{effective1}. 

On the other hand, by Theorem \ref{lower bound} and Theorem \ref{effective2},
$$
\pi_*\pair{\overline\omega_{X/S,a},\overline\omega_{X/S,a}}
= \frac{2}{3g-1}\CO(\ol \Phi_S)+\nef
= \frac{2}{39(3g-1)}\CO(\Delta_{\ol S})+\nef+\eff.
$$
Here ``$\nef$'' denotes a nef element of
$\wh\Pic(S)_\QQ$.

To cancel $\CO(\Delta_{\ol S})$, a positive linear combination of these equalities gives 
$$
\left(1+ \frac{39(3g-1)}{2} (2g-1) \right)\pi_*\pair{\overline\omega_{X/S,a},\overline\omega_{X/S,a}}
=12\lambda_{\overline S}+\nef+\eff.
$$
This implies the result after replacing the coefficient on the left-hand side by the slightly larger number
$20(3g-1)(2g-1)$. 
\end{proof}

\begin{remark}
The constant $\displaystyle\frac{3}{5(2g-1)(3g-1)}$ in the theorem can be improved to
$\displaystyle\frac{1}{12}$ (as in Theorem \ref{fiberwise2})
by the method of \cite[Prop. 6.1]{LSW} (originally from \cite[Thm. 1.2]{Wil3}) combining with the Hodge index theorem of \cite{Car} (originally from \cite{YZ3}). 
We refer to the proof of Theorem \ref{fiberwise2} for how to obtain this constant.
\end{remark}

It is easy to see how Theorem \ref{bigness2} implies Theorem \ref{bigness1} in the geometric case. The nefness part is already proved in Theorem \ref{isomorphism3}(1). For the bigness part, by passing to a finite extension of $S$, we can assume that 
$\pi: X\to S$ has a stable compactification $\ol\pi:\overline X\to \overline S$ over $k$.
The normality of $S$ can be kept by taking a further normalization. 
This process uses the fact that for any finite and surjective morphism $S'\to S$ from an integral scheme $S'$, a nef adelic line bundle on $S$ is big if and only if its pull-back to $S'$ is big, which follows from the projection formula for top intersection numbers in \cite[Prop. 4.1.2]{YZ2}.

Once we have the stable compactification, 
 the key fact is that the Hodge bundle $\lambda_{\ol S}$ is nef and big on $\ol S$. 
If so, then 
$$\ds\frac{3}{5(2g-1)(3g-1)} \lambda_{\ol S}+\ol A$$
is nef and big on $S$, and 
$$\frac{3}{5(2g-1)(3g-1)} \lambda_{\ol S}+\ol A+\CO(\ol B)$$
 is big on $S$.

The nefness and bigness of the Hodge bundle $\lambda_{\ol S}$ on $\ol S$ are well known to the experts. 
For lack of a direct reference, we explain a proof by modifying the setting of \cite{Knu} slightly. 
Denote by $A_g$ the coarse moduli scheme of principally polarized abelian varieties 
of dimension $g$ over $k$.
Denote by $\overline A_g$ the minimal compactification of $A_g$.
By \cite[V, Thm. 2.3]{FC}, the Hodge $\QQ$-line bundle $\omega_{\overline A_g}$ of 
$\overline A_g$ is ample. 
Denote by $t:\overline M_g\to \overline A_g$ the Torelli morphism induced by the functor of taking Jacobian schemes. 
The definition is similar to the map at the top of \cite[p. 211]{Knu}.
The Hodge bundle $\lambda_{g}$ on the stack $\overline \CM_g$ descends to a $\QQ$-line bundle $\lambda_{\overline M_g}$ on $\overline M_g$. 
As in \cite[p. 211]{Knu}, we have 
$\lambda_{\overline M_g}=t^*\omega_{\overline A_g}$.
Finally, let $\pi:\ol X\to \ol S$ be as in the theorem. 
This gives a morphism $\iota:\ol S\to \overline M_g$, whose restriction to $S$ is generically finite. Note that $t:\overline M_g\to \overline A_g$ is finite on 
$M_g$, so the composition $t\circ\iota:\ol S\to \overline A_g$ is generically finite. 
Then $\lambda_{\ol S}=(t\circ\iota)^*\omega_{\overline A_g}$ on $\ol S$ is nef and big.

\subsection{Bigness in the arithmetic case} \label{sec arithmetic bigness}

In the arithmetic case $k=\ZZ$, Theorem \ref{bigness1} is consequence of the geometric case $k=\QQ$ and the following result. 

\begin{thm} \label{positive lower bound}
For any integer $g>1$, there is a constant $c_0>0$ depending only on $g$ such that $\varphi(C)\geq c_0$ for any compact Riemann surface $C$ of genus $g$.
\end{thm}

Due to its arithmetic importance, asymptotic behavior of degeneration of the 
$\varphi$-invariant is widely studied in the literature. See \cite[Thm. 1.1]{dJo2} for a precise asymptotic formula for $g=2$, which also gives a conjectural formula for $g>2$. 
For $g\geq2$ and degeneration to isolated singularities, the asymptotic formula was later proved by \cite[Thm. 7.1]{JS2} and \cite[Cor. 1.2]{Wil2}. 
These imply Theorem \ref{positive lower bound} in the case $g=2$ and in the case of 1-parameter families of Riemann surfaces of genus $g\geq2$.
As we will see later, our proof of the theorem takes a different approach, which asserts that the $\varphi$-invariant goes to infinity under degeneration, but does not give a precise asymptotic behavior. 

Let us first prove the arithmetic case of Theorem \ref{bigness1} assuming its geometric case and Theorem \ref{positive lower bound}.
Let $\pi:X\to S$ be as in Theorem \ref{bigness1}.
Let $c_0>0$ be as in Theorem \ref{positive lower bound}.
Denote by $\CO(c_0)$  the trivial line bundle on $k=\ZZ$ with metric given by $\|1\|=e^{-c_0}$.
View $\CO(c_0)$ as an adelic line bundle on $S$ by pull-back. 
By Theorem \ref{lower bound} and Theorem \ref{positive lower bound}, we have 
$$
\pi_*\pair{\overline\omega_{X/S,a},\overline\omega_{X/S,a}}
= \frac{2}{3g-1}\CO(\ol\Phi_S)+\nef
= \frac{2}{3g-1}\CO(c_0)+\nef+\eff.
$$

Denote $\OL=\pi_*\pair{\overline\omega_{X/S,a},\overline\omega_{X/S,a}}$, and denote by $\wt L$ the image of $\OL$ under the canonical map 
$\wh\Pic(S/\ZZ)\to \wh\Pic(S_\QQ/\QQ)$. Recall that the map is defined in \cite[\S2.5.5]{YZ2}, as the limit of the canonical map $\wh\Div(\CX)\to \Div(\CX_\QQ)$ for projective models $\CX$ of $S$ over $\ZZ$. 
Denote $d=\dim S$. 
Since $\OL$ is nef, we have
$$
\OL^d\geq \OL^{d-1}\cdot \frac{2}{3g-1}\CO(c_0)= \frac{2c_0}{3g-1}  \wt L^{d-1}>0. 
$$
Here the last inequality uses Theorem \ref{bigness1} for the geometric case $X_\QQ\to S_\QQ$ over $\QQ$. 
This proves the arithmetic case of Theorem \ref{bigness1}.

\medskip

It remains to prove Theorem \ref{positive lower bound}. 
Our proof uses information of the arithmetic divisor $\ol\Phi_S$ in both the geometric case and the arithmetic case.

\begin{proof}[Proof of Theorem \ref{positive lower bound}]
Fix an integer $N\geq 3$, and denote by $\CM=\CM_{g,N,\QQ}$ the (fine) moduli scheme of smooth curves of genus $g$ over $\QQ$ with a full level-$N$ structure.
Then $\CM$ is a smooth quasi-projective variety over $\QQ$, which follows from the GIT construction in \cite[\S7.4]{MFK}.  
By \cite[Thm. 2.1]{GO}, there is a projective compactification $\overline\CM=\overline\CM_{g,N,\QQ}$ of $\CM$ together with a tautological stable relative curve $\ol\CCC\to \overline\CM$. 

The $\varphi$-invariant defines a function $\varphi:\CM(\CC)\to \RR$. 
Note the values of $\varphi$ are strictly positive from Remark 1 after \cite[Prop 2.5.3]{Zha3}.
We will see that $\varphi$ is continuous on $\CM(\CC)$ and tends to infinity around the boundary $\overline\CM(\CC)\setminus \CM(\CC)$. 
These are sufficient to imply a positive lower bound of $\varphi$.

By Theorem \ref{lower bound}, there is an adelic divisor $\ol\Phi\in \wh\Div(\CM/\ZZ)$ with underlying divisor $0$ on $\CM$, such that the value $\wt g_{\ol\Phi}(v)$ of the total Green's function $\wt g_{\ol\Phi}$ at any discrete or archimedean valuation $v\in S_v^\an$ is equal to $\varphi(X_{H_v})$.
Then $\varphi:\CM(\CC)\to \RR$ is continuous by the continuity of $\wt g_{\ol\Phi}$. 

Denote by $\wt\Phi$ the image of $\ol\Phi$ under the canonical composition 
$$\wh\Div(\CM/\ZZ)\to \wh\Div(\CM/\QQ)\to \wh\Div(\CM_\CC/\CC).$$ 
Recall that the map is defined in \cite[\S2.5.5]{YZ2}, as the limit of the canonical map $\wh\Div(\CX)\to \Div(\CX_\CC)$ for projective models $\CX$ of $\CM$ over $\ZZ$. 

Apply Theorem \ref{effective2} to $k=\CC$. We have that
$\wt\Phi-\frac{1}{39}\Delta_{\ol \CM_\CC}$ is an effective adelic divisor in 
$\wh\Div(\CM_\CC/\CC)$.
This is sufficient to imply that $\varphi$ tends to infinity around the boundary 
$\Delta_{\ol \CM_\CC}=\overline\CM(\CC)\setminus \CM(\CC)$. 
To see that, we will use the theory of adelic divisors on the pair $\OB=(\CC,|\cdot|)$ in 
\cite[\S2.7]{YZ2}, which can be avoided by still using the theory of adelic divisors on $\ZZ$, but its use is more natural and makes the situation clear. 

Let us first recall the theory of adelic divisors on $\OB=(\CC,|\cdot|)$. 
First, the group of model divisors are defined by 
$$\wh\Div(\CM_\CC/\OB)_{\rm mod,\QQ}=\varinjlim_\CX \wh\Div(\CX/\OB)_\QQ,$$
where the limit is over the system of projective models $\CX$ of $\CM_\CC$ over $\CC$. 
Here $\wh\Div(\CX/\OB)$ is the group of pairs $(D,g_D)$, where $D$ is a Cartier divisor on $\CX$ and $g_D$ is a Green's function of $D$ on $\CX$.

Second, take a boundary divisor $(\CMM_\CC, \OE_0)$ with $\OE_0=(E_0,g_0)$, where $E_0=\Delta_{\ol \CM_\CC}$, and $g_0$ is a strictly positive Green's function of $E_0$ on $\CMM_\CC$.

Third, $\wh\Div(\CM_\CC/\OB)_{\QQ}$ is defined to be the completion of 
$\wh\Div(\CM_\CC/\OB)_{\rmod,\QQ}$ with respect to the topology induced by $(\CMM_\CC, \OE_0)$. 

Return to $\varphi$. 
The image of $\ol\Phi$ under the canonical map $\wh\Div(\CM/\ZZ)\to \wh\Div(\CM_\CC/\OB)_{\QQ}$ is represented by the pair
$(\wt\Phi, \varphi)$. 
By the process of the completion, there is a sequence $\{(D_i,g_i)\}_{i\geq1}$ in $\wh\Div(\CM_\CC/\OB)_{\rmod,\QQ}$ such that 
$$
-a_i (E_0,g_0)\leq (\wt\Phi, \varphi)-(D_i,g_i)\leq  a_i (E_0,g_0)
$$
for a sequence of rational numbers $\{a_i\}_{i\geq1}$ converging to 0.
This gives a bound 
$$
\varphi\geq g_i- a_i g_0.
$$

Fix an $i$ with $a_i<1/78$. 
Note that $g_i- a_i g_0$ is a Green's function of 
$$
D_i-a_i E_0
\geq (\wt\Phi-a_i E_0)-a_i E_0
\geq (1/39-2a_i)E_0>0.
$$
As a consequence, $g_i- a_i g_0$ goes to infinity around $E_0$, and so does $\varphi$.
This finishes the proof.
\end{proof}

\begin{remark} \label{song's work}
By more detailed analysis of the above proof, we can describe the growth of $\varphi$ along the boundary of $\CM(\CC)$ by the graph-theoretic data encoded from $\wt\Phi$.
See the recent work of Song \cite{Son} in this direction, which confirms a variant of a conjecture of de Jong \cite[Conj. 1.2]{dJo2}. 
\end{remark}

\section{Uniform Bogomolov-type theorem} \label{sec 4}

The main goal of this section is to prove the uniform Bogomolov-type theorem. 
In \S\ref{sec potential}, we introduce a notion of potential bigness, prove a key criterion for it, and consider a quick consequence of this notion on small points. 
In \S\ref{sec big examples}, we introduce three examples of potential bigness, as consequences of the bigness of the admissible canonical bundle in Theorem \ref{bigness5}. The exposition of the remaining subsections are based on these examples.
In \S\ref{sec small points} and \S\ref{sec small points2}, we state and prove our uniform Bogomolov-type theorem.
In \S\ref{sec pointwise}, we prove a uniform fiberwise bigness of the admissible canonical bundle.
Finally, in \S\ref{sec non-degeneracy}, we deduce some consequences on non-degenerate subvarieties and on the relative Bogomolov conjecture.

\subsection{Potentially big line bundles} \label{sec potential}

\kkk
Let $S$ be a quasi-projective and flat normal integral scheme over $k$. 
Let $Y$ be an integral scheme with a projective and flat morphism $\pi:Y\to S$ over 
$k$. 
Assume that the generic fiber of $\pi:Y\to S$ is geometrically integral.
Let $\OL$ be an adelic line bundle on $Y/k$.

Recall that $\OL$ is \emph{big} on $Y$ if $\wh\vol(\OL)>0$ (cf. \cite[\S5.2]{YZ2}). If $\OL$ is nef, then we have the adelic Hilbert--Samuel formula 
$\wh\vol(\OL)=\OL^{\dim Y}$ (cf. \cite[Thm. 5.2.2]{YZ2}).
In this case, $\OL$ is big if and only if $\OL^{\dim Y}>0$.

We say that $\OL$ is \emph{potentially big on $Y/S$} if there is a positive integer $m$ such that the adelic line bundle 
$$
\OL^{\boxtimes m}:=(p_1^* \OL) \otimes (p_2^* \OL) \otimes \cdots \otimes (p_m^* \OL)
$$
is big on the $m$-fold fiber product
$$
Y_{/S}^m:=Y\times_SY\times_S\cdots \times_S Y.
$$
In additive notations, the above is written as 
$$
m_\boxtimes\OL=p_1^* \OL + p_2^* \OL + \cdots + p_m^* \OL.
$$

We remark that $Y_{/S}^m$ is always integral under the assumptions that $Y$ and $S$ are integral, and that $\pi:Y\to S$ is projective and flat with a geometrically integral generic fiber.
In fact, the generic fiber of $Y_{/S}^m\to S$ is geometrically integral. 
By flatness, $Y_{/S}^m\to S$ is equi-dimensional. Then $Y_{/S}^m$ is irreducible since its generic fiber is irreducible. 
To prove that $Y_{/S}^m$ is reduced, note that its generic fiber over $S$ is reduced, it suffices to check that it has no embedded components. 
Then it suffices to check that $Y_{/S}^m$ satisfies property $S_1$ (cf. \cite[\S8.2, 2.19, 2.20]{Liu}).
This follows from flatness. 
In fact, for any $y\in Y_{/S}^m$ with image $s\in S$, if $s$ is not the generic point of $S$, then the flatness implies that $\mathrm{depth}_y(Y_{/S}^m)\geq \mathrm{depth}_s(S)\geq 1$. 
If $s$ is the generic point of $S$, $y$ lies in the generic fiber of 
$Y_{/S}^m\to S$, which is integral and thus satisfies $S_1$.

\subsubsection{A quick criterion}

Concerning potential bigness, we have the following key criterion. 

\begin{thm}[potential bigness] \label{potentially big}
\kkk
Let $S$ be a quasi-projective and flat normal integral scheme over $k$. 
Let $Y$ be an integral scheme with a projective and flat morphism $\pi:Y\to S$ over 
$k$. 
Assume that the generic fiber of $\pi:Y\to S$ is geometrically integral.
Let $\OL$ be an adelic line bundle on $Y/k$.
Assume that $\OL$ is nef on $Y$, and the degree of $L$ on fibers of $\pi:Y\to S$ is strictly positive.
Then the following are equivalent:
\begin{enumerate}[(1)]
\item $\pi_*\pair{\OL,\cdots,\OL}$ is big on $S$;
\item $m_\boxtimes\OL$ is big on $Y_{/S}^m$ for all $m\geq \dim S$;
\item $\OL$ is potentially big on $Y/S$; i.e., $m_\boxtimes\OL$ is big on $Y_{/S}^m$ for some $m\geq 1$.
\end{enumerate}
\end{thm}

\begin{proof}
Denote $d=\dim S$ and $e=\dim Y-\dim S$, so $\dim Y_{/S}^m=d+em$.

We first prove that (1) implies (2). 
Since $\OL$ is nef, it is immediate that $m_\boxtimes\OL$ is nef. 
Its bigness is equivalent to the positivity of the intersection number
$$
(m_\boxtimes\OL)^{d+em}=(p_1^* \OL + p_2^* \OL + \cdots + p_m^* \OL)^{d+em}.
$$
Expanding the right-hand side, we obtain a term
$$
b_m=(p_1^* \OL)^{e+1}  \cdots (p_d^* \OL)^{e+1}
\cdot (p_{d+1}^* \OL)^e  \cdots (p_m^* \OL)^e.
$$
We claim that 
$b_m=a^{m-d} (\OM^d),$
where $a$ denotes the degree of $L$ on fibers of $\pi:Y\to S$,
and $\OM^d$ denotes the top self-intersection number of 
$\OM=\pi_*\pair{\OL,\cdots,\OL}$ on $S$.
The claim implies $b_m>0$ and thus proves that (1) implies (2). 

To prove the claim, we first treat the arithmetic case $k=\ZZ$.
By a limit process, we can assume that $S$ is projective and flat over $S_0=\Spec\ZZ$, and $Y\to S$ is projective and flat.
We can further assume that $\OL$ is a nef hermitian line bundle on $Y$, so that 
$$\OM=\pi_*\pair{\OL,\cdots,\OL}$$ 
is a nef hermitian line bundle on $S$.
Consider the projection $\psi_m:Y_{/S}^m\to S_0$. 
By the compatibility of the Deligne pairing with the intersection number, $b_m$ is the arithmetic degree of  
$$
\OB_m= \psi_{m*} \pair{(p_1^* \OL)^{e+1},\,   \cdots,\, (p_d^*\OL)^{e+1},\, 
(p_{d+1}^* \OL)^e,\,  \cdots,\, (p_m^* \OL)^e}.
$$
Here the notation $(p_i^* \OL)^{t}$ means that $p_i^* \OL$ appears $t$ times in the pairing. 

If $Y_{/S}^{m-1}$ is normal, 
apply \cite[Lem. 4.6.1(1)]{YZ2} to the composition
$$
Y_{/S}^m\stackrel{q}{\lra} Y_{/S}^{m-1} \stackrel{\psi_{m-1}}{\lra}  S_0.
$$
Here $q$ is the projection by forgetting the first component of $Y_{/S}^m$.
In the line bundles defining $\OB_m$, view the first $e+1$ terms, indicated by $(p_1^* \OL)^{e+1}$, as line bundles on $Y_{/S}^m$, and view the remaining $(d+me)-(e+1)$ terms as the pull-back of line bundles via $q$.
Then the lemma implies that
$$
\OB_m= \psi_{m-1,*} \pair{\OM_{Y_{/S}^{m-1}},\, (p_2^* \OL)^{e+1},\,   \cdots,\, (p_d^*\OL)^{e+1},\, 
(p_{d+1}^* \OL)^e,\,  \cdots,\, (p_m^* \OL)^e}.
$$
Here by abuse of notations, $p_2,\cdots, p_m:Y_{/S}^{m-1}\to Y$ denotes the projections to the $m-1$ components.

Here the normality of $Y_{/S}^{m-1}$ is required for the definition of the Deligne pairing via $q:Y_{/S}^m\to Y_{/S}^{m-1}$. However, if $Y_{/S}^{m-1}$ is not normal, denote by $W_{m-1}$ its normalization, and apply \cite[Lem. 4.6.1(1)]{YZ2} to the composition
$$
Y\times_S W_{m-1}  \stackrel{q'}{\lra} W_{m-1} \stackrel{\psi_{m-1}'}{\lra}  S_0
$$
induced by 
$$
Y_{/S}^m\stackrel{q}{\lra} Y_{/S}^{m-1} \stackrel{\psi_{m-1}}{\lra}  S_0.
$$
Then we obtain a similar formula for $\OB_m$, and this formula implies
 $$
b_m=\wh\deg\, \psi_{m-1,*} \pair{\OM_{Y_{/S}^{m-1}},\, (p_2^* \OL)^{e+1},\,   \cdots,\, (p_d^*\OL)^{e+1},\, 
(p_{d+1}^* \OL)^e,\,  \cdots,\, (p_m^* \OL)^e}
$$
by \cite[Prop. 4.1.2]{YZ2}.

Repeat the process on $(p_2^* \OL)^{e+1},\,   \cdots,\, (p_d^*\OL)^{e+1}$.
We obtain  
$$
b_m=\wh\deg\, \psi_{m-d,*} \pair{(\OM_{Y_{/S}^{m-d}})^d, \, 
(p_{d+1}^* \OL)^e,\,  \cdots,\, (p_m^* \OL)^e}.
$$
Apply \cite[Lem. 4.6.1(3)]{YZ2} to the composition
$$
Y_{/S}^{m-d}\stackrel{}{\lra} S \stackrel{\psi_0}{\lra}  S_0.
$$
Here we view the first $d$ terms as the pull-back of line bundles from $S$.
Then the lemma implies that
$$
b_m=\wh\deg \big(a^{m-d}\,\psi_{0*} \pair{\OM^d}\big).
$$
This proves that (1) implies (2) in the arithmetic case.

For the geometric case that $k$ is a field, we can first reduce it to the case that $S$ is projective over $k$. By blowing-up $S$ if necessary, we can find a fibration $S\to S_0$ with $S_0=\BP^1_k$. 
Then the proof is similar to the arithmetic case by taking Deligne pairings to $S_0$. 

Note that (2) implies (3) trivially. 
It remains to prove that (3) implies (1).
Assume that for some $m\geq 1$, $m_\boxtimes\OL$ is big on $Y_{/S}^m$.
By Lemma \ref{basic8}(1), the Deligne pairing 
$$\ON=(\pi_m)_*\pair{(m_\boxtimes\OL)^{em+1}}$$ 
is nef and big on $S$.
Here $\pi_m:Y_{/S}^m\to S$ is the structure morphism. 
It suffices to prove that $\ON$ is a positive multiple of $\OM=\pi_*\pair{\OL^{e+1}}$.

In fact, consider the expansion of 
$$\ON=(\pi_m)_*\pair{(p_1^* \OL + p_2^* \OL + \cdots + p_m^* \OL)^{em+1}}$$ 
by linearity.
It is a positive linear combination of
$$\ON(r_1,\cdots, r_m)=(\pi_m)_*\pair{(p_1^* \OL)^{r_1}, (p_2^* \OL)^{r_2}, \cdots, (p_m^* \OL)^{r_m}}$$ 
with $r_1+\cdots+r_m=em+1$. 
We claim that $\ON(r_1,\cdots, r_m)\neq 0$ only if $(r_1,\cdots, r_m)$ is a permutation of $(e+1,e,\cdots, e)$; in that case, $\ON(r_1,\cdots, r_m)= a^{m-1}\OM$.
Here $a$ denotes the degree of $L$ on fibers of $\pi:Y\to S$ as above.

In fact, by $r_1+\cdots+r_m=em+1$, some term $r_i\geq e+1$. 
By symmetry, we assume that $r_1\geq e+1$. 
Write 
$$\ON(r_1,\cdots, r_m)=(\pi_m)_*\pair{(p_1^* \OL)^{e+1}, (p_1^* \OL)^{r_1-e-1}, (p_2^* \OL)^{r_2}, \cdots, (p_m^* \OL)^{r_m}}$$ 
If $Y$ is normal, apply \cite[Lem. 4.6.1(3)]{YZ2} to the composition
$$
Y_{/S}^m\stackrel{p_1}{\lra} Y \stackrel{\pi}{\lra}  S.
$$
We see that 
$$\ON(r_1,\cdots, r_m)=c(r_1,\cdots, r_m) \pi_*\pair{ \OL^{e+1}},$$ 
where $c(r_1,\cdots, r_m)$ is the intersection number of 
$$(p_1^* L)^{r_1-e-1}, (p_2^* L)^{r_2}, \cdots, (p_m^* L)^{r_m}$$ 
over a fiber of $p_1:Y_{/S}^{m}\to Y$. 
Then $c(r_1,\cdots, r_m)\neq 0$ only if $r_1-e-1=0$, and $r_i=e$ for $i>1$;
in that case $c(r_1,\cdots, r_m)=a^{m-1}$. 

If $Y$ is not normal, we can apply a similar normalization trick as above to get the same result.
This finishes the proof.
\end{proof}

\subsubsection{Consequence on small points}

\kkk
By a \emph{global field over} $k$, we mean a field $K$ as follows. 
If $k=\ZZ$, then $K$ is a number field;
if $k$ is a field, then $K$ is a function field of one variable over $k$.

In the terminology of \cite[\S2.3]{YZ2}, a quasi-projective variety $V$ over $K$ is essentially quasi-projective over $k$. In \cite[\S2.4,\,\S2.5]{YZ2}, we can talk about adelic divisors and adelic line bundles on $V$ with base $k$, which are defined by direct limits over 
quasi-projective models of $V$ over $k$.

Let $V$ be a quasi-projective variety over $K$. Let $\OL$ be an adelic line bundle on $V$ over $k$. 
We define a height function 
$$h_\OL:V(\ol K)\lra \RR$$
 as follows.
For any $x\in V(\OK)$, denote by $\tilde x$ the closed point of $V$ corresponding to $x$.
By the construction in \cite[\S2.5.5]{YZ2}, the pull-back $\OL|_{\tilde x}$ gives an adelic line bundle on $\tilde x/k$. 
Since $\tilde x$ is the spectrum of a global field over $k$, we can take the arithmetic degree $\wh\deg(\OL|_{\tilde x})$. 
Here if $k$ is a field, the degree is normalized by multiplicity functions given by degrees over $k$.

If $k=\ZZ$, we define 
$$
h_\OL(x)=\frac{1}{\deg(\tilde x/\QQ)}\wh\deg(\OL|_{\tilde x}).
$$
Note that the normalizing factor has the degree of $\tilde x$ over $\QQ$ (instead of $K$), which is different from that in \cite[\S5.3]{YZ2}. 

If $k$ is a field, we define 
$$
h_\OL(x)=\frac{1}{\deg(\tilde x/K)}\wh\deg(\OL|_{\tilde x}).
$$
The normalizing factor depends on $K$ and agrees with that in \cite[\S5.3]{YZ2}.

The following is our key result relating potential bigness to distribution of small points. 
Due to different normalization of heights and arithmetic degrees, we state the result for the arithmetic case and for the geometric case separately.

\begin{thm} \label{small points}
The following statements hold.
\begin{enumerate}[(1)]
\item
Let $S$ be a quasi-projective variety over a number field $K$.
Let $\pi:X\to S$ be a smooth relative curve. 
Let $\OL$ be a nef adelic line bundle on $X/\ZZ$, and $\OM$ be an adelic line bundle on $S/\ZZ$.
If $\OL$ is potentially big on $X/S$, then there are a non-empty Zariski open subvariety $U\subset S$, and constants  
$c_1,c_2>0$, such that 
for any $y\in U(\ol K)$,
$$
\#\{x\in X(\ol K): \pi(x)=y, \, h_\OL(x)\leq c_1\, h_\OM(y)\} \leq c_2.
$$ 

\item 
Let $S$ be a quasi-projective variety over a field $k$ with $\dim S>0$.
Let $\pi:X\to S$ be a smooth relative curve. 
Let $\OL$ be a nef adelic line bundle on $X/k$, and $\OM$ be an adelic line bundle on $S/k$.
If $\OL$ is potentially big on $X/S$, then there are a non-empty Zariski open subvariety $U\subset S$, and constants  
$c_1,c_2>0$, such that 
for any 1-dimensional point $y\in U$,
$$
\#\{x\in X_y(\ol{k(y)}):   h_\OL(x)\leq c_1\,\wh\deg(\OM|_y)\} \leq c_2.
$$ 
Here $h_\OL:X_y(\ol{k(y)})\to \RR$ is the height function associated to $\OL|_{X_y}$ and normalized by degrees over $k(y)$, and $\wh\deg: \wh\Pic(k(y)/k)\to \RR$ is normalized by multiplicity functions given by degrees over $k$.
\end{enumerate}
\end{thm}

\begin{proof}
We first prove (1). 
By assumption, $\OL_m=m_\boxtimes\OL$ is big on $X_m=X_{/S}^m$ for some $m\geq 1$. Fix such an $m$. 
By \cite[Thm. 5.3.5(1)]{YZ2}, there are a closed subset $Z$ of codimension one in $X_m$ and $\epsilon>0$ such that 
$$
\{x\in X_m(\ol K):h_{\OL_m}(x)\leq \epsilon h_\OM(\pi_m(x))  \} \subset Z(\ol K).
$$

For any $y\in S(\ol K)$, denote 
$$
\Sigma(y):=\{x \in X(\ol K): \pi(x)=y,\, h_{\OL}(x)\leq \frac{\epsilon}{m} h_\OM(y)\}.  
$$
The key is the inclusion 
$$\Sigma(y)^m \subset Z(\ol K).$$
This follows from the height identity
$$
h_{\OL_m}(x)=h_{\OL}(x_1)+\cdots +h_{\OL}(x_m)
$$
for any $x\in X_m(\ol K)$ represented by $(x_1,\cdots, x_m)$ under the expression 
$$X_m(\ol K)=\{(x_1,\cdots, x_m)\in X(\ol K)^m: \pi(x_1)=\cdots=\pi(x_m)\}.$$ 

Let $U$ be a non-empty open subscheme of $S$ such that $Z$ is flat over $U$. 
By Lemma \ref{counting} below, the inclusion $\Sigma(y)^m \subset Z_y(\ol K)$ 
forces $\Sigma(y)$ to be finite with
$$
\#\Sigma(y)\leq  \deg_{m_\boxtimes N}(Z_\eta),\quad
\forall y\in U(\overline K).
$$
Here $Z_\eta$ is the generic fiber of $Z\to S$, 
 $N$ is a line bundle of degree 1 on the generic fiber $X_\eta$ of $X\to S$, 
so that $m_\boxtimes N$ is a line bundle on $X_\eta^m$. 
This proves (1).

The proof of (2) is similar to that of (1), so we will only emphasize the difference.
We still have that $\OL_m=m_\boxtimes\OL$ is big on $X_m=X_{/S}^m$ for some 
$m\geq 1$. 
We are not able to apply \cite[Thm. 5.3.5(1)]{YZ2} directly, but we can modify its proof slightly. 
In fact, as in the proof of the loc. cit., $\OL_m-\epsilon \pi^*\OM$ is big for some positive rational number $\epsilon>0$.
Then some multiple $n(\OL_m-\epsilon\OM)$ of $\OL_m-\epsilon\OM$ is an (integral) adelic line bundle with a nonzero effective section $s$. 

Denote $Z=\div(s)$. 
We claim that for any 1-dimensional point $y$ of $S$,
$$
\{x\in (X_m)_y(\overline{k(y)}):
h_{\OL_m}(x)\leq \epsilon \,\wh\deg(\OM|_y)  \} \subset Z_y(\overline{k(y)}).
$$ 
In fact, denote by $\tilde x$ the closed point of $(X_m)_y$ corresponding to $x$, which is a 1-dimensional point of $X_m$.
If $x\notin Z_y(\OK)$ or equivalently
$\tilde x\notin Z_y$, then the restriction of $s$ gives a nonzero effective section of 
$$
n(\OL_m-\epsilon \pi^*\OM)|_{\tilde x}
=n(\OL_m|_{\tilde x}-\epsilon \pi_x^*(\OM|_{y})),
$$ 
where $\pi_x:\tilde x\to y$ is the induced morphism, which is finite. 
Taking arithmetic degrees, we have
$$
0\leq \wh\deg(\OL_m|_{\tilde x})-\epsilon \wh\deg(\pi_x^*(\OM|_{y}))
=\deg(\pi_x)\big(h_{\OL_m}( x)-\epsilon \wh\deg(\OM|_{y})\big).
$$
This prove the claim. 
The remaining part of the proof is similar to that of (1).
\end{proof}

In the above proof, we have used the following basic result, which improves \cite[Lem. 6.3]{DGH1} and \cite[Lem. 7.3]{Gao3}.

\begin{lem} \label{counting}
Let $C$ be a smooth projective curve over an algebraically closed field $F$. 
Let $Z\subsetneq C^m$ be a Zariski closed subset of codimension $d>0$.
Let $\Sigma$ be a finite subset of $C(F)$ such that $\Sigma^m\subset Z(F)$.
Then 
$$
\#\Sigma\leq  (\deg_{m_\boxtimes N}(Z))^{1/d}. 
$$
Here $N$ is a line bundle on $C$ of degree 1, and
$$
\deg_{m_\boxtimes N}(Z)=\sum_{i=1}^r \deg_{m_\boxtimes N}(Z_i), 
$$
where $Z_1,\cdots, Z_r$ are all the irreducible components of $Z$. 
\end{lem}

\begin{proof}
Denote $M=m_\boxtimes N$ and $n=\#\Sigma$. 
We claim that for every $i=1,\cdots, m$
$$
\deg_{nM}(Z)\geq  \deg_{nM}(Z\cap p_i^{-1}\Sigma).
$$
Here $p_i:C^m\to C$ is the projection to the $i$-th component, and
$Z\cap p_i^{-1}\Sigma$ is endowed with the reduced structure.

If the claim holds, then successively applying the inequality gives 
$$
\deg_{nM}(Z)\geq  \deg_{nM}(Z\cap p_1^{-1}\Sigma\cap 
\cdots\cap p_m^{-1}\Sigma)
= \#(\Sigma^m)
=n^m.
$$
This implies the result since
$$
\deg_{nM}(Z)\leq n^{m-d}\deg_{M}(Z)
$$
by $\dim Z_i\leq m-d$. 

It remains to prove the claim.
By writing in terms of irreducible components, we can assume that $Z$ is irreducible. 
If $Z$ is contained in $p_i^{-1}\Sigma$, the degree does not change. 
If $Z$ is not contained in $p_i^{-1}\Sigma$, then it intersects 
$p_i^{-1}\Sigma$ properly, and thus 
$$
\deg_{nM}(Z)
=(nM)^{\dim Z}\cdot Z
\geq (nM)^{\dim Z-1}\cdot (np_i^*N) \cdot Z
\geq \deg_{nM}(Z\cap p_i^{-1}\Sigma).
$$
This proves the claim. 
\end{proof}

\begin{remark}
Here we describe a path to find effective constants $(c_1,c_2)$ in Theorem \ref{small points}, which depends on many numerical invariants of $(X/S, \OL,\OM)$. We focus on the arithmetic case.
Instead of applying \cite[Thm. 5.3.5(1)]{YZ2}, we take the idea sketched right after Theorem \ref{bigness77}.
First, by \cite[Lem. 5.1.6]{YZ2}, we can find a nef adelic line bundle $\OM'$ on $S$ such that $\OM'-\OM$ is effective. 
Then we can assume that $\OM$ is nef by 
replacing $\OM$ by $\OM'$ everywhere in the problem.
Second, we have
$$
\wh\vol(\OL_m-\epsilon\OM)\geq \OL_m^{d+m}-(d+m)\OL_m^{d+m-1}\cdot \epsilon\OM>0.
$$ 
Here $d$ is the dimension of a quasi-projective model of $S$ over $k$,  $m$ is any fixed integer with $m\geq d$, and $\epsilon$ is chosen to be any positive rational number satisfying the second inequality. 
Then the constant $c_1=\epsilon/m$. 
Third, some multiple $n(\OL_m-\epsilon\OM)$ of $\OL_m-\epsilon\OM$ has a nonzero effective section $s$. 
By Lemma \ref{counting}, we can express $c_2$ in terms of intersection numbers involving $Z=\div(s)$. Note that $Z$ is linearly equivalent to $n(\OL_m-\epsilon\OM)$.  
To find an effective $c_2$, it suffices to find an effective positive integer $n$ such that $n(\OL_m-\epsilon\OM)$ is an (integral) adelic line bundle with a nonzero effective section. 
\end{remark}

\subsection{Three potentially big examples} \label{sec big examples}

The goal of this subsection is to introduce three potentially big examples built on fibered powers of families of curves.

\kkk
Let $S$ be a flat and quasi-projective normal integral scheme over $k$.
Let $\pi:X\to S$ be a smooth relative curve over $S$ of genus $g>1$.
Let $\pi_J:J\to S$ be the relative Jacobian scheme of $X$ over $S$.

Here we review three canonical morphisms from $X$ to $J$ (up to base changes) treated in Theorem \ref{isomorphism5} and Theorem \ref{isomorphism3}. 
Only the third morphism will be crucially used in our treatments later, but we include the other two for their own value. For example, the second morphism is essentially the one used in the proofs of \cite{DGH1, Kuh}. 

The first morphism is the basic one. 
Let $\alpha$ be a line bundle on $X$ with degree $d> 0$ on the fibers of $X\to S$. 
This gives a finite $S$-morphism 
$$
i_\alpha: X\lra J, \quad x\longmapsto dx-\alpha.
$$

The second one is the morphism
$$
i_\Delta: X\times_S X\lra  X\times_SJ ,\quad (x,y)\longmapsto (x,y-x).
$$
This agrees with the $X$-morphism 
$$
i_\Delta: X_X\lra  J_X,\quad x\longmapsto x-\Delta
$$
defined in the proof of Theorem \ref{isomorphism3}(2). 
Namely, denote $X_X=X\times_S X$ and $J_X=X\times_S J$, viewed as $X$-schemes via the first projections $p_1:X_X\to X$ and $q_1:J_X\to X$. 
Then $q_1:J_X\to X$ is canonically isomorphic to the Jacobian scheme of $p_1:X_X\to X$.
View $\Delta:X\to X\times_SX=X_X$ as a section of   
$p_1:X_X\to X$. This defines an $X$-morphism 
$$
i_\Delta=i_{\CO(\Delta)}: X_X\lra  J_X,\quad x\longmapsto x-\Delta.
$$

The third one is the morphism 
$$\tau:J\times_SX \lra J\times_SJ, \quad (y,x)\longmapsto (y,y+(2g-2)x-\omega_{X/S}).$$
As in the proof of Theorem \ref{isomorphism3}, if $X\to S$ has a section, then $J\times_SX$ has a universal line bundle $Q$, and $\tau$ agrees with
$$
i_{\omega-Q}: X_J\lra J_J,\quad x\longmapsto (2g-2)x-(\omega_{X_J/J}-Q).
$$

Recall that from Definition \ref{qqq} we have $\Theta=\Delta_J^*(P^\vee)$ on $J$. 
We also have a line bundle $\Theta_X=p_2^*\Theta$ on $J_X$.
Then we have the canonical adelic extensions $\OTheta$ on $J$ and its base changes 
$\OTheta_X$ to $J_X$ and $\OTheta_J$ to $J_J$.
Finally, we have the following bigness result. 

\begin{thm}  \label{bigness7}
\kkk
Let $S$ be a flat and quasi-projective normal integral scheme over $k$.
Let $\pi:X\to S$ be a smooth relative curve over $S$ of genus $g>1$ with maximal variation.
Then the following hold:
\begin{enumerate}[(1)]
\item  The adelic line bundle $\pi_*\pair{i_\alpha^*\OTheta,i_\alpha^*\OTheta}$ is nef and big on $S$. Therefore, $i_\alpha^*\OTheta$ is potentially big on $X/S$. 
\item The adelic line bundle $p_{1*}\pair{i_\Delta^*(\OTheta_X), i_\Delta^*(\OTheta_X)}$ is nef and big on $X$. Therefore, $i_\Delta^*(\OTheta_X)$ is potentially big on $X_X/X$. 
\item  The adelic line bundle $q_{1*}\pair{\tau^*(\OTheta_J), \tau^*(\OTheta_J)}$ is nef and big on $J$. Therefore, $\tau^*(\OTheta_J)$ is potentially big on $X_J/J$. 
\item  Let $Y$ be a quasi-projective variety over $k$ with a generically finite morphism $Y\to J$.
Denote by $T$ the image of $Y\to S$. 
Assume that the composition $T\to S\to M_{g,k}$ is generically finite.
Then the adelic line bundle $q_{1*}\pair{\tau^*(\OTheta_J), \tau^*(\OTheta_J)}|_{Y}$ is nef and big on $Y$. Therefore, $(\tau^*(\OTheta_J))|_{Y\times_SX}$ is potentially big on $Y\times_SX/Y$. 
\end{enumerate}

\end{thm}

\begin{proof}

By Theorem \ref{potentially big},  potential bigness follows from  bigness of the Deligne pairing. Thus in each part, we only need to prove the first statement.

For (1), by Theorem \ref{isomorphism5}(2), 
$$
\pi_*\pair{i_\alpha^*\OTheta,i_\alpha^*\OTheta}
= \frac{gd^4}{g-1}\pi_*\pair{\ol\omega_{X/S,a}\, \ol\omega_{X/S,a}}
+\frac{d^2}{g-1}\iota_\alpha^*\OTheta
$$
in $\wh\Pic(S)_\QQ$.
By Theorem \ref{bigness1}, $\pi_*\pair{\overline\omega_{X/S,a},\overline\omega_{X/S,a}}$ is nef and big on $S$. 
Note that $\iota_\alpha^*\OTheta$ is nef on $S$, since 
$\OTheta$ is nef on $J$.
It follows that $\pi_*\pair{i_\alpha^*\OTheta,i_\alpha^*\OTheta}$ is nef and big on $S$.

For (2), by Theorem \ref{isomorphism3}(2), 
$$
p_{1*}\pair{i_\Delta^*(\OTheta_X),i_\Delta^*(\OTheta_X)}=
4g\, \overline\omega_{X/S,a}+\pi^*\pi_*\pair{\overline\omega_{X/S,a},\overline\omega_{X/S,a}}
$$
in $\wh\Pic(X)_\QQ$.
By Theorem \ref{bigness5}, $\overline\omega_{X/S,a}$ is nef and big on $X$.
Then $p_{1*}\pair{i_\Delta^*(\OTheta_X),i_\Delta^*(\OTheta_X)}$ is nef and big on $X$. 

For (3), by Theorem \ref{isomorphism3}(3),   
$$
q_{1*}\pair{\tau^*(\OTheta_J), \tau^*(\OTheta_J)}= 
16(g-1)^3\OTheta+16g(g-1)^3\pi_J^*\pi_*\pair{\ol\omega_{X/S,a}\, ,\ol\omega_{X/S,a}}
$$
in $\wh\Pic(J)_\QQ$.
Then $q_{1*}\pair{\tau^*(\OTheta_J), \tau^*(\OTheta_J)}$ is nef on $J$.
For the bigness, it suffices to note that the binomial expansion of the self-intersection
of $q_{1*}\pair{\tau^*(\OTheta_J), \tau^*(\OTheta_J)}$ on $J$ has a term
$$
(\pi_J^*\pi_*\pair{\ol\omega_{X/S,a}\, ,\ol\omega_{X/S,a}})^{\dim S}\cdot 
\OTheta^{\dim J-\dim S}
=n(\pi_*\pair{\ol\omega_{X/S,a}\, ,\ol\omega_{X/S,a}})^{\dim S}
>0.
$$
Here $n$ is the degree of $\Theta$ on a fiber of $J\to S$, and the identity is the limit version of that on projective varieties over $k$.

For (4), replacing $Y$ by the Zariski closure of the image of $Y$ in $J$, we can assume that $Y$ is a closed subvariety of $J$ and thus projective over $S$.
This process is based on the fact that bigness does not change under generically finite and surjective base changes, which was also mentioned in the proof of
Theorem \ref{bigness2}.
Then we need to check that
$$
q_{1*}\pair{\tau^*(\OTheta_J), \tau^*(\OTheta_J)}|_Y= 
16(g-1)^3\OTheta|_Y+16g(g-1)^3\pi_J^*\pi_*\pair{\ol\omega_{X/S,a}\, ,\ol\omega_{X/S,a}}|_Y
$$
is big in $\wh\Pic(Y)_\QQ$.
The proof is similar to (3), as $\Theta|_Y$ is ample on the generic fiber of $Y\to S$.
\end{proof}

\begin{remark}
Only (3) and (4) will be used in the proof of our uniform Bogomolov-type theorem.
Note that (4) essentially implies (1)-(3), but we include all of them here for their own simplicity.
There are also counterparts of (4) for (1) and (2), but we omit them in this paper. 
\end{remark}

\subsection{The uniform Bogomolov-type theorem: statement} \label{sec small points}

Now we state the following uniform version of the Bogomolov conjecture, which generalizes the new gap principle proved in \cite[Prop. 7.1]{DGH1} and \cite[Thm. 3]{Kuh} and summarized in \cite[Thm. 4.1]{Gao3}.

\begin{thm}[Theorem \ref{small points22}] \label{small points2}
Let $g>1$ be an integer. 
Then there are constants $c_1,c_2>0$ depending only on $g$ satisfying the following properties. 
Let $K$ be either a number field or a function field of one variable over a field $k$. 
Then for any geometrically integral, smooth and projective curve $C$ of genus $g$ over $K$, and for any 
line bundle $\alpha\in \Pic(C_{\ol K})$ of degree 1, with the extra assumption that $(C_{\ol K},\alpha)$ is non-isotrivial over $\overline k$ in the case that $K$ is a function field of one variable over a field $k$, 
one has
\small
$$
\#\left\{x\in C(\ol K): \wh h(x-\alpha)\leq c_1\big( \max\{h_\Fal(C),1\}+
\wh h((2g-2)\alpha-\omega_{C/K})\big) \right\} \leq c_2.
$$ 
\normalsize
\end{thm}

The theorem will be proved in \S\ref{sec small points2}. In the following, we explain the precise definition of the heights in the theorem, and explore the non-isotriviality condition in terms of moduli spaces.

\subsubsection{Normalization of the heights} \label{sec normalization}

Since Theorem \ref{small points2} asserts uniformity when varying $K$, it is necessary to explain how the heights are normalized in terms of $K$. 
The normalization here is compatible with that in Theorem \ref{small points}. 

Let $K$ be either a number field or a function field of one variable over a field $k$. 
We will refer these two cases as ``the arithmetic case'' and ``the geometric case''. 

Let us first explain the N\'eron--Tate height.
Let $C$ be a geometrically integral, smooth and projective curve of genus $g>1$ over $K$. 
Take a point $x_0\in C(\OK)$. 
Denote by $J$ the Jacobian variety of $C$ over $K$.
The theta divisor $\theta_{x_0}$ is defined as the image of the morphism 
$$C_{\OK}^{g-1}\lra J_\OK,\quad (x_1,\cdots, x_{g-1})\longmapsto 
x_1+\cdots+ x_{g-1}-(g-1)x_0.$$
Then $\Theta=\CO(\theta_{x_0}+[-1]^*\theta_{x_0})$ is a symmetric and ample line bundle on $J_\OK$. 
By Theorem \ref{isomorphism}(4), the isomorphism class of $\Theta$ is independent of the choice of $x_0$ and actually descends to a line bundle on $J$. 
Moreover, the current definition of $\Theta$ is also compatible with that of Definition \ref{qqq}. 

The N\'eron--Tate height $\wh h:J(\ol K)\to \RR$
is defined by $\ds\wh h=\frac12 \wh h_{\Theta}$;
the canonical height 
$\wh h_\Theta:J(\ol K)\to \RR$ is defined by Tate's limiting argument from a Weil height of $\Theta$; the Weil height is induced by the classical height 
$h:\BP^n(\ol K) \to \RR$ as follows.

If $K$ is a number field, the classical height 
$h:\BP^n(\ol K) \to \RR$ is
given by 
$$
h(x_0,\cdots, x_n)=\frac{1}{[K':\QQ]} \sum_{v\in M_{K'}} \log \max\{\|x_0\|_v,\cdots, \|x_n\|_v\}.
$$
Here $K'$ is a finite extension of $K$ containing the coordinates $x_0,\cdots, x_n$,  $M_{K'}$ is the set of places of $K'$, and $\|\cdot\|_v$ is normalized as 
follows:
\begin{enumerate}[(1)]
\item If $v$ is real, it is the usual absolute value on $\RR$.
\item If $v$ is complex, it is the square of the usual absolute value on $\CC$. 
\item If $v$ is non-archimedean, then $\|\cdot\|_v=N_v^{-\ord_v(\cdot)}$.
Here $N_v$ is the cardinality of the residue field of $v$ in $O_{K'}$. 
\end{enumerate}

If $K$ is a function field of one variable over a field $k$, the classical height 
$h:\BP^n(\ol K) \to \RR$ is
given by 
$$
h(x_0,\cdots, x_n)=\frac{1}{[K':K]} \sum_{v\in M_{K'}} \log \max\{\|x_0\|_v,\cdots, \|x_n\|_v\}.
$$
Here $K'$ is a finite extension of $K$ containing the coordinates $x_0,\cdots, x_n$,  $M_{K'}$ is the set of places of $K'$, and $\|\cdot\|_v$ is normalized 
by $\|\cdot\|_v=e^{-[k(v):k]\,\ord_v(\cdot)}$.
Here $k(v)$ is the residue field of the place $v$ of $K'$. Note that even if $K'$ (or even $K$) contains a bigger constant field $k'$, the above factor is still $[k(v):k]$ instead of 
$[k(v):k']$.

Now we explain the Faltings height, which is actually the stable Faltings height. 
Let $C$ and $J$ be as above.
Let $K'$ be a finite extension of $K$ such that $C$ (or equivalently $J$) has \emph{semistable reduction} everywhere over $K'$. 
Denote by $C'=C_{K'}$ and $J'=J_{K'}$ the base change to $K'$.

In the number field case, denote $B'=\Spec O_{K'}$; 
in the function field case, denote by $B'$ the unique projective and regular curve over $k$ with function field $K'$. 
Denote by $\CJ'$ the N\'eron model of $J'$ over $B'$, and denote by $e:B'\to \CJ'$ the  identity section.
Denote the Hodge bundle 
$$\underline\omega_{B'}=e^*\Omega_{\CJ'/B'}^g,$$ 
which is a line bundle on $B'$. 

In the function field case, define the Faltings height
$$
h_\Fal(C)=h_\Fal(J)= \frac{1}{[K':K]} \deg_{B'/k}(\underline\omega_{B'}).
$$
Here the degree $\deg_{B'/k}(\underline\omega_{B'})$ comes from the map $\deg_{B'/k}:\Div(C')\to \ZZ$ given by $\deg_{B'/k}(v)=[k(v):k]$, where $v\in B'$ is a closed point and $k(v)$ is the residue field of $v$ in $B'$.

Via the stable relative curve $\CCC'\to B'$ extending $C_{K'}$, we obtain a Hodge bundle $\lambda_{B'}$ over $B'$ defined in \S\ref{sec bundle on stable}. 
We claim that
$$
h_\Fal(C)= \frac{1}{[K':K]} \deg_{B'/k}(\lambda_{B'}).
$$
This follows from a special case of Lemma \ref{hodge bundle equivalence}. 
In fact, by \cite[\S7.4, Prop. 4]{BLR}, the group scheme $\Pic_{\CCC'/B'}^0$  is the relative identity component of the N\'eron model of $J_{K'}$ over $B'$.
Then Lemma \ref{hodge bundle equivalence} gives a canonical isomorphism $\lambda_{B'}\to \underline\omega_{B'}$.

In the number field case, endow $\underline\omega_{B'}$ with the Faltings metric such that 
$$
\|\alpha\|_{\rm Fal}^2=\frac{i^{g^2}}{2^g} \int_{J'_\sigma(\CC)} \alpha\wedge\bar \alpha
$$
for any embedding $\sigma:K'\to \CC$ and any element $\alpha$ of
$$
(\underline\omega_{B'})_\sigma(\CC)
=e^*\Omega_{J'_\sigma(\CC)/ \CC}^g
\simeq
\Gamma(J'_\sigma(\CC), \Omega_{J'_\sigma(\CC)/ \CC}^g).
$$
Define the Faltings height
$$
h_{\Fal}(C)=h_\Fal(J)=\frac{1}{[K':\QQ]} \wh\deg(\underline\omega_{B'}, \|\cdot\|_{\rm Fal}).
$$
Here the degree $\wh\deg(\underline\omega_{B'}, \|\cdot\|)$ comes from the map 
$\wh\deg:\wh\Div(O_{K'})\to \RR$ given by $\wh\deg(v)=\log N_v$, where $N_v$ is
the order of the residue field of $v$.

\subsubsection{The non-isotriviality condition} \label{sec isotrivial}

Now let us explain the non-isotriviality condition and interpret it in terms of moduli spaces of curves. For convenience, we focus on ``isotrivial'', and view ``non-isotrivial'' as the negation of ``isotrivial''.

Fix $g>1$ and a function field $K$ of one variable over $k$ in the following. 
Fix an embedding $\bar k\hookrightarrow \OK$ extending $k\hookrightarrow K$. 
For the sake of level structures, fix an integer $N\geq 3$ non-divisible by $\charr(k)$. 

Consider the pair $(C,\alpha)$ consisting of a smooth projective curve $C$ of genus $g$ over $\OK$, and a line bundle 
$\alpha\in \Pic(C)$ of degree 1.
We say that $(C,\alpha)$ is \emph{isotrivial} over $\overline k$ if $(C,\alpha)$ is isomorphic to the base change from $\overline k$ to $\overline K$ of some pair $(C_0,\alpha_0)$ consisting of a smooth projective curve $C_0$ over $\overline k$ and a line bundle $\alpha_0\in \Pic(C_{0})$ of degree 1. 

By a \emph{triple of degree $d$ and level $N$ over a field $F$}, we mean a triple $(C,\eta, \alpha)$ consisting of the following data:
\begin{enumerate}[(1)]
\item a geometrically integral, smooth and projective curve $C$ of genus $g$ over $F$, 
\item an element $\alpha\in \Pic(C)$ of degree $d$ on $C$,
\item an isomorphism $\eta:(\ZZ/N\ZZ)^{2g}\to \Pic(C)[N]$ compatible with the symplectic forms in the following sense.
Denote by $\mu_N(k)$ the group of $N$-th roots of unity in $k$. 
By the Weil pairing $e_N:\Pic(C)[N]\times \Pic(C)[N] \to \mu_N(k)$, 
the existence of $\eta$ implies that $\mu_N(k)$ contains all $N$-th roots of unity of $\bar k$.
Fix an isomorphism $i:\ZZ/N\ZZ\to \mu_N(k)$ of abelian groups.
Denote by $\pair{\cdot, \cdot}:(\ZZ/N\ZZ)^{2g}\times (\ZZ/N\ZZ)^{2g} \to \ZZ/N\ZZ$   the standard symplectic form. 
Then we require that there exists $c\in (\ZZ/N\ZZ)^\times$ (depending on $i$) such that $\eta$ and $i$ transfer
$c \pair{\cdot, \cdot}$ to $e_N$.
\end{enumerate}
By a \emph{pair of degree $d$ over $F$}, we mean a pair $(C,\alpha)$ of $C$ and 
$\alpha$ described in (1) and (2).

We say that a triple of degree $d$ and level $N$ (resp. a pair of degree $d$) over $\ol K$ is \emph{isotrivial} if it is isomorphic to the base change of from $\ol k$ to $\ol K$ of a triple of degree $d$ and level $N$ (resp. a pair of degree $d$) over $\ol k$. 
We have the following easy result.

\begin{lem} \label{isotrivial equiv}
Let $(C,\eta, \alpha)$ be a triple of degree $d$ and level $N$ over $\ol K$. 
Then the following are equivalent:
\begin{enumerate}[(a)]
\item $(C, \alpha)$ is isotrivial as a pair of degree $d$; 
\item $(C,\eta, \alpha)$ is isotrivial as a triple of degree $d$ and level $N$, 
\item $(C,\eta, (2g-2)\alpha-d\omega_{C/\OK})$ is isotrivial as a triple of degree $0$ and level $N$.
\end{enumerate}
\end{lem}

\begin{proof}
It is trivial to have (b) $\Rightarrow$ (a) and (b) $\Rightarrow$ (c). 
For (a) $\Rightarrow$ (b), note that if $(C,\alpha)$ is the base change from $\overline k$ to $\overline K$ of some pair $(C_0,\alpha_0)$, then 
$\Pic(C)[N]\simeq \Pic(C_0)[N]$, which can be seen in terms of torsion points of Jacobian varieties.

For (c) $\Rightarrow$ (b), assume that $(C,\eta, (2g-2)\alpha-d\omega_{C/\OK})$
is the base change from $\overline k$ to $\overline K$ of some triple $(C_0,\eta_0, (2g-2)\alpha_0-d\omega_{C_0/\ol k})$.
Then $(2g-2)\alpha=(2g-2)\alpha_0$ in $\Pic(C)$.
Then $\alpha\in \alpha_0+\Pic(C)[2g-2]$. Similar to the above, we have 
$\Pic(C)[2g-2]=\Pic(C_0)[2g-2]$ and thus $\alpha\in \Pic(C_0)$. 
This proves (c) $\Rightarrow$ (b).
\end{proof}

As before, denote by $M_{g,N,k}$ the moduli scheme of smooth curves of genus $g$ with a level-$N$ structure over $k$, which is exactly the moduli space of pairs $(C,\eta)$ over $k$. 
Denote by $\CCC\to M_{g,N,k}$ the universal curve. 
Denote by $\Pic_{\CCC/M_{g,N,k}}$ the relative Picard functor, and denote by 
$\Pic_{\CCC/M_{g,N,k}}^d$ the sub-functor of line bundles of degree $d$ on fibers.
Then both $\Pic_{\CCC/M_{g,N,k}}$ and $\Pic_{\CCC/M_{g,N,k}}^d$ are representable, and $\CJ=\Pic_{\CCC/M_{g,N,k}}^0$ is the relative Jacobian scheme of $\CCC\to M_{g,N,k}$. 
We refer to \cite[\S8.1,\S8.2]{BLR} for these standard results. 

Consider the moduli space of triples $(C,\alpha,\eta)$ of degree $d$ and level $N$ over $k$. 
Then it is represented by $\Pic_{\CCC/M_{g,N,k}}^d$ over $k$, which can be checked by considering relative representability of the functor 
$(C,\alpha,\eta)\mapsto (C,\eta)$.
For our purpose, it suffices to note that $\Pic_{\CCC/M_{g,N,k}}^d(F)$
is bijective to the isomorphism classes of 
triples $(C,\alpha,\eta)$ of degree $d$ and level $N$ over $F$ for any algebraically closed field $F$ over $k$. 

Finally, our key result in terms of the Picard functor is as follows. 
A triple $(C,\alpha,\eta)$ of degree $d$ and level $N$ over $\ol K$ is isotrivial if and only the point of $\Pic_{\CCC/M_{g,N,k}}^d(\ol K)$ representing the triple $(C,\alpha,\eta)$ lies in the image of the natural map 
$$\Pic_{\CCC/M_{g,N,k}}^d(\ol k)\lra
\Pic_{\CCC/M_{g,N,k}}^d(\ol K).$$
This is further equivalent to the property that the point of $\Pic_{\CCC/M_{g,N,k}}^d(\ol K)$ representing the triple $(C,\alpha,\eta)$ corresponds to a closed point of $\Pic_{\CCC/M_{g,N,k}}^d$.
The proof is trivial.

\subsubsection{An easy lower bound} 

Here we present a non-optimal but easy uniform lower bound of the stable Faltings height in the function field case. 

\begin{lem}\label{lower bound function field}
Let $K$ be a function field of one variable over a field $k$. 
Let $C$ be a geometrically integral, smooth and projective curve over $K$ of genus $g> 1$. Assume that $C_{\OK}$ is non-isotrivial over $\overline k$.
Then
$$
h_\Fal(C) >  
\begin{cases}
3^{-4g^2} & \text{ if } \charr(K)\neq 3,\\
4^{-4g^2} & \text{ if } \charr(K)= 3.
\end{cases}
$$ 
\end{lem}
\begin{proof}
Let $K'$ be a finite extension of $K$ such that $C$ has semistable reduction over $K'$, and let $B'$ be the unique projective regular curve over $k$ with function field $K'$.
Then we have
$$
h_\Fal(C)= \frac{1}{[K':K]} \deg_{B'/k}(\underline\omega_{B'})
= \frac{1}{[K':K]} \deg_{B'/k}(\lambda_{B'}) \geq \frac{1}{[K':K]}.
$$
The inequality is based on two facts. 
First, the non-isotriviality condition implies that $h_\Fal(C)>0$. 
This well-known result is a consequence of the bigness of $\lambda_{\ol S}$ (for $\ol S=B'$) described at the end of \S\ref{sec bigness geometric}. 
Alternatively, it can be derived from \cite[XI, Thm. 4.5]{MB1} and Torelli's theorem. 
Second, $\deg_{B'/k}(\underline\omega_{B'})>0$ implies $\deg_{B'/k}(\underline\omega_{B'})\geq1$ because it is an integer.

With the claim, it suffices to find $K'$ such that $[K':K]$ is small enough for the bound.
Let $N\geq3$ be an integer non-divisible by $\charr(K)$, and set $K'$ to be the smallest subfield of $\ol K$ such that all points of $J(\OK)[N]$ are defined over $K'$. Here $J$ is the Jacobian variety of $X$ over $K$. 
By Raynaud's theorem (cf. \cite[Prop. 4.7]{Gro}), $J$ has semistable reduction over $K'$.
By \cite[Thm. 2.4]{DM}, $C$ has semistable reduction over $K'$. 
So $K'$ satisfies the requirement.

It remains to bound $[K':K]$. Note that $\Gal(K'/K)$ has a natural injection to
$\GSp(J(\OK)[N])\simeq \GSp_{2g}(\ZZ/N\ZZ)$. 
This has a very rough upper bound $N^{4g^2}$.
\end{proof}

\subsection{The uniform Bogomolov-type theorem: proof} \label{sec small points2}

Now we prove Theorem \ref{small points2}. 
Our proof consists of three parts, which treat respectively the number field case, the function field case, and the uniformity of the constants $c_1,c_2$ on $K$.
It is possible to merge them into a single part, but we introduce them case by case for the sake of readability.

\subsubsection{Part 1: number fields}  \label{sec number field}

In this part, we treat the number field case.
Fix an integer $N\geq 3$, and denote by $M_{g,N,\QQ}$ the (fine) moduli space of smooth curves of genus $g$ over $\QQ$ with a full level-$N$ structure. 
As before, $M_{g,N,\QQ}$ is a smooth quasi-projective variety over $\QQ$, which follows from the GIT construction in \cite[\S7.4]{MFK}. 
Set $S=M_{g,N,\QQ}$ and let $\pi:X\to S$ be the universal curve.

As above, consider the situation of Theorem \ref{isomorphism3}(3). 
Then we have a morphism
$$\tau:J\times_SX \lra J\times_SJ, \quad (y,x)\longmapsto (y,y+(2g-2)x-\omega_{X/S}).$$
Rewrite it as a $J$-morphism 
$$
\tau:X_J\lra J_J.
$$
Here we write
$X_J=J\times_S X$ and $J_J=J\times_S J$, viewed as $J$-schemes via the first projections $q_1:J\times_S X\to J$ and $p_1:J\times_S J\to J$.

Denote 
$$\OL=\tau^*(\OTheta_J),\quad
\OTheta_J=p_2^*\OTheta.$$ 
Here $p_2:J\times_S J\to J$ is the second projection. 
We claim that for any adelic line bundle $\OM$ on $J$, there are constants  
$c_1,c_2>0$ such that 
for any $y\in J(\ol K)$,
$$
\#\{z\in X_J(\ol K): q_1(z)=y, \, h_\OL(z)\leq c_1\, h_\OM(y)\} \leq c_2.
$$

To prove the claim, by Theorem \ref{bigness7}(3), $\OL=\tau^*(\OTheta_J)$ is potentially big on $X_J/J$. 
Apply Theorem \ref{small points} to the family $q_1:X_J\to J$. 
Then there are a non-empty Zariski open subvariety $U\subset J$ and constants  
$c_1,c_2>0$ such that the claim holds 
for all $y\in U(\ol K)$.  
This is already very close to the claim, except that it restricts $y$ to be in the open subset $U$ of $J$.
To cover all algebraic points of $J$, we we need to apply a similar result to every irreducible component of $J\setminus U$, and repeat the process finitely many times.
For that, it suffices to prove that, for any (non-empty) closed subvariety $Y$ of $J$, there is a non-empty open subvariety $V$ of $Y$, and constants  
$c_1,c_2>0$ such that for any $y\in V(\ol K)$,
$$
\#\{z\in X_J(\ol K): q_1(z)=y, \, h_\OL(z)\leq c_1\, h_\OM(y)\} \leq c_2.
$$ 
This is again a consequence of Theorem \ref{small points}, applied to the base change
$X_Y \to Y$ of $q_1:X_J\to J$ by $Y\to J$. 
Here the potential bigness condition is obtained by Theorem \ref{bigness7}(4).
This proves the claim. 

Now we verify that the claim implies the theorem for number fields.
In fact, denote by $\tau_1:X_J\to J$ the composition 
$X_J\stackrel{\tau}{\to} J_J\stackrel{p_2}{\to} J$.
This gives $\tau_1^*\OTheta=\OL$. 
Under the identification 
$$
X_J(\ol K)=(J\times_SX)(\ol K)=\{(y,x)\in J(\ol K)\times X(\ol K):\pi_J(y)=\pi(x)\},
$$
$\tau_1:X_J\to J$ maps
a point $z=(y,x)\in X_J(\ol K)$ with $s=\pi(x)\in S(\ol K)$ to the point of $J(\ol K)$ corresponding to the divisor $y+(2g-2)x-\omega_{X_s}$ on $X_s$.
The projection formula for $\tau_1:X_J\to J$ gives
$$
h_\OL(z)=h_\OTheta((2g-2)x-(\omega_{X_s}-y))
=2\,\wh h((2g-2)x-(\omega_{X_s}-y)). 
$$
Then the claim implies
for any $y\in J(\ol K)$ and $s=\pi_J(y)\in S(\ol K)$, 
$$
\#\{x\in X_{s}(\ol K): 2\,\wh h((2g-2)x-(\omega_{X_s}-y))\leq c_1\, h_\OM(y)\} \leq c_2.
$$ 

We claim that there is an adelic line bundle $\OM$ on $J$ such that 
for any $y\in J(\ol K)$ with $s=\pi_J(y)\in S(\ol K)$,
$$
h_\OM(y)\geq 
c_3\big( \max\{h_{\rm Fal}(X_s),1\}+
\wh h(y)\big)
$$
for some constant $c_3>0$.

If this holds, then for any $s\in S(\ol K)$ and $y\in J_s(\ol K)$, 
$$
\#\left\{x\in X_{s}(\ol K): 2\,\wh h((2g-2)x-(\omega_{X_s}-y))\leq c_1c_3\big( \max\{h_{\rm Fal}(X_s),1\}+
\wh h(y)\big)\right\}$$
is less than or equal to $c_2.$
For any $\alpha\in \Pic(X_s)$ of degree 1, set 
$$y=\omega_{X_s}-(2g-2)\alpha.$$
It gives the theorem for number fields.

For the construction of $\OM$, set
$$
\OM= \OTheta +\pi_J^* \ol\lambda_{S} +\CO(c).
$$
Here $c>0$ and $\CO(c)$ is the trivial line bundle on $\ZZ$ with metric given by $\|1\|=e^{-c}$, and viewed as an adelic line bundle on $J$ by pull-back, 
and $\ol\lambda_{S}$ is the adelic Hodge bundle associated to $X\to S$ described in \S\ref{sec hodge bundle}. 
In particular, $\ol\lambda_{S}$ is an adelic line bundle on $S/\ZZ$ computing the Faltings height by Lemma \ref{hodge bundle equivalence} and \cite[\S5.5]{YZ2}.

This choice gives
$$
h_{\OM}(y)
= 2\, \wh h(y)
+ h_{\rm Fal}(X_s) +c.
$$
By a result of Bost \cite{Bos2} (cf. \cite[App.]{GR} or \cite[\S1.3]{JS1}), we have an explicit uniform lower bound 
$$h_{\rm Fal}(X_s)\geq - g\log (\sqrt2 \pi).$$
Note that the constant might be different in some references due to different normalization of the metric defining the Faltings height.   
Then it suffices to take 
$$c= g\log (\sqrt2 \pi)+1.$$ 
This proves the theorem for number fields. 


\subsubsection{Part 2: function fields}

In this part, we prove Theorem \ref{small points2} in the case that $K$ is a function field of one variable over $k$. The proof is similar to the number field case, so we sketch it here and emphasize on the difference.

Fix an integer $N\geq 3$ invertible in $k$, and denote by $M_{g,N,k}$ the (fine) moduli space of smooth curves of genus $g$ over $k$ with a full level-$N$ structure. 
Set $S=M_{g,N,k}$ and let $\pi:X\to S$ be the universal curve.
By \cite[Thm. 2.1]{GO} again, $\pi:X\to S$ has a stable compactification $\ol\pi:\ol X\to \ol S$.

As in the number field case, we still have a $J$-morphism 
$$
\tau: X_J\lra  J_J,
\quad (y,x)\longmapsto (y,y+(2g-2)x-\omega_{X/S}).
$$
Denote $\OL=\tau^*(\OTheta_J)$ as an element of $\wh\Pic(X_J/k)$. 

We claim that for any adelic line bundle $\OM$ on $J/k$, there are constants  
$c_1,c_2>0$ such that 
for any 1-dimensional point $y\in J$,
$$
\#\{z\in (X_J)_y(\ol{k(y)}): h_\OL(z)\leq c_1\, \wh\deg(\OM|_y)\} \leq c_2.
$$ 
Here $(X_J)_y$ is the fiber of $q_1:X_J\to J$ above $y$.
This is similar to the number field case and is obtained by successively applying Theorem \ref{bigness7}(4) and Theorem \ref{small points} to the family $q_1:X_J\to J$ and the adelic line bundle $\OL=\tau^*(\OTheta_J)$.

Denote $s=\pi_J(y)\in S$, viewed as a schematic point. Denote by
$X_s$ and $J_s$ the fibers of $\pi:X\to S$ and $\pi_J:J\to S$ above $s$.
As $z$ lies in the fiber of $J\times_S X\to S$ above $s\in S$,  we can write $z=(y,x)$ with $y\in J_s(\ol{k(y)})$ and $x\in X_s(\ol{k(y)})$. 
Similarly to the number field case,  the claim implies that
for any 1-dimensional point $y\in J$ with $s=\pi_J(y)\in S$, 
$$
\#\{x\in X_{s,k(y)}(\ol{k(y)}): 2\,\wh h((2g-2)x-(\omega_{X_s}-y))\leq c_1\, \wh\deg(\OM|_y)\} \leq c_2.
$$ 
Note that $s=\pi_J(y)$ is either 1-dimensional or 0-dimensional in
$S$, but $X_{s,k(y)}=X_{s}\times_s y$ is viewed as a curve over the function field $k(y)$ such that the N\'eron--Tate height is defined with respect to $k(y)/k$.

Now we take the adelic line bundle $\OM$ on $J$ by
$$
\OM= \OTheta +\pi_J^* \ol\lambda_{S} + N.
$$
Here
$\ol\lambda_{S}$ is the Hodge bundle $\lambda_{\ol S}$ on $\ol S$ defined in
in \S\ref{sec bundle on stable}, viewed as an adelic line bundle on $S/k$;
$N$ is an ample line bundle on a projective compactification of $J$ over $k$, viewed as an adelic line bundle on $J/k$. 
It remains to compute
$$
\wh\deg(\OM|_y)
=  \wh\deg(\OTheta|_y)
+\wh\deg((\pi_J^* \ol\lambda_{S})|_y)
+\wh\deg(N|_y).
$$
For this, resume the above notations $s=\pi_J(y)$ and $z=(y,x)$ with $y\in J_s(\ol{k(y)})$ and $x\in X_s(\ol{k(y)})$. 

First,  
$\OTheta|_y=(\OTheta|_{J_s})|_y$ is the pull-back of the adelic line bundle $\OTheta|_{J_s}$ to $y\in J_s$, and its arithmetic degree is just $2 \wh h(y)$ for $y\in J_{s, k(y)}(k(y))$ over the function field $k(y)$.

Second, 
$(\pi_J^* \ol\lambda_{S})|_y
= \pi_y^*((\ol\lambda_{S})|_s)$.
Here $\pi_y:y\to s$ denotes the natural morphism.
If $s$ is a 0-dimensional point of $S$, then $\pi_y^*((\ol\lambda_{S})|_s)$ has arithmetic degree 0. 
Otherwise, $s$ is a 1-dimensional point of $S$, and thus 
$$
\wh\deg((\pi_J^* \ol\lambda_{S})|_y)
=\deg(\pi_y) \wh\deg((\ol\lambda_{S})|_s)
=\deg(\pi_y) h_{\rm Fal}(X_s)
=h_{\rm Fal}(X_{s, k(y)}).
$$
Third, as $N$ is an ample line bundle on a compactification $\ol J$ of $J$, 
$\wh\deg(N|_y)$ is the degree of $N$ on the Zariski closure of $y$ in $\ol J$. 
This gives a trivial bound $\wh\deg(N|_y)\geq1$.

Hence, for any 1-dimensional point $y\in J$ with $s=\pi_J(y)\in S$,
$$
\wh\deg(\OM|_y)
\geq  2\,\wh h(y)
+ h_{\rm Fal}(X_{s,k(y)})+1.
$$
Then we have proved that  there are constants  
$c_1,c_2>0$ such that
for any 1-dimensional point $y\in J$ with $s=\pi_J(y)\in S$, 
$$
\#\{x\in X_{s,k(y)}(\ol{k(y)}): \wh h((2g-2)x-(\omega_{X_s}-y))\leq c_1 \big( h_{\rm Fal}(X_{s,k(y)})+\wh h(y)+1\big)\} \leq c_2.
$$ 

We can further modify the statement as follows.
For any function field $K$ of one variable over $k$, and for any non-isotrivial $k$-morphism $y':\Spec K\to J$ with image $s'=\pi\circ y':\Spec K\to S$, we have 
$$
\#\{x\in X_{s'}(\ol{K}): \wh h((2g-2)x-(\omega_{X_{s'}}-y'))\leq c_1 \big( h_{\rm Fal}(X_{s'})+\wh h(y')+1\big)\} \leq c_2.
$$ 
Here $y':\Spec K\to J$ is said to be \emph{non-isotrivial} if its image is a 1-dimensional point of $J$. 
The heights are normalized over $K/k$, and this statement is a consequence of the original one by the relation $y=\Im(y')$. It allows $K$ to be any arbitrary finite extension of $k(y)$ (via $y'$), which does not change the truth of the inequality, since a finite extension replaces the heights by multiplying the same constant.

By the moduli property, the above result further implies that for any function field $K$ of one variable over $k$, and for any non-isotrivial triple $(C,\eta, y'')$ of degree 0 and level N over $K$, 
$$
\#\{x\in C(\ol{K}): \wh h((2g-2)x-(\omega_{C}-y''))\leq c_1 \big(h_{\rm Fal}(C)+\wh h(y'')+1
\big)\} \leq c_2.
$$ 
We refer to \S\ref{sec isotrivial} for the equivalence of the non-isotriviality conditions. 
By Lemma \ref{isotrivial equiv}, we can remove $\eta$ from the non-isotriviality condition, so the result holds for any non-isotrivial pair $(C, y'')$ of degree 0 over $K$. 
Moreover, by the relation $y''=\omega_{C}-(2g-2)\alpha$, the result implies that for 
 any non-isotrivial pair $(C,\alpha)$ of degree 1 over $K$, 
$$
\#\{x\in C(\ol{K}):(2g-2)^2 \wh h(x-\alpha)\leq c_1 \big(h_{\rm Fal}(C)+\wh h(\omega_{C}-(2g-2)\alpha)+1
\big)\} \leq c_2.
$$ 
We can further allow $\alpha\in \Pic^1(C_{\OK})$ (instead of $\alpha\in \Pic^1(C)$) by extending $K$.

Finally, we can change $h_\Fal(C)+1$ to $\max\{h_\Fal(C),1\}$ in the result by the fact that $h_\Fal(C)\geq 0$. This non-negativity is easier than Lemma \ref{lower bound function field}, and is discussed in the proof of the lemma.
This proves the theorem for function fields.

\subsubsection{Part 3: uniformity of constants}
The constants $c_1,c_2$ obtained above a priori depend on $(g,K)$, but we are going to prove that they can be chosen uniformly for all $K$. It suffices to treat the uniformity for $K$ varying as a function field.
The key is the following uniform version of Theorem \ref{small points}.

\begin{thm} \label{small points3}
Let $S$ be a flat and quasi-projective integral scheme over $\ZZ$.
Let $\pi:X\to S$ be a smooth relative curve of genus $g>1$.
Let $\OL$ be a nef adelic line bundle on $X/\ZZ$, and $\OM$ be an adelic line bundle on $S/\ZZ$.
If $\OL$ is potentially big on $X/S$, 
then there is a non-empty open subscheme $U$ of $S$, and constants  
$c_1,c_2>0$, such that for any function field $K$ of one variable over a field $k$, 
and for any non-isotrivial morphism $y:\Spec K\to U$,
$$
\#\{x\in X_y(\ol{K}):   h_{\OL|_{X_y}}(x)\leq c_1\,\wh\deg(\OM|_y)\} \leq c_2.
$$ 
\end{thm}

In the theorem, $\OL$ is an adelic line bundle on $X/\ZZ$, so it is a limit of hermitian line bundles $\CLL_i$ on projective models of $X$ over $\ZZ$. However, the hermitian metrics of $\CLL_i$ do not play any essential role in the theorem, since the uniformity comes from the underlying line bundles $\CL_i$ over $\ZZ$. 
However, we state the theorem in the above form to avoid an extra long list of confusing notations. 

Before the proof, we first explain many terms of the statement in the theorem.
Let $K$ be a function field of one variable over a field $k$ as in the theorem. 
Denote by $k'$ the algebraic closure of $k$ in $K$. 
The morphism $y:\Spec K\to U$ is said to be \emph{non-isotrivial} if 
it does not factor through the natural morphism $\Spec K\to \Spec k'$.

In the following, we refer to \cite[\S2.5.5]{YZ2} for various functorial maps on the groups (or categories) of adelic line bundles. 
Denote $p=\charr(k)$.
If $p>0$, the restriction $\OM|_y$ denotes the image via the composition
$$
\wh\Pic(S/\ZZ)\lra \wh\Pic(S_{\FF_p}/\FF_p) \lra \wh\Pic(S_{k}/k)
\stackrel{y'^*}{\lra} \wh\Pic(K/k).
$$
Here $y':\Spec K\to S_k=S\times_\ZZ k$ is induced by $y:\Spec K\to S$ via the fiber product. 
Similarly, if $p=0$, the restriction $\OM|_y$ denotes the image via the composition
$$
\wh\Pic(S/\ZZ)\lra \wh\Pic(S_{\QQ}/\QQ) \lra \wh\Pic(S_{k}/k) \stackrel{y'^*}{\lra} \wh\Pic(K/k).
$$
In both cases, $\wh\deg(\OM|_y)$ is the image of the degree map $\wh\deg: \wh\Pic(K/k)\to \RR$
normalized by multiplicity functions given by degrees over $k$. 

The term $X_y=X\times_S (\Spec K,y)$ is a smooth projective curve over $K$. 
If $p>0$, then
$\OL|_{X_y}$
denotes the image via the composition
$$
\wh\Pic(X/\ZZ)\lra \wh\Pic(X_{\FF_p}/\FF_p) \lra \wh\Pic(X_{k}/k)
\lra \wh\Pic(X_y/k).
$$
If $p=0$, then
$\OL|_{X_y}$
denotes the image via the composition
$$
\wh\Pic(X/\ZZ)\lra \wh\Pic(X_{\QQ}/\QQ) \lra \wh\Pic(X_{k}/k)
\lra \wh\Pic(X_y/k).
$$
In both cases, the height function $h_{\OL|_{X_y}}: X_y(\ol K)\to \RR$ is normalized by $K/k$.

\begin{proof}[Proof of Theorem \ref{small points3}]
The proof is similar to that of Theorem \ref{small points}(2). 
In fact, we still have that $\OL_m=m_\boxtimes\OL$ is big on $X_m=X_{/S}^m$ for some 
$m\geq 1$. 
Then $\OL_m-\epsilon \pi^*\OM$ is big for some positive rational number $\epsilon>0$, and thus some multiple $n(\OL_m-\epsilon\OM)$ of $\OL_m-\epsilon\OM$ is an (integral) adelic line bundle on $X_m/\ZZ$ with a nonzero effective section $s$. 
Denote $Z=\div(s)$. 
Similarly, for any non-isotrivial point $y:\Spec K\to S$,  
$$
\{x\in (X_m)_y(\overline{K}):
h_{\OL_m|_{(X_m)_y}}(x)\leq \epsilon \,\wh\deg(\OM|_y)  \} \subset Z_y(\overline{K}).
$$ 
The remaining part of the proof is also similar.
\end{proof}

Now we prove that we can choose $(c_1,c_2)$ uniformly for all $K/k$ in Theorem \ref{small points}.
Fix an integer $N\geq 3$, and denote by $M_{g,N}$ the (fine) moduli space of smooth curves of genus $g$ over $\ZZ[1/N]$ with a full level-$N$ structure. 
Set $S=M_{g,N}$ and let $\pi:X\to S$ be the universal curve. 
As in the proof of the number field case, we have data $(X, S, J, X_J, J_J, \OTheta, \OTheta_J, \OL)$ over $\ZZ[1/N]$. 
For example, $\OL=\tau^*(\OTheta_J)$ now lies in $\wh\Pic(X_J/\ZZ)$. 

Let $\OM$ be the adelic line bundle over $J/\ZZ$ by
$$
\OM= \OTheta +\pi_J^* \ol\lambda_{S}+\CNN .
$$
Here $\ol\lambda_{S}$ is the Hodge bundle associated to $X\to S$ over $\ZZ$ introduced in \S\ref{sec hodge bundle}, 
and $\CNN$ is a nef hermitian line bundle on a projective model $\CJ'$ of $J$ over $\ZZ$ with an ample underlying line bundle $\CN$ on $\CJ'$, viewed as an adelic line bundle on $J/\ZZ$. 
Note that the the metrics at archimedean places will not play an essential role here.

Apply Theorem \ref{bigness7}(4) and Theorem \ref{small points3} to the family $q_1:X_J\to J$ and the adelic line bundles $\OL$ and $\OM$.  
As in the proof in the function field case, we conclude that 
there are a Zariski closed subset $Z\subsetneq S$ and  constants  
$c_1,c_2>0$, such that
for any function field $K$ of one variable over $k$, and for any non-isotrivial point $y:\Spec K\to J$ whose image $s=\pi\circ y:\Spec K\to S$ lies in $S\setminus Z$, we have 
$$
\#\{x\in X_{s}(\ol{K}): \wh h((2g-2)x-(\omega_{X_{s}}-y))\leq c_1 \big( h_{\rm Fal}(X_{s})+\wh h(y)+1\big)\} \leq c_2.
$$ 
The heights are normalized over $K/k$.
Here we have applied the trivial bound $\wh\deg(\CNN|_y)\geq1$, 
since $\CN$ is an ample line bundle on the projective model $\CJ'$ of $J$ over $\ZZ$.

If the fiber $Z_\QQ$ of $Z$ above $\Spec \QQ$ is non-empty, denote by $Z'$ the Zariski closure of $Z_\QQ$ in $Z$, and
apply Theorem \ref{small points3} to the irreducible components of $Z'$ (instead of $S$). 
This lowers $\dim Z_\QQ$ in the above inequality. 
Therefore, we can assume that $Z_\QQ$ is empty in the above inequality. 

As a consequence, $Z$ lies in finitely many fibers of $S$ over $\Spec \ZZ$.
Then we conclude that there is a positive integer $N'$ such that 
the constants $(c_1,c_2)$ in Theorem \ref{small points} are uniform for all $K/k$ with $\charr(k)\nmid N'$. 

It remains to treat the finitely many characteristics given by $\charr(k)\mid N'$. 
Then we can fix a prime number $p$, and treat all $K/k$ with $\charr(k)=p$. 
This is proved similarly, except that we work over a fine moduli space $S=M_{g,N'', \FF_p}$ over $\FF_p$ (with $p\nmid N''$). Accordingly, Theorem \ref{small points3} can be modified to fit the setting that $S$ is quasi-projective over $\FF_p$ and $\OL$ is adelic over $X/\FF_p$. 
We omit the details here.

\subsection{Uniform fiberwise bigness} \label{sec pointwise}

The goal of this section is to treat Theorem \ref{fiberwise11}. 

Let us first sketch an idea to prove Theorem \ref{fiberwise11} by Theorem \ref{bigness55}
(or more directly Theorem \ref{bigness1}). 
We only take the number field case for example, and the exposition is similar to that in \S\ref{sec number field}. 
Take 
 $S$ to be a fine moduli space $M_{g,N,\QQ}$ over $\QQ$ with $N\geq 3$, and take $\pi:X\to S$ to be the universal curve.
By Theorem \ref{bigness1}, the adelic line bundle
$$\ON=\pi_*\pair{\ol\omega_{X/S,a}\, ,\ol\omega_{X/S,a}}$$
is big on $S/\ZZ$.
On the other hand, the adelic line bundle 
$$
\OM= \ol\lambda_{S} +\CO(c)
$$
is big on $S/\ZZ$ for sufficiently large constant $c>0$. 
This is a consequence of \cite[Lem. 5.2.10]{YZ2}, by the fact that the Hodge bundle $\lambda_{\ol S}$ on $\ol S/\QQ$ for a stable compactification $\ol\pi:\ol X\to\ol S$ is big on $\ol S/\QQ$. 
Once the adelic line bundles are big, the height inequality in \cite[Thm. 5.3.5(1)]{YZ2} implies that there are a Zariski closed subset $Z\subsetneq S$ and constants $a_1, a_2>0$ such that 
$$
a_1 h_{\OM}(s) \leq h_{\ON}(s)\leq a_2 h_{\OM}(s),\quad \forall s\in (S\setminus Z)(\ol \QQ).
$$ 
Applying the height inequality to $(\OM|_Z,\ON|_Z)$ repeatedly, we can assume that $Z$ is empty. 
Now Theorem \ref{fiberwise11} follows from the equalities
$$
h_{\ON}(s)= [\ol\omega_{X_s,a}^2],\quad
h_{\OM}(s)
= h_{\rm Fal}(X_s) +c, \quad \forall s\in S(\ol K).
$$
Here we also need Bost's lower bound $h_{\rm Fal}(X_s)\geq - g\log (\sqrt2 \pi).$

It is hard to get explicit constants $a_1, a_2$ from the above proof, so the goal of this section is to prove Theorem \ref{bigness1} in a different way to obtain constants as explicit as possible. 
Due to different settings, we divide the treatment into the function field case (Theorem \ref{fiberwise2}) and the number field case (Theorem \ref{fiberwise1}). 
To avoid extra complication of the proof, we do not aim for the most optimal constants.

Note that Wilms \cite[Cor. 1.4]{Wil2} claims a similar lower bound of 
$\ol\omega_{C/K,a}^2$ by $h_\Fal(C)$ in the number field case. 
But his proof is incomplete, since \cite[Cor. 1.2]{Wil2} does not cover all degeneration types of the $\varphi$-invariant of Riemann surfaces. However, the idea of his proof works well in the function field case (due to absence of the archimedean invariants). We will essentially follow his idea in the function field case, and our approach to the number field case also uses his idea to get relatively optimal constants.

\subsubsection{Function fields}

We start with the following explicit bounds in the case of function fields, which are previously known to experts. 

\begin{thm}[Theorem \ref{fiberwise11}: function fields]\label{fiberwise2}
Let $K$ be a function field of one variable over a field $k$. 
Let $C$ be a geometrically integral, smooth and projective curve of genus $g> 1$ over $K$. 
Let $K'$ be a finite extension of $K$ such that $C_{K'}$ has a semistable minimal regular model $\CCC'$ over $S'$, where $S'$ is the unique regular projective curve over $k$ with function field $K'$.
Then
$$
\frac{1}{12} \, h_\Fal(C)\leq \ol\omega_{C/K,a}^2 \leq \frac{1}{[K':K]}\omega_{\CCC'/S'}^2 \leq 12 \, h_\Fal(C).
$$ 
If furthermore $C_{\ol K}$ is non-isotrivial over $\ol k$, then 
$$\frac{1}{12} \max\{ h_\Fal(C), 4^{-4g^2}\}\leq 
\ol\omega_{C/K,a}^2 \leq 12 \, \max \{ h_\Fal(C), 1\}.$$ 
\end{thm}
\begin{proof}
The non-isotrivial case is a consequence of the general case by Lemma \ref{lower bound function field}. It suffices to prove the first chain of inequalities.
By base change, we can assume that $K=K'$, and we further denote $(\CCC,S)=(\CCC',S')$. 

Recall the global $\varphi$-invariant 
$$
\varphi(C)  = \sum_{v} \varphi_v(C) \deg(v),
$$
where the summation is over all places of $K$.
Denote the global $\epsilon$-invariant $\epsilon(C)$, and the global $\delta$-invariant 
$\delta(C)=\deg(\Delta_{ S})$ similarly. 
As recalled in \S\ref{sec 3}, we have Zhang's formula (from \cite[Thm. 4.4]{Zha1})
$$
\ol\omega_{C/K,a}^2 = \omega_{\CCC/ S}^2 -\epsilon(C)
$$
and the Noether formula
$$
\omega_{\CCC/ S}^2 = 12\deg(\lambda_{S})-\delta(C)
= 12\, h_\Fal(C)-\delta(C).
$$
It follows that 
$$
\ol\omega_{C/K,a}^2 
= 12\, h_\Fal(C)-\delta(C)-\epsilon(C).
$$

By the positivity $\epsilon(C)\geq 0$ and $\delta(C)\geq0$, 
$$
\ol\omega_{C/K,a}^2\leq \omega_{\CCC/S}^2 \leq 12\, h_\Fal(C).
$$
This gives two of the inequalities of the theorem.

For the inequality in the other direction, the case $\dim S=1$ of Theorem \ref{bigness2} gives
$$
\ol\omega_{C/K,a}^2\geq \frac{3}{5(2g-1)(3g-1)} \, h_\Fal(C).
$$ 
This is much weaker than the one in the theorem, but it is also easy to get stronger bounds from the literature following the idea of Wilms \cite{Wil2}. 

First, we have the inequality
$$
\epsilon(C)\leq \frac12\delta(C)+2\varphi(C),
$$
which is more optimal than $\epsilon(C)\leq (2g-2)\delta(C)$ (cf. Lemma \ref{effective1}) for our purpose here.
To see the inequality, note that de Jong \cite[Prop. 9.2]{dJo3} proves
$$
\delta+\epsilon-2\varphi=\frac32\delta-6\tau
$$
for invariants of polarized metrized graphs. 
Here the $\tau$-invariant is introduced by 
Baker--Rumely in \cite [Thm. 14.1]{BR}, and its positivity $\tau\geq 0$ is a consequence of 
\cite [Lem. 14.4]{BR}. 

Second, let $c(g)>0$ be a constant depending only on $g$ satisfying Cinkir's bound 
$$
\varphi(C) \geq c(g)\delta(C).
$$
We will discuss about the value of $c(g)$ later. 
Then the above inequality gives 
$$
\delta(C)+\epsilon(C)\leq  \big(\frac{3}{2c(g)}+2\big)\varphi(C).
$$
It follows that
$$
\ol\omega_{C/K,a}^2 
= 12\, h_\Fal(C)-\delta(C)-\epsilon(C)
\geq 12\, h_\Fal(C)-  \big(\frac{3}{2c(g)}+2\big)\varphi(C).
$$

Third, by \cite[Prop. 6.1]{LSW}, we have
$$
\ol\omega_{C/K,a}^2 \geq \frac{g-1}{2g+1}  \varphi(C). 
$$
Note that this gives a stronger constant than \cite[Cor. 1.4]{dJo3}. 
Take a linear combination of the last two bounds. We have
$$
\left(1+   \big(\frac{3}{2c(g)}+2\big)\frac{2g+1}{g-1} \right)\ol\omega_{C/K,a}^2 \geq 12\, h_\Fal(C). 
$$

Now we estimate the value of 
$$
A(g)=1+   \big(\frac{3}{2c(g)}+2\big)\frac{2g+1}{g-1}. 
$$
By \cite[Thm. 2.11]{Cin1} and \cite[p. 318, (iv)]{Cin2} (for $g=3$), we can take 
$$
c(2)=\frac{1}{27}, \quad
c(3)=\frac{17}{288}, \quad
c(4)= \frac{3}{88},\quad
c(g)\geq \frac{1}{25},\, \forall g\geq5.
$$
This gives
$$
A(2)=213.5, \quad
A(3)\approx 96.94, \quad
A(4)= 139,\quad
A(g)\leq A(5)=109.625,\, \forall g\geq5.
$$
Therefore, for $g> 2$, we have $A(g)\leq 139< 144$, and thus the inequality implies
$$
139\, \ol\omega_{C/K,a}^2 \geq 12\, h_\Fal(C),  \quad
12\, \ol\omega_{C/K,a}^2 \geq  h_\Fal(C).
$$

If $g=2$, then $C$ is hyperelliptic. Zhang \cite[Cor. 1.3.3]{Zha3}  actually gives
the equality
$$
\ol\omega_{C/K,a}^2 = \frac{2(g-1)}{2g+1}  \sum_{v} \varphi_v(C) \deg(v).
$$
Replace the bound of  \cite[Prop. 6.1]{LSW} by this equality to run the above process. 
The eventually result gives 
$$
\left(1+   \big(\frac{3}{2c(g)}+2\big)\frac{2g+1}{2(g-1)} \right)\ol\omega_{C/K,a}^2 \geq 12\, h_\Fal(C). 
$$
For $g=2$, the coefficient on the left-hand side is $107.25$, and thus 
$$
9\,\ol\omega_{C/K,a}^2 \geq h_\Fal(C). 
$$
This finishes the proof.
\end{proof}

\begin{remark} \label{fiberwise remark}
\begin{itemize}
\item[(1)]
The above proof gives the inequality 
$$ \omega_{C/K,a}^2 \geq \frac{12}{A(g)}  h_\Fal(C).$$ 
The constant of \cite[Thm. 2.11]{Cin1} satisfies $\displaystyle c(g)\to \frac{1}{14}$ as $g\to \infty$, so the coefficient
$\displaystyle\frac{12}{A(g)}\to \frac{12}{47}$ as $g\to \infty$. 

\item[(2)]
There is an alternative approach of the theorem in terms of the slope inequality
$$\omega_{\CCC/ S}^2 \geq \big(4-\frac4g\big) h_\Fal(C),$$
which is proved by Cornalba--Harris \cite[Thm. 1.3]{CH} and Xiao \cite{Xia} in characteristic zero and by Moriwaki \cite{Mor} in positive characteristics. 
See also \cite[Thm. 3.1]{YZt2} for a different proof. 
Then the inequalities used in the proof of Theorem \ref{fiberwise2} can be used to transfer the slope inequality to an inequality of type $ \omega_{C/K,a}^2 \geq  B(g)\, h_\Fal(C)$. 
It turns out that the constant $B(g)$, expressed in terms of $B(g)$ by this method, is weaker than the constant $\displaystyle\frac{12}{A(g)}$ in the theorem. 
\end{itemize}
\end{remark}

\subsubsection{Admissible Noether formula}

In the above proof, the Noether formula plays an important role. To prepare a parallel proof in the number field case, let us first introduce an adelic version of Faltings' arithmetic Noether formula. 

Let $C$ be a smooth projective curve of genus $g>1$ over a number field $K$. 
Assume that $C$ has semistable reduction over $O_K$, and let $\pi:\CCC\to \Spec O_K$ be the minimal regular model of $C$. Denote by $\ol\omega_{\CCC/O_K, \Ar}$ the relative dualizing sheaf over $O_K$ endowed with the Arakelov metric at archimedean places. 
Recall that the arithmetic Noether formula of Faltings \cite[Thm. 6]{Fal2} asserts  
$$
12\, \wh\deg(\ol\lambda_{O_K})= 
\ol\omega_{\CCC/O_K, \Ar}^2+ \delta(C). 
$$
Here the Hodge bundle
$\lambda_{O_K}=\det\pi_*\omega_{\CCC/O_K}$
is endowed with the determinant metric in Lemma \ref{hodge bundle equivalence}, which particularly gives
$$
\wh\deg(\ol\lambda_{O_K})=[K:\QQ]h_\Fal(C).
$$ 
The global $\delta$-invariant 
$$
\delta(C)= \sum_{v} \delta_v(C) \log N_v,
$$
where the summation is over all places $v$ of $K$, and we take the convention that $\log N_v=1$ for real $v$ and $\log N_v=2$ for complex $v$.

Now let us consider the possibility of globalizing the above formula to an adelic version, 
where we replace $C\to \Spec K$ by a more general smooth relative curve $X\to S$ over $\ZZ$.
An obstacle is that the Arakelov canonical bundle $\omega_{\CCC/O_K, \Ar}$ cannot be globalized to form an adelic line bundle on $X/\ZZ$ (even assuming existence of a stable compactification). 
Our remedy of this issue is to replace $\omega_{\CCC/O_K, \Ar}$ by the admissible $\omega_{C/K, a}$, 
which can be globalized by construction. 
For this purpose, recall that
Zhang's formula (cf. \cite[Thm. 4.4]{Zha1}) gives
$$
\ol\omega_{\CCC/O_K, \Ar}^2=\ol\omega_{C/K, a}^2+ \epsilon(C),
$$
where the global $\epsilon$-invariant 
$$
\epsilon(C)= \sum_{v} \epsilon_v(C) \log N_v. 
$$
Note that $\epsilon_v(C)$ is generally only defined for non-archimedean $v$, but we take the convention $\epsilon_v(C)=0$ for archimedean $v$.
Then the arithmetic Noether formula
becomes
$$
12\, \wh\deg(\ol\lambda_{O_K})= 
\ol\omega_{C/K, a}^2+ \delta(C)+\epsilon(C). 
$$
This is the formula we will globalize in the following.

Let $S$ be a flat and quasi-projective normal integral scheme over
$\ZZ$ or over $\QQ$.
Let $\pi:X\to S$ be a smooth relative curve of genus $g>1$.
As in \S\ref{sec bundle on stable}, 
the Noether formula gives a semi-canonical isomorphism
$$
12\lambda_S\lra  \pi_{*}\pair{\omega_{ X/ S},\omega_{ X/ S}}
$$
of line bundles on $S$. 
The semi-canonical isomorphism is unique up to multiplication by $\{\pm 1\}$, and thus compatible with base change up to $\{\pm 1\}$. 
We will not worry about this ambiguity by $\{\pm 1\}$, since it does not contribute to divisors or metrics in our consideration. 
Note that the arithmetic Noether formula is actually defined based by a semi-canonical isomorphism on the generic fiber. 

As in \S\ref{sec hodge bundle}, the Hodge bundle $\lambda_{S}$ and the Faltings metric form an adelic line bundle $\ol\lambda_{S}$ on $S/\ZZ$.
The inverse of the semi-canonical isomorphism defines a section $t$ of 
$$
12\lambda_S-  \pi_{*}\pair{\omega_{ X/ S},\omega_{ X/ S}}.
$$
Define 
$$\ol\Xi_S=\wh\div(t)$$ 
to be the adelic divisor on $S/\ZZ$ respect to the adelic line bundle
$$
12\ol\lambda_S-  \pi_{*}\pair{\ol\omega_{ X/ S,a},\ol\omega_{ X/ S,a}}.
$$
Then the underlying divisor $\Xi_S=0$ on $S$.
By definition, the semi-canonical isomorphism induces an isomorphism
$$
12\ol\lambda_S\lra  \pi_{*}\pair{\ol\omega_{X/ S,a},\ol\omega_{ X/ S,a}} + \CO(\ol\Xi_S).
$$
This is our globalized Noether formula. 
For convenience, we will call it the \emph{admissible Noether formula}.

Compare the admissible Noether formula with the original Noether formula and the arithmetic Noether formula. 
The conclusion is that the total Green's function $\wt g_{\ol\Xi_S}:S^\an\to \RR$ satisfies 
$$
\wt g_{\ol\Xi_S}(v)=  \delta_v(X_{H_v})+\epsilon_v(X_{H_v})
$$
for any discrete or archimedean valuation $v\in S^\an$ with $e_v=e$
(which normalizes the valuation). 
In this sense, the adelic divisor $\ol\Xi_S$ globalizes the invariants
$\{\delta_v+\epsilon_v\}_{v\nmid\infty}$ and $\{\delta_v\}_{v\mid\infty}$ of the setting $C/K$.

We further assume that $\wt\pi:X_\QQ\to S_\QQ$ has a stable compactification 
$\ol\pi:\ol X_\QQ\to \ol S_\QQ$ over $\QQ$. 
By the construction in \S\ref{sec 3}, the semi-canonical isomorphism induces an isomorphism 
$$
12\wt\lambda_{S_\QQ}\lra \tilde \pi_{*}\pair{\wt\omega_{X_\QQ/ S_\QQ,a},\wt\omega_{ X_\QQ/ S_\QQ,a}} + 
\CO(\wt E_{S_\QQ}+\Delta_{\overline S_\QQ}).
$$
Here we use another set of notations to avoid confusion;
namely, we write $\wt E_{S_\QQ}$ for $\ol E_{S_\QQ}$ on $S_\QQ/\QQ$, and write 
$\wt\omega_{ X_\QQ/ S_\QQ,a}, \wt\lambda_{S_\QQ},  \wt \Xi_{S_\QQ}$ for the geometric parts of $\ol\omega_{ X/ S,a}, \ol\lambda_S, \ol \Xi_S$, i.e., their images under the functorial maps (or functors)
$$\wh\Picc(Y/\ZZ)\lra \wh\Picc(Y_\QQ/\QQ), \quad 
\wh\Div(Y/\ZZ)\lra \wh\Div(Y_\QQ/\QQ)$$
defined in \cite[\S2.5.5]{YZ2}.
As a consequence, we have an equality
$$
\wt \Xi_{S_\QQ}=\wt E_{S_\QQ}+\Delta_{\overline S_\QQ}
$$
of adelic divisors on $S/\QQ$.

\begin{remark}
\begin{itemize}
\item[(1)]
As the $\epsilon$-invariant  has non-trivial non-archimedean component but trivial archimedean component, we can prove that it cannot be globalized in the above sense. 
Then the $\delta$-invariant cannot be globalized either. 
\item[(2)]
As in Remark \ref{song's work}, our proof of Theorem \ref{positive lower bound} and the work of Song \cite{Son} demonstrate that the growth of the archimedean $\varphi$-invariant along the boundary of the moduli space of curves is controlled by graph-theoretic data from $\wt\Phi$. 
Replace the adelic divisor $\ol\Phi$ by the adelic divisor $\ol\Xi$ in this process. 
Then the growth of Faltings' archimedean $\delta$-invariant 
 is controlled by graph-theoretic data from the non-archimedean invariant $\delta+\epsilon$. 
This partially explains the results of de Jong \cite[Thm. 1.1]{dJo4} and Faltings \cite[Thm. 7]{Fal4}. 
\end{itemize}
\end{remark}

\subsubsection{Number fields}

Now we are ready to prove the following more explicit version of the number field case of Theorem \ref{fiberwise11}. Note that the main coefficients $1/12$ and $12$ are the same as those in Theorem \ref{fiberwise2}, and that the constants $c_0,c_5$ both come from lower bounds of natural continuous functions on the (complex) moduli space of curves of genus $g$ over $\CC$.

\begin{thm} [Theorem \ref{fiberwise11}: number field case] \label{fiberwise1}
Let $g>1$ be an integer. 
Let $c_0>0$ be the constant depending only on $g$ as in Theorem \ref{positive lower bound}. 
Then there is a constants $c_5$ depending only on $g$ such that for any geometrically connected, smooth and projective curve $C$ of genus $g$ over a number field $K$, one has
\begin{align*}
 \max\left\{ \frac{1}{12} \, h_\Fal(C)+c_5,\, \frac{g-1}{2g+1} c_0\right\}
 \leq &\, [\ol\omega_{C/K,a}^2] \\
 \leq &\, [\ol\omega_{\CCC'/O_{K'},\Ar}^2]
 \leq   12 h_\Fal(C)+ 2g\log (2\pi^4).
\end{align*}
Here $K'$ is any finite extension of $K$ such that $C_{K'}$ has a semistable minimal regular model $\CCC'$ over $O_{K'}$, $\ol\omega_{\CCC'/O_{K'}, \Ar}$ is the relative dualizing sheaf endowed with the Arakelov metric, and 
$$[\ol\omega_{C/K,a}^2]=\frac{1}{[K:\QQ]} \ol\omega_{C/K,a}^2,\quad
[\ol\omega_{\CCC'/O_{K'},\Ar}^2]=\frac{1}{[K':\QQ]} \ol\omega_{\CCC'/O_{K'},\Ar}^2.$$ 
\end{thm}

\begin{proof}
The proof is an arithmetic version of that of Theorem \ref{fiberwise2}, and it
uses many pieces of the proof of Theorem \ref{bigness55} (or Theorem \ref{bigness1}). 
By base change, we can assume $K'=K$, and we write $\CCC$ for $\CCC'$. 
We need to prove the following three inequalities:
\begin{enumerate}[(a)]
\item $\displaystyle [\ol\omega_{C/K,a}^2] \geq \frac{g-1}{2g+1} c_0$,
\item $\displaystyle [\ol\omega_{C/K,a}^2] \geq \frac{1}{12} \, h_\Fal(C)+c_5,$
\item $\displaystyle [\ol\omega_{C/K,a}^2] \leq  [\ol\omega_{\CCC/O_{K},\Ar}^2]
\leq 12 h_\Fal(C)+  2g\log (2\pi^4)$.
\end{enumerate}

For (a), Wilms \cite[Thm. 1.2]{Wil3} has proved
$$
\ol\omega_{C/K,a}^2 \geq \frac{g-1}{2g+1}  \varphi(C)
$$
for the global $\varphi$-invariant
$$
\varphi(C)= \sum_{v} \varphi_v(C) \log N_v.
$$
By Cinkir \cite[Thm. 2.11]{Cin1}, $\varphi_v(C)\geq 0$ for non-archimedean $v$.
By Theorem \ref{positive lower bound}, $\varphi_v(C)\geq c_0$ for archimedean $v$.
This gives (a). 

For (c), as in the last subsection, Zhang's formula (cf. \cite[Thm. 4.4]{Zha1}) gives
$$
\ol\omega_{C/K, a}^2=\ol\omega_{\CCC/O_{K},\Ar}^2- \epsilon(C) \leq \ol\omega_{\CCC/O_{K},\Ar}^2.
$$
And the arithmetic Noether formula gives
$$
\ol\omega_{\CCC/O_{K},\Ar}^2=12\, \wh\deg(\ol\lambda_K)
- \delta(C). 
$$
Note that  $\delta_v(C)\geq 0$ for non-archimedean $v$. 
On the other hand, by Wilms \cite[Cor. 1.2]{Wil1}, the archimedean $\delta$-invariant of a compact Riemann surface $Y$ of genus $g$ satisfies 
$\delta(Y)\geq -2g\log (2\pi^4).$
It follows that 
$$
\ol\omega_{\CCC/O_{K},\Ar}^2
\leq 12\, \wh\deg(\ol\lambda_K)+ 2g\log (2\pi^4)[K:\QQ]. 
$$
This proves inequality (c).

It remains to prove inequality (b). 
We go through the arithmetic analogue of the relevant part of the proof of Theorem \ref{fiberwise2}. 
Combine the formula 
$$
12\, \wh\deg(\ol\lambda_K)= 
\ol\omega_{C/K, a}^2+ \delta(C)+\epsilon(C)
$$
 with Wilms' bound
$$
\ol\omega_{C/K,a}^2 \geq \frac{g-1}{2g+1}  \varphi(C). 
$$
Following the proof of Theorem \ref{fiberwise2}, we take the linear combination
\begin{align*}
& \left(1+   \big(\frac{3}{2c(g)}+3 \big)\frac{2g+1}{g-1} \right) \ol\omega_{C/K,a}^2\\
\geq &\,12  \, \wh\deg(\ol\lambda_K)
-\delta(C)-\epsilon(C)
+ \big(\frac{3}{2c(g)}+3\big) \varphi(C).
\end{align*}
Here $c(g)$ is still Cinkir's constant. 
Note that the coefficient on the left-hand side is taken to be slightly bigger than the constant
$$A(g)=\left(1+   \big(\frac{3}{2c(g)}+2 \big)\frac{2g+1}{g-1} \right)$$
used in the proof of Theorem \ref{fiberwise2}. The purpose of this will be seen later. 

We need to estimate the global invariant
$$
d(C):= \big(\frac{3}{2c(g)}+3\big) \varphi(C)-\delta(C)-\epsilon(C).
$$
For any place $v$, denote the $v$-component
$$
d_v(C)= \big(\frac{3}{2c(g)}+3\big) \varphi_v(C)-\delta_v(C)-\epsilon_v(C). 
$$
Recall the convention $\epsilon_v(C)=0$ for archimedean $v$. 
For any non-archimedean place $v$, the proof of Theorem \ref{fiberwise2}
already gives 
$$
d_v(C)
\geq \varphi_v(C) \geq 0.
$$
It remains to treat the archimedean components.
We claim that there is a constant $c_5'$ depending only $g$ such that for every archimedean $v$, 
$$
d_v(C) \geq c_5'.
$$

Our proof of the claim follows the idea of Theorem \ref{positive lower bound} of using adelic divisors over moduli spaces of curves. 
Let $S$ be a flat and quasi-projective normal integral scheme over $\QQ$.
Let $\pi:X\to S$ be a smooth relative curve of genus $g$ with a stable compactification $\bar\pi:\bar X\to \bar S$. 
Recall that the adelic divisor $\ol \Xi_S$ on $S/\ZZ$ globalizes $\delta+\epsilon$.  
Then the adelic divisor 
$$
\ol D_S=\big(\frac{3}{2c(g)}+3\big) \ol\Phi_S-\ol \Xi_S
$$
on $S/\ZZ$ globalizes the invariant $d(C)$. 
The geometric part, as an adelic divisor on $S/\QQ$, is given by 
$$
\wt D_S=\big(\frac{3}{2c(g)}+3\big) \wt\Phi_S-\wt \Xi_S
=\big(\frac{3}{2c(g)}+3\big) \wt\Phi_S-\wt \Delta_S-\wt E_S.
$$
Here we write $\wt E_{S}, \wt \Delta_{S}$ for the adelic divisors $\ol E_{S}, \Delta_{\ol S}$ on $S_\QQ/\QQ$ defined in \S\ref{sec 3}, and write 
$\wt D_S, \wt \Xi_{S}, \wt\Phi_S$ for the images of $\ol D_S, \ol \Xi_S, \ol\Phi_S$
under the functorial map
$\wh\Div(S/\ZZ)\to \wh\Div(S/\QQ)$.

The key is that $\wt D_S-\wt \Delta_{S}$ is effective on $S/\QQ$. 
In fact, the total Green's function $\wt g_{\wt D_S}: (S/\QQ)^\an\to \RR$ at any discrete valuation $v\in S^\an$ normalized by $e_v=e$ satisfies
\begin{align*}
\wt g_{\wt D_S}(v)
=&\, \big(\frac{3}{2c(g)}+3\big) \varphi_v(X_{H_v})-\delta_v(X_{H_v})-\epsilon_v(X_{H_v})\\
\geq &\,  \varphi_v(X_{H_v})\\
\geq &\,  c(g)\delta_v(X_{H_v})\\
= &\, c(g) \wt g_{\wt \Delta_S}(v).
\end{align*}
Here the
first inequality follows from the corresponding inequality in the proof of Theorem \ref{fiberwise2}
as a consequence of de Jong's identity.
By the density of discrete valuations $v\in S^\an$, the inequality holds for every 
$v\in S^\an$, and thus
$\wt D_S-c(g)\wt \Delta_{S}$ is effective on $S/\QQ$. 

Once $\wt D_S-c(g)\wt \Delta_{S}$ is effective on $S/\QQ$, the method of 
the proof of Theorem \ref{positive lower bound} works here.
In fact, take $S=\CM_{g,N,\QQ}$  to be the fine moduli scheme over $\QQ$ of curves of genus $g$ with a full level-$N$ structure (with $N\geq3$), and take $\pi:X\to S$ to be the universal curve.
By the method, 
the effectivity of the geometric parts implies 
that the archimedean Green's function $g_{\ol D_S,\infty}:S(\CC)\to \RR$ of $\ol \Xi_S$ increases to infinity along the boundary. 
Then  $g_{\ol D_S,\infty}\geq c_5'$ for a constant $c_5'$ depending only on $g$.
This proves the claim. 

By the claim, the original inequality implies
$$ \left(1+   \big(\frac{3}{2c(g)}+3 \big)\frac{2g+1}{g-1} \right) \ol\omega_{C/K,a}^2
\geq 12\, \wh\deg(\ol\lambda_K)+ [K:\QQ] c_5'.
$$
By the computation in the proof of Theorem \ref{fiberwise2}, 
if $g>2$, the coefficient on the left-hand side is
$$
A(g)+\frac{2g+1}{g-1} \leq 139+ \frac{7}{2} <144.
$$
Then we have 
$$
144\, \ol\omega_{C/K,a}^2
\geq 12\, \wh\deg(\ol\lambda_K)+ [K:\QQ] c_5'.
$$
This gives inequality (b) for $g>2$. 

Finally, if $g=2$, as in the proof of Theorem \ref{fiberwise2}, the constant $144$ can be obtained by applying the formula
$$
\ol\omega_{C/K,a}^2 = \frac{2(g-1)}{2g+1}  \varphi(C)
$$
of \cite[Cor. 1.3.3]{Zha3}. 
This finishes the proof.
\end{proof}

\begin{remark}
\begin{itemize}
\item[(1)]  
An arithmetic version of the approach in Remark \ref{fiberwise remark} can also be sued to prove inequality (b). 
In fact, 
Bost \cite[Thm. IV]{Bos1} proved an arithmetic analogue of the slope inequality, where the archimedean function $\psi$ is bounded below by the growth condition in the theorem.
See also \cite[\S3.3.2]{Bur}, \cite[Prop. 5.6]{dJo1} and 
\cite[Cor. E]{YZt1} for similar inequalities without proving growth conditions at archimedean places. 
This can be converted to prove an version of inequality (b). 
As in the geometric case, the coefficient obtained in this way is weaker than our original approach, while the constant $c_5$ is still implicit. 

\item[(2)] 
By \cite[Chap. V, Prop. 4.6]{FC}, the Faltings heights of abelian varieties with semistable reduction over number fields
satisfy the Northcott property. By Torelli's theorem (cf. \cite[Thm. 12.1]{Mil}), the Northcott property transfers to the Faltings heights of smooth projective curves with semistable reduction over number fields. Then the theorem further implies that the invariant $\ol\omega_{C/K,a}^2$ satisfies the Northcott property when varying the smooth projective curve $C$ and the number field $K$. 
\end{itemize}

\end{remark}

\subsection{Non-degeneracy and relative Bogomolov conjecture} 
\label{sec non-degeneracy}

The goal of this subsection is to introduce the consequences of our potential bigness on the non-degeneracy problem and relative Bogomolov conjecture related to the three examples in \S\ref{sec big examples}.
These results are key ingredients in the approaches of
\cite{Gao2, GH, DGH1, DGH2, Kuh} to uniform Bogomolov-type results, but we do not use them directly in our approach.

\subsubsection{Generality on non-degeneracy}

Let $k$ be a \emph{field}.
Let $S$ be a quasi-projective normal variety over $k$.
Let $\psi:A\to S$ be an abelian scheme over $S$ of relative dimension $g$, with identity section $e:S\to A$. 
Let $Y$ be a quasi-projective variety over $k$ with a generically finite morphism
$\iota: Y\to A$ over $k$. 
Let us first review the notion of non-degeneracy of $Y$ in $A$.
For convenience, we do not assume that $\iota$ is a closed immersion.

Let $L$ be a symmetric and $\psi$-ample line bundle on $A$, 
rigidified by an isomorphism $e^*L\simeq \CO_S$.  
Let $\OL$ be the unique adelic line bundle on $A$ extending $L$ and satisfying $[2]^*\OL\simeq 4\OL$ as constructed in \cite[Thm. 6.1.1]{YZ2}.

Following \cite[\S6.2.2]{YZ2}, we say that the morphism $\iota:Y\to A$ is \emph{non-degenerate} if
$\OL|_{Y}=\iota^*\OL$ is a big adelic line bundle on $Y$.
The definition is independent of the choice of $L$, as any two ample line bundles can bound each other up to positive multiples. 
As explained in the loc. cit., if there is an embedding $k\to \CC$, the definition is equivalent to the definition of \cite[Def. 1.5]{DGH1}, which requires the Betti map $Y(\CC)_V \to (\RR/\ZZ)^{2g}$ to have a full rank at some point of $Y(\CC)_V$ for some simply connected open subset $V$ of $S(\CC)$. We refer to Andr\'e--Corvaja--Zannier \cite{ACZ} for a systematic study of the rank of the Betti map with a view towards Diophantine applications. 

As an easy consequence of Theorem \ref{potentially big}, we have the following quick criterion. 

\begin{prop}  \label{potentially non-degenerate1}
Assume that $Y\to S$ is projective and flat with a geometrically integral generic fiber.
If $L|_Y$ is big on the generic fiber of $Y\to S$, and 
$(\pi|_Y)_*\pair{\OL|_Y,\cdots,\OL|_Y}$ is big on $S/k$, then 
$\iota^m:Y_{/S}^m\to A_{/S}^m$ is non-degenerate for all $m\geq \dim S$.
\end{prop}

\subsubsection{Three non-degeneracy examples}

Let $k$ be a \emph{field}.
Let $S$ be a quasi-projective normal variety over $k$.
Let $\pi:X\to S$ be a smooth relative curve over $S$ of genus $g>1$. 
Let $\alpha$ be a line bundle on $X$ of degree $d>0$ on fibers of $\pi:X\to S$. 

Recall from \S\ref{sec big examples} that we have the following three natural morphisms:
\begin{enumerate}[(1)]
\item
the $S$-morphism
$i_\alpha: X\to J,$
\item 
the $X$-morphism 
$i_\Delta: X_X\to J_X,$
\item
the $J$-morphism 
$\tau=i_{\omega-Q}: X_J\to J_J.$
\end{enumerate}
Note that the definition of $Q$ requires a section of $X\to S$, but
the definition of the morphism $\tau$ does not. 
We explain this discrepancy slightly. 
In fact, we usually require a section of $X\to S$ to rigidify $Q$ and thus determine $Q$ as a  unique class in $\Pic(J\times_S X)$. However, without the section, $Q$ can be considered a unique class in $\Pic(J\times_S X)/\Pic(J)$, which still gives a well-defined $i_{\omega-Q}$. 

The morphisms induce for $m\geq1$ three natural morphisms:
\begin{enumerate}[(a)]
\item
the $S$-morphism
$$
i_\alpha^m=(i_\alpha,\cdots, i_\alpha): X^m_{/S}\lra J^m_{/S},
$$
\item 
the $X$-morphism 
$$
i_\Delta^m=(i_\Delta,\cdots, i_\Delta): (X_X)^m_{/X}\lra (J_X)^m_{/X},
$$
\item
the $J$-morphism 
$$
i_{\omega-Q}^m=(i_{\omega-Q},\cdots, i_{\omega-Q}): (X_J)^m_{/J}\lra (J_J)^m_{/J}.
$$
\end{enumerate}

The morphism in (b) is closely related to the Faltings--Zhang morphism
$$
i_{\FZ,m}: X^{m+1}_{/S}\lra J^m_{/S} ,\quad
(x_0,\cdots, x_m)\longmapsto (x_1-x_0,\cdots, x_m-x_0). 
$$ 
The morphism in (c) is closely related to the morphism
\begin{eqnarray*}
\tau_m:(X^{m}_{/S})\times_SJ &\lra& J^m_{/S} , \\
(x_1,\cdots, x_m,y) &\longmapsto& ((2g-2)x_1-\omega_{X/S}+y, x_2-x_1,\cdots,x_m-x_1).
\end{eqnarray*}
Again, the definition of $i_{\omega-Q}^m$ does not require a section of $X\to S$.

Now we have the following theorem, whose parts (1), (2), (4) in the case $\charr(k)=0$ were proved in \cite[Thm. 1.2(i), Thm 1.2']{Gao2}. 

\begin{thm}  \label{potentially non-degenerate2}
Let $k$ be a field.
Let $S$ be a quasi-projective variety over $k$.
Let $\pi:X\to S$ be a smooth relative curve over $S$ of genus $g>1$ with maximal variation.
Let $\alpha$ be a line bundle on $X$ of positive degree on fibers of $\pi:X\to S$. 
Then the following morphisms are non-degenerate:
\begin{enumerate}[(1)]
\item  the $S$-morphism $i_\alpha^m:X^m_{/S}\to J^m_{/S}$ for any $m\geq \dim S$,
\item  the $X$-morphism $i_\Delta^m: (X_X)^m_{/X}\to (J_X)^m_{/X}$ for any $m\geq \dim S+1$, 
\item  the $J$-morphism $i_{\omega-Q}^m: (X_J)^m_{/J}\to (J_J)^m_{/J}$ for any $m\geq \dim S+g$, 
\item  the $S$-morphism $i_{\FZ,m}:X^{m+1}_{/S}\to J^m_{/S}$ for any $m\geq \dim S+1$, 
\item the $S$-morphism $\tau_m:(X^{m}_{/S})\times_SJ\to J^m_{/S}$ for any $m\geq \dim S+g$. 
\end{enumerate}

\end{thm}

\begin{proof}
We can assume that $S$ is normal by replacing it by an open subscheme.
Parts (1)-(3) follow from Theorem \ref{bigness7} and Proposition \ref{potentially non-degenerate1}. 
We will see that (2) and (4) are equivalent, and (3) and (5) are equivalent.

For (4), we need to interpret $i_{\FZ,m}:X^{m+1}_{/S}\to J^m_{/S}$ in terms of the  morphism in (2).
For this purpose, let us recall the $X$-morphism 
$$
i_\Delta: X_X\lra  J_X,\quad (x_0,x_1)\longmapsto (x_0,x_1-x_0).
$$
Under the identifications $(X_X)^m_{/X}=X^{m+1}_{/S}$ and 
$(J_X)^{m}_{/X}=X\times_S J^m_{/S}$, the morphism $i_\Delta^m: (X_X)^m_{/X}\to (J_X)^m_{/X}$
becomes 
$$
i_\Delta^m: X^{m+1}_{/S}\lra  X\times_S J^m_{/S},\quad
(x_0,\cdots, x_m)\longmapsto (x_0,x_1-x_0,\cdots, x_m-x_0).
$$
In the notation, $X\times_S(J^{m}_{/S})$ is viewed as an abelian scheme over $X$ via projection to the first component.  
By forgetting the first component, we obtain the non-degeneracy of $i_{\FZ,m}$.

For (5), taking a generically finite base change $S'\to S$ if necessary, we can assume that $X\to S$ has a section. Then we have a universal line bundle $Q$ on $J\times_SX$. 
Consider the morphism
$$
i_{\omega-Q}: X_J\lra J_J,\quad x\longmapsto (2g-2)x-(\omega_{X_J/J}-Q).
$$
The morphism 
$i_{\omega-Q}^m: (X_J)^m_{/J}\to (J_J)^m_{/J}$
 is just 
$$
i_{\omega-Q}^m:J\times_S(X^{m}_{/S})\lra J\times_S(J^{m}_{/S}),$$ 
$$(y,x_1,\cdots, x_m) \longmapsto (y,(2g-2)x_1-\omega_{X/S}+y,\cdots,(2g-2)x_m-\omega_{X/S}+y).$$
In the notation, $J\times_S(J^{m}_{/S})$ is viewed as an abelian scheme over $J$ via projection to the first component.  
By forgetting the first component, we obtain a non-degenerate morphism 
$$
(i_{\omega-Q}^m)':J\times_S(X^{m}_{/S})\lra J^{m}_{/S},$$ 
$$(y,x_1,\cdots, x_m) \longmapsto ((2g-2)x_1-\omega_{X/S}+y,\cdots,(2g-2)x_m-\omega_{X/S}+y).$$
Compose it further with the endomorphism 
$$
\iota:J^{m}_{/S}\to J^{m}_{/S},\quad
(y_1,\cdots, y_m) \longmapsto (y_1,y_2-y_1,\cdots,y_m-y_1).
$$
We obtain a non-degenerate morphism 
$$
(i_{\omega-Q}^m)'':J\times_S(X^{m}_{/S})\lra J^{m}_{/S},$$ 
$$(y,x_1,\cdots, x_m) \longmapsto ((2g-2)x_1-\omega_{X/S}+y,(2g-2)(x_2-x_1)\cdots,(2g-2)(x_m-x_1)).$$
Removing the factors $2g-2$ in the last $m-1$ components, we end up with the morphism $\tau_m$. 
This finishes the proof.
\end{proof}

\subsubsection{Generality on the relative Bogomolov conjecture}

\kkk
Let $K$ be a {global field over} $k$. 
Namely, $K$ is a number field if $k=\ZZ$; $K$ is a function field of one variable over $k$ if $k$ is a field. 

Let $S$ be a quasi-projective normal variety over $K$.
Let $\psi:A\to S$ be an abelian scheme over $S$ of relative dimension $g$, with identity section $e:S\to A$. 
Let $Y$ be a quasi-projective variety over $K$ with a generically finite morphism
$\iota: Y\to A$ over $K$. 
Assume that the composition $Y\to S$ is surjective, $Y_{\bar \eta}$ is irreducible and generates the algebraic group $A_{\bar \eta}$, where $\bar\eta$ is the geometric generic point of $S$.
Recall that the relative Bogomolov conjecture of \cite[Conj.1.2]{DGH2} asserts that, if $\dim Y<g$, then there is a constant $\epsilon>0$ such that 
$$
Y(L,\epsilon)=\{y\in Y(\ol K):\wh h_{L}(\iota(y))\leq \epsilon\}
$$
is not Zariski dense in $Y$. 

Here $L$ is a symmetric and $\psi$-ample line bundle on $A$, 
rigidified by an isomorphism $e^*L\simeq \CO_S$.
The canonical height $\wh h_L:A(\ol K)\to \RR$ is defined fiberwise over $S$, which is also the height function associated to the adelic line bundle
 $\OL$ in $\wh\Pic(A/k)$ (instead of $\wh\Pic(A/K)$) extending $L$ and satisfying $[2]^*\OL\simeq 4\OL$ as constructed in \cite[Thm. 6.1.1]{YZ2}.
 
The validity of the relative Bogomolov conjecture is independent of the choice of $L$. 
In fact, if $L'$ is another such line bundle on $A$, then there is a positive integer $n$ such that $nL'-L$ and $nL-L'$ are again such line bundles on $A$. 
Then $\wh h_{nL'-L}\geq 0$ and $\wh h_{nL-L'}\geq 0$, and thus $n^{-1}\wh h_{L}\leq \wh h_{L'} \leq n\, \wh h_{L}$.

By the height inequality in \cite[Thm. 5.3.5(1)]{YZ2}, if $\OL|_Y$ is big on $Y$, then $Y$ satisfies the relative Bogomolov conjecture.
As an easy consequence of Theorem \ref{potentially big}, we have the following result on the relative Bogomolov conjecture. 

\begin{prop}  \label{relative bogomolov1}
Assume that $Y\to S$ is projective and flat with a geometrically integral generic fiber.
If $L|_Y$ is big on the generic fiber of $Y\to S$, and 
$(\pi|_Y)_*\pair{\OL|_Y,\cdots,\OL|_Y}$ is big on $S/k$, then 
$\iota^m:Y_{/S}^m\to A_{/S}^m$ satisfies the relative Bogomolov conjecture for all $m\geq \dim S+1$.
\end{prop}

The proposition needs the bigness of adelic line bundles over $k$ instead of over $K$.
Projective models of $S$ over $k$ have dimension $\dim S+1$, which explains the bound of $m$ in the proposition.

\subsubsection{Three examples of the relative Bogomolov conjecture}

The following theorem is analogous to
Theorem \ref{potentially non-degenerate2}.

\begin{thm}  \label{relative bogomolov2}
Let $K$ be either a number field or a function field of one variable. 
Let $S$ be a quasi-projective variety over $K$.
Let $\pi:X\to S$ be a smooth relative curve over $S$ of genus $g>1$ with maximal variation.
Let $\alpha$ be a line bundle on $X$ of positive degree on fibers of $\pi:X\to S$. 
 Then the following morphisms 
satisfy the relative Bogomolov conjecture:
\begin{enumerate}[(1)]
\item  the $S$-morphism $i_\alpha^m:X^m_{/S}\to J^m_{/S}$  for any $m\geq \dim S+1$. 
\item  the $X$-morphism $i_\Delta^m: (X_X)^m_{/X}\to (J_X)^m_{/X}$  for any $m\geq \dim S+2$. 
\item  the $J$-morphism $i_{\omega-Q}^m: (X_J)^m_{/J}\to (J_J)^m_{/J}$  for any $m\geq \dim S+g+1$. 
\item  the $S$-morphism $i_{\FZ,m}:X^{m+1}_{/S}\to J^m_{/S}$  for any $m\geq \dim S+2$. 
\item the $S$-morphism $\tau_m:(X^{m}_{/S})\times_SJ\to J^m_{/S}$  for any $m\geq \dim S+g+1$. 
\end{enumerate}

\end{thm}

\begin{proof}

The proof is very similar to Theorem \ref{potentially non-degenerate2}.
In fact, parts (1)-(3) follow from Theorem \ref{bigness7} and 
Proposition \ref{relative bogomolov1}. 
We also have that (2) and (4) are equivalent, and (3) and (5) are equivalent.
Similar to Proposition \ref{relative bogomolov1}, we need the bigness of the adelic line bundles over $k$ instead of over $K$, where $k=\ZZ$ if $K$ is a number field and $k$ is the field of constants if $K$ is a function field of one variable.
As a consequence, the bounds for $m$ are increased by 1 from those in Theorem \ref{potentially non-degenerate2}.
\end{proof}

\appendix
\section{Admissible metrized line bundles} \label{sec appendix}

The goal of this appendix is to review the theory of Zhang \cite{Zha1} on admissible pairings on curves, and we will use the terminology of adelic line bundles and metrized line bundles introduced in \cite{Zha2}.
Note that the treatment of \cite{Zha1} was written before \cite{Zha2} and was based on graph theory, so we think it is necessary to survey a detailed transfer of the terminology. 
We will use Berkovich spaces of \cite{Ber} for metrics on line bundles, and our approach of the main result (Theorem \ref{admissible1}) works for both the archimedean case and the non-archimedean case.

\subsection{Arakelov metrics in the complex case} \label{sec ara}

We first recall the complex case introduced in Arakelov's original work \cite{Ara}. 
We refer to \cite[II.1, II.2, IV.5]{Lan} for more details on this subject. 
Let $C$ be a compact Riemann surface of genus $g>0$. 

Define a natural hermitian pairing on $\Gamma(C, \omega_{C})$
by 
$$
\pair{\alpha, \beta}=\frac{i}{2}\int_{C} \alpha\wedge \overline\beta.
$$
Let $\alpha_1,\cdots, \alpha_g$ be an orthonormal basis of this pairing.
The \emph{Arakelov K\"ahler form} on $C$ is defined by
$$
d\mu_\Ar=\frac{i}{2g} \sum_{j=1}^g \alpha_j\wedge \overline\alpha_j.
$$
The definition is independent of the choice of the orthonormal basis.

A smooth hermitian metric $\|\cdot\|$ on a line bundle $L$ on $C$ is called \emph{admissible}
if the Chern form $c_1(L,\|\cdot\|)$ satisfies
$$
c_1(L,\|\cdot\|)= \deg(L) \, d\mu_\Ar.
$$

A smooth Green's function of a divisor $D$ on $C$ is a smooth function $g_D:C\setminus |D|\to \RR$ such that for any meromorphic function $f_U$ on an open subset $U$ of $C$ with $D|_U=\div(f_U)$, the function $g_D|_U+\log |f_U|:U\setminus |D|\to \RR$ extends to a smooth function on $U$.   
The smooth Green's function $g_D$ is called \emph{admissible} if the hermitian metric on $\CO(D)$ defined by $\|1\|=\exp(-g_D)$ is admissible. The condition is equivalent to
$$
\frac{i}{\pi}\partial\overline\partial g_D
=d\mu_\Ar-\delta_{D} 
$$
as currents on $C.$
Note that our normalization is different from that of \cite[II.1]{Lan}, where the Green's functions correspond to twice of the Green's functions here.

A smooth hermitian metric on a line bundle on $C^2=C\times C$ is called \emph{admissible} if for any point $x\in C$, the pull-back metrics of the metric to $C\times x$ and $x\times C$ are both admissible. 

Any line bundle on $C$ (resp. $C^2$) admits an admissible metric, which is unique up to multiplicative constants.

The \emph{Arakelov Green's function} on $C^2$ is the unique symmetric and smooth Green's function 
$$g_\Ar:C^2\setminus \Delta\lra \RR$$
of the diagonal $\Delta$ in $C^2$ satisfying the following properties:
\begin{enumerate}[(1)]
\item For any point $x_0\in C$, the Green's function $g_\Ar(x_0,\cdot)$ of $x_0$ on 
$C$ is admissible.
\item For any point $x_0\in C$, the total integral
$$
\int_{C} g_\Ar(x_0,\cdot) d\mu_\Ar=0.
$$
\end{enumerate}
The Green's function gives a hermitian metric $\|\cdot\|_{\Delta,\Ar}$ of $\CO(\Delta)$ on $C^2$ by $\|1\|=\exp(-g_\Ar)$. 
It is actually admissible.

Finally, the \emph{Arakelov metric $\|\cdot\|_\Ar$} of $\omega_{C}$ on $C$ is the unique smooth metric such that the residue map
$$
(\omega_{C} \otimes_{\CO_C} \CO(x_0) )|_{x_0}\lra \CC 
$$
is an isometry. 
Here $\CC$ is endowed with the usual absolute value, 
and $\CO(x_0)$ is endowed with the metric given by $\|1\|=\exp(-g_\Ar(x_0,x))$.

The choice of $d\mu_\Ar$ implies that the metric $\|\cdot\|_\Ar$ is admissible.
This is the reason to choose $d\mu_\Ar$, and determines $d\mu_\Ar$ uniquely.

\subsection{Zhang metrics in the non-archimedean case}

Let $K$ be a non-archimedean field, i.e. a complete field with a non-trivial non-archimedean valuation $|\cdot|$. 

For any projective variety $X$ over $K$, denote by $X^\an$ the Berkovich analytic space of $X$ over $K$.  
We refer to \cite[App. 1]{YZ1} for the notion of continuous (resp. semipositive, integrable) \emph{metrics} of line bundles of $X$ over $K$. 

Let $C$ be a smooth projective curve of genus $g>0$ over $K$.
As a convention, all curves are assumed to be geometrically integral in this appendix.

We start with some terminology of metrics on $C^2$.
Denote by $\Delta$ the diagonal divisor of $C^2$. 
A metric $\|\cdot\|_{\Delta}$ of $\CO_{C^2}(\Delta)$ on $(C^2)^\an$ is called 
\emph{symmetric}
if the Green's function 
$$g_\Delta=-\log\|1\|_{\Delta}: (C^2)^\an\setminus \Delta^\an\to \RR$$ 
is symmetric in the sense that it is invariant under the action on $(C^2)^\an$ induced by the transposition action on $C^2$ by switching the two components.

The Green's function induces a function 
$$g_\Delta: C(\overline K)^2\setminus \Delta(\overline K)\lra \RR$$ 
by the natural maps
$$
C(\overline K)^2 \lra |C^2|_0 \lra (C^2)^\an.$$
Here $|C^2|_0$ denotes the set of closed points of $C^2$.
Then the symmetry can also be understood in the usual sense.

For any finite extension $K'$ of $K$ and any point $x\in C(K')$, denote
$$
(\CO(x), \|\cdot\|_{x}): = i_{x}^* (\CO(\Delta), \|\cdot\|_{\Delta})
$$
as metrized line bundles on $C_{K'}$.
Here we view $x$ as a closed point of $X_{K'}$,  
$\CO(x)$ is the line bundle on $X_{K'}$ corresponding to $x\in X_{K'}$, and
$$i_{x}=(x, \mathrm{id}): \Spec K'\times C \lra C\times C$$
is the natural morphism.
It follows that 
$$
g_x=-\log\|1\|_x:(C_{K'})^\an \setminus \{x\}\lra \RR$$ 
is equal to 
the pull-back of $g_\Delta$ via the map
$i_{x}^\an:(C_{K'})^\an \to (C^2)^\an$.
We may also write $g_x=g_\Delta(x,\cdot)$ by abuse of notations. 
Finally, we re-organize the construction of \cite{Zha1} as the following statement. 

\begin{thm} \label{admissible1}
Let $K$ be a non-archimedean field, and $C$ be a smooth projective curve of genus $g>0$ over $K$.
There is a unique pair $(\|\cdot\|_{\Delta,a}, \|\cdot\|_a)$, where 
\begin{enumerate}[(a)]
\item $\|\cdot\|_a$ is an integrable metric of $\omega_{C/K}$ on $C^\an$,
\item $\|\cdot\|_{\Delta,a}$ is a symmetric integrable metric of $\CO_{C^2}(\Delta)$ on $(C^2)^\an$,
\end{enumerate}
satisfying the following properties
 for all finite extensions $K'/K$ and all points $x,y\in C(K')$:
\begin{enumerate}[(1)]
\item the equalities
$$
c_1(\CO(x), \|\cdot\|_{x})=c_1(\CO(y), \|\cdot\|_{y})
$$
and
$$
(2g-2)c_1(\CO(x), \|\cdot\|_{x})
=c_1(\omega_{C_{K'}/K'}, \|\cdot\|_a),
$$
hold as Chambert-Loir measures on $(C_{K'})^\an$;
\item the integral
$$
\int_{(C_{K'})^\an} g_{\Delta,a}(x,\cdot)\, c_1(\CO(x), \|\cdot\|_{x})=0.
$$
\item the residue map
$$
(\omega_{C/K} \otimes_{\CO_C} \CO(x) )|_{x}\lra K'
$$
is an isometry, where $K'$ is endowed with the absolute value extending that of $K$.
\end{enumerate}
\end{thm}

Here the Chambert-Loir measure $c_1(\omega_{C_{K'}/K'}, \|\cdot\|_a)$ uses the base change of the metric 
$\|\cdot\|_a$ from to $\omega_{C/K}$ to $\omega_{C_{K'}/K'}$.
Note that if $g>1$, the first two properties are equivalent to the following:
\begin{enumerate}[(1)]
\item the equality
$$
(2g-2)c_1(\CO(x), \|\cdot\|_{x})
=c_1(\omega_{C_{K'}/K'}, \|\cdot\|_a)
$$
holds as Chambert-Loir measures on $(C_{K'})^\an$;
\item the integral
$$
\int_{(C_{K'})^\an} g_{\Delta,a}(x,\cdot) c_1(\omega_{C_{K'}/K'}, \|\cdot\|_a)=0.
$$
\end{enumerate}

We call $\|\cdot\|_{a}$ the \emph{Zhang metric} of $\omega_{C/K}$ on $C^\an$,
and call $\|\cdot\|_{\Delta,a}$  the \emph{Zhang metric} of $\CO(\Delta)$ on $(C^2)^\an$. We also call them the \emph{canonical admissible metrics}.

Note that the statement of the theorem also works for archimedean $K$, which concerns the Arakelov metrics recalled in the previous subsection. 
Moreover, our proof of the theorem in the following works in both the archimedean case and the non-archimedean case. For simplicity, we will still assume that $K$ is non-archimedean in the following.

We provide two approaches of the uniqueness in the theorem.
For the first approach, by base change, we can assume $C(K)\neq \emptyset$ and take a point $x_0\in X(K)$.
Assume that there is another pair $(\|\cdot\|_{\Delta,a}', \|\cdot\|_a')$ of metrics satisfying the properties, and denote $\|\cdot\|_{x}'$ similarly.  
Write $\|\cdot\|_{x_0}'=\|\cdot\|_{x_0}e^{\varphi}$ for a continuous function $\varphi:C^\an\to \RR$.
By property (1), 
$$
c_1(\CO_{C}, \|\cdot\|_a'/\|\cdot\|_a)=
(2g-2)c_1(\CO_{C}, \|\cdot\|_{x_0}'/\|\cdot\|_{x_0})=
c_1(\CO_{C}, e^{(2g-2)\varphi})
$$
and 
$$
c_1(\CO_{C_{K'}}, \|\cdot\|_{x}'/\|\cdot\|_{x})=
c_1(\CO_{C_{K'}}, \|\cdot\|_{x_0}'/\|\cdot\|_{x_0})=
c_1(\CO_{C_{K'}}, e^{\varphi}).
$$
By the non-archimedean Calabi theorem in \cite[Cor. 2.2]{YZ1}, 
$$\|\cdot\|_{a}'/\|\cdot\|_{a}=e^{(2g-2)\varphi+c}$$ 
over $C^\an$ for some constant $c\in \RR$,
and 
$$
\|\cdot\|_{x}'/\|\cdot\|_{x}=e^{\psi(x)} e^{\varphi}$$ 
over $C_{K'}^\an$ for some function $\psi:C(\overline K)\to \RR$. 
This gives 
$$g_\Delta(x,y)-g_\Delta'(x,y)=\psi(x)+\varphi(y),\quad
x,y\in C(\overline K).$$ 
Since $g_\Delta'$ and $g_\Delta$ are symmetric, we have $\psi=\varphi+c'$
for a constant $c'\in\RR$. 
By property (3), 
$$
\frac{\|\cdot\|_a'(x)}{\|\cdot\|_a(x)} \cdot \frac{\|\cdot\|_{x}'(x)}{\|\cdot\|_{x}(x)}=1,\qquad
\frac{\|\cdot\|_a'(x)}{\|\cdot\|_a(x)} \cdot e^{\psi(x)} e^{\varphi(x)}=1,
$$
and thus
$$
\|\cdot\|_a'=\|\cdot\|_ae^{-2\varphi-c'}.
$$
Compare with $\|\cdot\|_{a}'/\|\cdot\|_{a}=e^{(2g-2)\varphi+c}$ (for $g>0$). 
We see that $\varphi$ is a constant on $C^\an$.  
It follows that $\|\cdot\|_{x}'/\|\cdot\|_{x}$ is a constant.
Then 
$$c_1(\CO(x), \|\cdot\|_{x}')=c_1(\CO(x), \|\cdot\|_{x}).$$ 
Property (2) implies $\|\cdot\|_{x}'=\|\cdot\|_{x}$, and property (3) implies 
$\|\cdot\|_{a}'=\|\cdot\|_a$. 
This proves the uniqueness.
 
Our second proof of the uniqueness in the theorem follows from \cite[Thm. 4.6]{Zha1}.
This also gives the connection to the construction in the loc. cit. 
In fact, for any $D,E\in \Div(C_{\overline K})$ with disjoint supports, define
$$
(D,E)_a=g_{\Delta,a}(D,E),
$$
by bi-linearity.
We can check that this pairing satisfies all properties of \cite[Thm. 4.6]{Zha1}, and thus the pairing $g_{\Delta,a}: C(\overline K)^2\setminus \Delta(\overline K)\to \RR$ 
is unique up to additive constants. 
By continuity, the function $g_{\Delta,a}: (C^2)^\an\setminus \Delta^\an\to \RR$ is unique to additive constants. 
This determines $c_1(\omega_{C/K},\|\cdot\|_a)$ uniquely by property (1).
Then $g_x$ is uniquely determined by property (2), and $\|\cdot\|_a$ is uniquely determined by property (3).
Note that \cite[Thm. 4.6]{Zha1} assumes that $K$ has a discrete valuation, but it also holds for non-discrete valuations.

There are also two approaches of the existence in the literature. 
The first one is in terms of reduction graphs, which is worked out in details by Zhang \cite{Zha1}. 
The second one is in terms of admissible metrics of line bundles on the Jacobian variety, which is outlined in \cite{Zha1}, and worked out partially by Heinz \cite{Hei}. 
The first approach has the advantage of being explicit, while the second approach has the advantage of being functorial and includes all (non-classical) points of the Berkovich spaces. 
We will prove the theorem by the second approach, and then relate it to the first approach.

\subsection{Admissible metrics} \label{sec appendix admissible}

Let $K$ be a non-archimedean field. As in \cite[\S3]{Hei}, we introduce admissible metrics of line bundles on curves (resp. square of curves) by means of Jacobian varieties. 
 
Let $A$ be an abelian variety over $K$, and $M$ be a line bundle on $A$. 
Let $\|\cdot\|$ be a continuous metric of $M$ on $A^\an$, and denote $\OM=(M,\|\cdot\|)$. 
Define admissibility in the following three cases:
\begin{enumerate}[(1)]
\item If $M$ is even in the sense that $[-1]^*M\simeq M$, there is an isomorphism $[2]^*M \otimes M^{\otimes (-4)}\simeq \CO_A$. 
The metric $\|\cdot\|$ is called \emph{admissible} if there is a constant metric on $\CO_A$ such that the isomorphism is an isometry. 

\item If $M$ is odd in the sense that $[-1]^*M\simeq M^{\otimes (-1)}$,  there is an isomorphism $[2]^*M \otimes M^{\otimes (-2)}\simeq \CO_A$. 
The metric $\|\cdot\|$ is called \emph{admissible} if there is a constant metric on $\CO_A$ such that the  isomorphism is an isometry. 

\item In general, the metric $\|\cdot\|$ of $M$ is called \emph{admissible} if the induced metrics on the even line bundle
$M\otimes [-1]^*M$ and the odd line bundle $M\otimes [-1]^*M^{\otimes (-1)}$ are both admissible. 
\end{enumerate}

By Tate's limiting argument in \cite{Zha2}, any line bundle on $A$ has an admissible metric, which is unique up to multiplicative constants. 
Note that all these arguments were originally written for metrics on $A(\overline K)$ instead of on $A^\an$, but they can be modified to $A^\an$ without essential difficulty.
Moreover, admissible metrics are integrable for general $L$ and semipositive for ample $L$.

Let $C$ be a smooth projective curve of genus $g>0$ over $K$.   
Let $\alpha\in \Pic^1(C)$ be a line bundle of degree 1 on $C$. 
If $\alpha$ does not exist, we will need to pass to a finite extension of $K$, and we will come back to this issue later. 
Denote by $J=\Pic^0_{C/K}$ the Jacobian variety of $C$ over $K$.
Denote the canonical embedding
$$i_\alpha: 
C\lra J,\quad x\longmapsto (x)-\alpha.
$$
The theta divisor $\theta_\alpha\subset J$ is the image of the morphism 
$$C^{g-1}\lra J, \quad (x_1,\cdots, x_{g-1})\longmapsto i_\alpha(x_1)+\cdots+ i_\alpha(x_{g-1}).$$

Let $L$ be a line bundle on $C$. 
Following \cite[\S4]{Hei}, a continuous metric $\|\cdot\|$ of $L$ on $C^\an$  
is called \emph{admissible} if there exist a line bundle 
$M$ on $J$, and an admissible metric $\|\cdot\|_M$ of $M$ on $J^\an$, such that $M$ is algebraically equivalent to an integer multiple of $\theta_\alpha$ on $J$, and such that
 $(L,\|\cdot\|)^{\otimes m}$ is isometric to $i_\alpha^*(M,\|\cdot\|_M)$ for some positive integer $m$. 

This definition looks very random, but any line bundle $L$ on $C$ has an admissible metric, which is unique up to multiplicative constants. 
Moreover, the definition is independent of the choice of $\alpha$ and stable under base change. 

In the case $\Pic^1(C)= \emptyset$, we say that a metric $\|\cdot\|$ of $L$ on $C^\an$ is \emph{admissible} if the induced metric of $L_{K'}$ on $(C_{K'})^\an$
is admissible for some finite extension $K'/K$ of $K$ with $\Pic^1(C_{K'})\neq \emptyset$. Then in this general case, any line bundle $L$ on $C$ still has an admissible metric, unique up to multiplicative constants. 

We have the following interpretation in terms of Monge--Amp\`ere equations, where an explicit form of $d\mu_a$ will be given in Proposition \ref{measure2}. 

\begin{prop} \label{measure1}
There exists a unique probability measure $d\mu_a$ on $C^\an$ such that for any metrized line bundle $(L,\|\cdot\|)$ on $C$, the metric $\|\cdot\|$ is admissible if and only if it is integrable and satisfies
$$
c_1(L,\|\cdot\|) = \deg(L)\, d\mu_a
$$
for the Chambert-Loir measure on $C^\an$.
\end{prop}
\begin{proof}
It is easy to reduce the problem to the case $\Pic^1(C)\neq \emptyset$. Then we take $\alpha\in \Pic^1(C)$.

Denote  
$(L,\|\cdot\|)=i_\alpha^*(\CO(\theta_\alpha),\|\cdot\|_\alpha)$, where $\|\cdot\|_\alpha$ 
is an admissible metric of $\CO(\theta_\alpha)$ on $J^\an$.
Define the measure $d\mu_a$ by the equation 
$$
c_1(L,\|\cdot\|) = \deg(L)\, d\mu_a.
$$

To prove the equation for any admissible $(L,\|\cdot\|)$, it suffices to prove that if $M$ is a line bundle on $J$ algebraically equivalent to 0, then the measure
$c_1(i_\alpha^*(M,\|\cdot\|_M))=0$ for any admissible metric $\|\cdot\|_M$ of $M$ on $J^\an$. 
This is exactly \cite[Thm. 5.16]{YZ1}, which is essentially a result of Gubler. 

It remains to prove that any metric satisfying the equation is admissible. 
Note that admissible metrics already satisfy the equation, so it suffices to prove that 
the equation determines the metric up to a constant multiple. 
This follows from the non-archimedean Calabi theorem in \cite[Cor. 2.2]{YZ1}.
We only need the theorem for curves, in which case the ampleness assumption can be removed by additivity. 
\end{proof}

Let $L$ be a line bundle on $C^2$. 
A continuous metric of $L$ on $(C^2)^\an$ is called \emph{admissible} if its pull-back metrics to $(C_{K'})^\an$ via all morphisms $(x,\id):C_{K'}\to C\times C$ and $(\id,x):C_{K'}\to C\times C$, for all points $x\in C(K')$ over all finite extensions $K'/K$, are admissible. 
This definition follows from \cite[Def. 4.2]{Hei}, except that the loc. cit. used the term ``bi-admissible'' instead of ``admissible''.

Any line bundle on $C^2$ admits an admissible metric, unique up to multiplicative constants.
The uniqueness is an easy consequence of the definition.  
For the existence, it suffices to assume that $C(K)$ contains an element $x_0$ by extending $K$ if necessary.
Note that admissible metrics of line bundles in $p_1^*\Pic(C)$ and $p_2^*\Pic(C)$ can be obtained as pull-back of admissible metrics.
One is reduced to treat the subgroup $\Pic^-(C^2)$ of $\Pic(C^2)$ consisting of line bundles on $C^2$ whose restrictions to $C\times\{x_0\}$ and $\{x_0\}\times C$ are both trivial. 
Denote by $\Pic^-(J^2)$ the subgroup of $\Pic(J^2)$ consisting of line bundles on $J^2$ whose restrictions to $J\times\{0\}$ and $\{0\}\times J$ are both trivial. 
It turns out that the canonical map $\Pic^-(J^2)\to \Pic^-(C^2)$ is an isomorphism, and thus elements of $\Pic^-(C^2)$ have admissible metrics by pull-back of admissible metrics from $J^2$.
See \cite[Lem. 2.2.1, Lem. 2.2.2, Lem. 2.2.3]{Zha3} for all details in this construction. 

By our construction below, we will see that in Theorem \ref{admissible1},
the metrics $\|\cdot\|_a$ and $\|\cdot\|_{\Delta,a}$ are both admissible.

\subsection{Construction of the Zhang metrics}

The goal is to construct the Zhang metrics in Theorem \ref{admissible1} in terms of Jacobian varieties. 

For the existence in the theorem, we remark that for any finite Galois extension $K'/K$, if a pair $(\|\cdot\|_{\Delta,a}, \|\cdot\|_a)$ exists over $K'$, then it is Galois invariant by the uniqueness, and thus descends to $K$.
Therefore, by replacing $K$ by a finite Galois extension, we can assume that $C$ has some convenient properties. For example, we can assume that $C(K)\neq \emptyset$. 

We will need the following theorem here, and it will also be used in the main part of this article. 
For clarity, we will state all the notations and conditions of this theorem separately. 

\begin{thm} \label{isomorphism}
Let $K$ be any field, and $C$ be a smooth projective curve of genus $g>0$ over $K$. Denote by $J$ the Jacobian variety of $C$ over $K$.
Let $\alpha$ be a divisor on $C$ of degree 1.
Let $i_\alpha$ be the canonical embedding
$$i_\alpha: 
C\lra J,\quad x\longmapsto (x)-\alpha.
$$
Let $\theta_\alpha\subset J$ be the theta divisor, i.e., the image of the morphism 
$$C^{g-1}\lra J, \quad (x_1,\cdots, x_{g-1})\longmapsto i_\alpha(x_1)+\cdots+ i_\alpha(x_{g-1}).$$ 
Denote by $P$ the Poincar\'e line bundle on $J\times J$. 
Then the following are true:
\begin{enumerate}[(1)]
\item There are isomorphisms
$$i_\alpha^*\CO([-1]^*\theta_\alpha)= g\alpha,$$
$$i_\alpha^*\CO(\theta_\alpha)= \omega_{C/K}+(2-g)\alpha.$$
\item There is an isomorphism
$$P=  m^*\CO(-\theta_\alpha)\otimes p_1'^*\CO(\theta_\alpha)\otimes p_2'^*\CO(\theta_\alpha),$$
where $m:J\times J \to J$ is the addition law, and $p_1', p_2':J\times J\to J$ are the projections.
\item Denote by $(i_\alpha,i_\alpha):C\times C\to J\times J$ the natural morphism, and by  $p_1, p_2:C\times C\to C$ the projections. 
Then there is an isomorphism
$$
(i_\alpha,i_\alpha)^*P = \CO(\Delta)-p_1^*\alpha-p_2^*\alpha,
$$ 
where $\Delta\subset C\times C$ is the diagonal. 
\item Denote by $\Delta_J:J\to J\times J$ the diagonal morphism. Then 
 there is an isomorphism
$$
\Delta_J^*P = \CO(-(\theta_\alpha+[-1]^*\theta_\alpha)).
$$ 
\end{enumerate}
\end{thm}

\begin{proof}
One can derive these results from either \cite[\S5.6]{Ser} or \cite[\S A.8.2]{HS}. We will follow 
\cite[\S5.6]{Ser} for example.
Note that our $(\theta_\alpha, [-1]^*\theta_\alpha)$ is $(\Theta,\Theta')$ in \cite[\S5.6]{Ser} under the relation $\alpha=a$.
The first equation in \cite[p. 75]{Ser} gives
$$
i_{\alpha-c}^*([-1]^*\theta_\alpha)=g\alpha-c, \quad \forall c\in \Pic^0(C). 
$$
It gives the first isomorphism of (1) by taking $c=0$,
and it implies the second isomorphism of (1) by combining with \cite[p. 74, eq. (1)]{Ser}.
Part (2) of our theorem is exactly \cite[p. 76, eq. (3)]{Ser}.
Part (3) of our theorem is \cite[p. 76, eq. (4)]{Ser}, which assumes $(2g-2)\alpha=\omega_{C/K}$, but the proof works without this assumption by some extra computation using the first equation in \cite[p. 75]{Ser}.
The result can also be checked using the universal property of $P$.
For  (4), note that (2) implies
$$\Delta_J^*P= [2]^*\CO(-\theta_\alpha)\otimes \CO(2\theta_\alpha).$$ 
The result follows from the basic formula
$$
[2]^*\CO(\theta_\alpha)= \CO(3\theta_\alpha+[-1]^*\theta_\alpha).
$$
See the corollary in \cite[p. 33]{Ser}. 
\end{proof}

Part (1) of the lemma gives a concrete way to construct admissible metrics of $\alpha$ and $\omega_{C/K}$ on $C^\an$. 
Part (2) gives a concrete way to construct admissible metrics of $\CO(\Delta)$ on $(C^2)^\an$.

Now we are ready to prove the existence in Theorem \ref{admissible1}. 
We will prove that there is a unique admissible metric $\|\cdot\|_a$ of $\omega_{C/K}$ and a unique admissible metric $\|\cdot\|_{\Delta,a}$ of $\CO(\Delta)$ satisfying the requirements.
For that purpose, fix an admissible metric $\|\cdot\|_a'$ of $\omega_{C/K}$ and an admissible metric $\|\cdot\|_{\Delta}'$ of $\CO(\Delta)$.
We will find scalar multiples of these two metrics to satisfy the conditions. 
By Proposition \ref{measure1}, part (1) of the theorem holds for $(\|\cdot\|_\Delta',\|\cdot\|_a')$. 
It suffices to prove that (2) holds for $(\|\cdot\|_\Delta',\|\cdot\|_a')$ up to an additive constant independent of $x$, and that 
(3) holds for $(\|\cdot\|_\Delta',\|\cdot\|_a')$ up to a multiplicative constant independent of $x$.
Then we can easily modify $\|\cdot\|_\Delta'$ to satisfy (2) and modify $\|\cdot\|_a'$ to satisfy (3).

We first prove that (3) holds for $(\|\cdot\|_\Delta',\|\cdot\|_a')$ up to a constant. 
By definition, there is a canonical isomorphism
$$
\omega_{C/K}\otimes_{\CO_C} (\CO(\Delta)|_\Delta)  \lra \CO_C.
$$
We claim that the isomorphism is an isometry up to constants under the metrics $(\|\cdot\|_\Delta',\|\cdot\|_a')$ and the trivial metric of $\CO_C$.
By uniqueness of admissible metrics, it suffices to prove that $(\CO(\Delta)|_\Delta,\|\cdot\|_\Delta')$ is admissible on $C$. 
This is \cite[Prop. 4.3]{Hei}. 
It is also easy to prove it in our setting. In fact,
construct an admissible metric on $\CO(\Delta)$ by Theorem \ref{isomorphism}(3).
Note that the composition 
$$C\stackrel{\Delta}{\lra} C^2 \stackrel{(i_\alpha,i_\alpha)}{\lra} J^2$$ 
is equal to the composition 
$$C\stackrel{i_\alpha}{\lra} J \stackrel{\Delta_J}{\lra} J^2,$$ 
where $\Delta_J:J\to J^2$ is the diagonal morphism. 
Thus it suffices to check that $\Delta_J^*P$ is algebraically equivalent to a multiple of 
$\theta_\alpha$, but this follows from Theorem \ref{isomorphism}(4).

Therefore, the isomorphism
$$
\omega_{C/K}\otimes_{\CO_C} (\CO(\Delta)|_\Delta)  \lra \CO_C
$$
is an isometry up to a constant multiple. 
Take the pull-back of this isometry via $x:\Spec K'\to C$, we get an isomorphism
$$
x^*\omega_{C/K}\otimes_{K'} x^*(\CO(\Delta)|_\Delta)  \lra K',
$$
which is an isometry up to a constant independent of $x$.
Note that the composition 
$$\Spec K' \stackrel{x}{\lra} C \stackrel{\Delta}{\lra} C^2$$ 
is equal to the composition 
$$\Spec K' \stackrel{x}{\lra} C_{K'} \stackrel{(x,\id)}{\lra} C^2,$$ 
which induces canonical isomorphisms 
$$x^*(\CO(\Delta)|_\Delta)\lra x^*\CO((x,\id)^*\Delta)\lra  x^*\CO(x).$$ 
Under this composition, the above isomorphism becomes exactly the residue map. 
This proves that (3) holds for $(\|\cdot\|_\Delta',\|\cdot\|_a')$ up to a constant.

It remains to prove that (2) holds for $(\|\cdot\|_\Delta',\|\cdot\|_a')$ up to constants. 
Let $L$ be a line bundle on $C$ endowed with an admissible metric. 
There are canonical isomorphisms 
$$
p_{1*}\pair{\CO(\Delta), p_2^*L} \lra \Delta^*p_2^*L\lra L,
$$
where the left-hand side is the Deligne pairing with respect to the morphism 
$p_1:C^2\to C$. 
The metrics of $\CO(\Delta)$ and $p_2^*L$ induce a canonical continuous 
metric of the Deligne pairing. 
In the complex case, this process is explicitly treated in \cite{Del,Elk} and \cite[\S4.2]{YZ2}. The non-archimedean case is written in \cite[\S4.6.2]{YZ2}. 

With the natural metric of the Deligne pairing, we claim that the isomorphism 
$$
p_{1*}\pair{\CO(\Delta), p_2^*L} \lra L
$$
is an isometry up to a constant. 
If this holds, taking the base change of $p_1:C^2\to C$ by the morphism $x:\Spec K'\to C$, we have the structure morphism $\pi':C_{K'}\to \Spec K'$.
The Deligne pairing under this base change gives an isomorphism 
$$
\pi'_* \pair{\CO(x), L_{K'}} \lra L(x),
$$
which is an isometry up to a constant. 
By definition of the metrics, the logarithm of the norm of this isomorphism is actually
$$
-\int_{(C_{K'})^\an} g_{x}' c_1(L,\|\cdot\|_L).
$$
The constancy of this integral with $L=\omega_{C/K}$ is exactly what we need for (2) of Theorem \ref{admissible1}. 

It remains to prove the claim that the isomorphism 
$$
p_{1*}\pair{\CO(\Delta), p_2^*L} \lra L
$$
is an isometry on $C$ up to a constant. 
By the uniqueness of admissible metrics, it suffices to prove that the metric of the 
left-hand side is admissible. 

As a consequence of \cite[Thm. 4.6.2]{YZ2}, the Deligne pairings and the metrics are compatible with base change.
To apply this compatibility, denote by $q_1:J\times C\to J$ and $q_2:J\times C\to C$ the projections.
View $p_1:C^2\to C$ as the base change of $q_1:J\times C\to J$ by $i_\alpha:C\to J$. 
Denote by $Q=(\id, i_\alpha)^* P$ the pullback of the Poincar\'e bundle $P$ via 
$(\id, i_\alpha):J\times C\to J\times J$.
By Theorem \ref{isomorphism}(3), there are isomorphisms 
$$
(i_\alpha,\id)^*Q\lra 
(i_\alpha,i_\alpha)^*P\lra 
\CO(\Delta)-p_1^*\alpha-p_2^*\alpha.
$$
Then the compatibility with base change gives an isometry (up to a constant)
$$
p_{1*}\pair{\CO(\Delta)-p_1^*\alpha-p_2^*\alpha,\ p_2^*L}
\lra i_\alpha^*\big(q_{1*}\pair{Q, q_2^*L}\big),
$$
where $\alpha$ is endowed with an admissible metric, and $Q=(\id, i_\alpha)^* P$ is endowed with the pull-back of an admissible metric of $P$.

On the other hand, by a local version of \cite[Lem. 4.6.1(2)]{YZ2}, 
$p_{1*}\pair{p_1^*\alpha,\ p_2^*L}$
is isometric to $\alpha^{\otimes \deg(L)}$, and thus admissible. 
Moreover,
$p_{1*}\pair{p_2^*\alpha,\ p_2^*L}$ is isometric to $\CO_C$ with a constant metric. 
To see this, by additivity and base change, we can assume that $\alpha=\CO(x_0)$ for a point $x_0\in C(K)$. Then we can apply the integration formula in \cite[\S4.6.2]{YZ2}
to get the result.

Therefore, it suffices to prove that $M=q_{1*}\pair{Q, q_2^*L}$ is algebraically equivalent to 0, and its metric is admissible on $J$.
Consider the base change of $q_1:J\times C\to J$ by $[2]:J\to J$. 
This gives an isometry of Deligne pairings
$$
[2]^* q_{1*}\pair{Q, q_2^*L} \lra q_{1*}\pair{[2]_C^* Q, [2]_C^*q_2^*L}.
$$
Here $[2]_C: J\times C \to J\times C$ is the base change of $[2]: J \to J$. 
Note that there is a natural isometry $[2]_C^*q_2^*L\simeq q_2^*L$.
We also have an isometry (up to constant) $[2]_C^*Q\simeq Q^{\otimes 2}$, which can be obtained by the argument right before \cite[Lem. 6.5.3]{YZ2}. 
It follows that there is an isometry (up to constant) 
$[2]^*M \simeq M^{\otimes 2}$ on $J$. 
This proves that $M$  is algebraically equivalent to 0, and its metric is admissible on $J$.

Therefore, part (2) of Theorem \ref{admissible1} also holds up to multiplicative constants.
Then the proof of the theorem is complete.

\subsection{Zhang's construction by graph theory}
\label{sec graph theory}

Now we review the construction of \cite{Zha1} to give explicit descriptions of the metrics in Theorem \ref{admissible1}.

Let $K$ be a non-archimedean field with a discrete valuation.
Denote $e_K=|\varpi|^{-1}$, where $\varpi$ is a generator of the maximal ideal of the valuation ring $O_K$. 

Let $C$ be a smooth projective curve of genus $g>0$ over $K$.
As before, we can replace $K$ by a finite Galois extension in the definitions. 
Then we assume that $C(K)\neq \emptyset$, and $C$ has \emph{split semistable reduction} over $O_K$, i.e., the minimal regular model $\CCC$ of $C$ over $O_K$ has a semistable special fiber $\CCC_s$ over the residue field $\kappa$ of $O_K$, every node $x$ of $\CCC_s$ is defined over $\kappa$, and the two tangent lines of $\CCC_s$ at every node $x$ are also defined over $\kappa$. 

Denote by $\|\cdot\|_{\Ar}$ the metric of $\omega_{C/K}$ on $C^\an$ induced by the integral model $(\CCC, \omega_{\CCC/O_K})$. 
In the following, we will modify $\|\cdot\|_{\Ar}$ to get $\|\cdot\|_{a}$.

Denote by $\Gamma=\Gamma(C)$ the \emph{reduction graph} of $C$, which is the dual graph of $\CCC_s$. 
Then every irreducible component $D$ of $\CCC_s$ is represented by a vertex $v(D)$ of $\Gamma(C)$; every node $N$ of $\CCC_s$ is represented by an edge $e(N)$ in $\Gamma(C)$
connecting the vertices $v(D_1), v(D_2)$, where $D_1,D_2$ are the (possible equal) irreducible components of $\CCC_s$ containing $N$.   
By fixing an isomorphism from each edge $e(N)$ to a closed interval of length 1 in $\RR$, 
 $\Gamma(C)$ becomes a metrized graph. We refer to \cite{BF} and \cite[\S4.1]{Bak} for some basics of metrized graphs. 

The \emph{canonical divisor} $K_C$ of $\Gamma$ is the formal linear combination 
$$
K_C = \sum_{\xi\in V(\Gamma)} \deg(\omega_{\CCC/O_K}|_{F_\xi})\, \xi,
$$
where $V(\Gamma)$ denotes the vertex set of $\Gamma$, and $F_\xi$ denotes the irreducible component of $\CCC_s$ corresponding to $\xi$. 
It is easy to see that $\deg(K_C)=2g-2$.

Following \cite[Def. a.3]{Zha1}, 
denote by $F(\Gamma)$ the space of continuous and piecewise smooth functions $f:\Gamma\to \RR$ such that all the one-sided directional derivatives $d_{\vec{v}}f(P)$ of $f$ exist. 
The Laplacian operator $\Delta$ sends each $f\in F(\Gamma)$ to a signed measure on $\Gamma$ given by 
$$
\Delta f=-f''(x)dx- \sum_{P\in\Gamma} \sum_{\vec v\in T_P(\Gamma)} d_{\vec v}f(P) \,\delta_P.
$$
Here $x$ represents a canonical coordinate on each edge of $\Gamma$, which is uniquely defined up to a sign, and $T_P(\Gamma)$ denotes the (finite) set of tangent directions of $\Gamma$ at $P$.
The summation for $P\in\Gamma$ on the right-hand side has only finitely many nonzero terms. 
See also \cite[Def. 5]{BF} and \cite[Prop. 4.2.6]{Bak}.

By \cite[Thm. 3.2, Lem. 3.7]{Zha1}, there is a unique positive measure 
$\mu=\mu_{K_C}$ on $\Gamma$ of total volume 1, which is a finite linear combination of the uniform measures on the edges of $\Gamma$ and the Dirac measures supported on the vertices of $\Gamma$, such that the function 
$g_{\mu}: \Gamma^2 \to \RR$ uniquely defined by 
$$
g_{\mu}(x,\cdot) \in F(\Gamma),\quad \forall x\in \Gamma,
$$
$$
\Delta g_{\mu}(x,\cdot)=\delta_x-\mu,
$$
$$
\int_\Gamma g_{\mu}(x,\cdot) \mu_{K_C}=0,
$$
satisfies the condition that 
$$
c+g_{\mu}(K_C, x)+ g_{\mu}(x,x)=0
$$
for a constant $c$ independent of $x\in \Gamma$. 
Here we set $g_{\mu}(D, x)=\sum_{i=1}^r a_i g_{\mu}(x_i, x)$ for a divisor $D=\sum_{i=1}^r a_i x_i$ on $\Gamma$ with $a_i\in\RR$ and $x_i\in V(\Gamma)$. 
The function  $g_{\mu}: \Gamma^2 \to \RR$ is actually continuous and symmetric. 

By \cite[\S4.3]{Ber}, there is a canonical injection
$$
i:  \Gamma(C)\lra C^\an 
$$
and a canonical retraction map
$$
r: C^\an \lra \Gamma(C).
$$
By composition, it induces a retraction map 
$$
r: C(\overline K) \lra \Gamma(C).
$$
In terms of this map, we can view $g_{\mu}(K_C, x)$ and $g_{\mu}(x,x)$
as functions of $x\in C^\an$. 

Finally, we have explicit formulas for the metrics $\|\cdot\|_a$ and $\|\cdot\|_{\Delta,a}$ in Theorem \ref{admissible1}. 
Following \cite[(4.1)]{Zha1}, we have
$$
\|\cdot\|_a(x)
=\|\cdot\|_\Ar(x) \cdot e_K^{-c-g_{\mu}(K_C,r(x))}
=\|\cdot\|_\Ar(x) \cdot e_K^{g_{\mu}(r(x),r(x))},\quad x\in C^\an
$$
and
$$
g_{\Delta,a}(x,y)=i(x,y)\log e_K + g_{\mu}(r(x), r(y))\log e_K, \quad x,y\in C(\overline K),\, x\neq y.  
$$
Here for $x,y\in C(\overline K)$, if $K'$ is a finite extension of $K$ such that $x,y\in C(K')$, then $i(x,y)$ is the intersection number of the corresponding sections in the minimal regular model of $C$ over $O_{K'}$, divided by the normalizing factor $[K':K]$.

Note that the above expression for $\|\cdot\|_a$ is already defined over $C^\an$, but the expression for $g_{\Delta,a}$ is defined over the $\overline K$-points, which determines its value over $(C^2)^\an$ by continuity.  

To verify the above expressions, we first reduce to the case that $e_K=|\varpi|^{-1}$ is equal to $e$ (the base of the natural logarithm).
In fact, if we replace the absolute value $|\cdot|$ of $K$ by $|\cdot|^a$ for some $a>0$, then  the canonical admissible metrics are changed to their powers of exponent $a$ by the uniqueness of the metrics.

For the standard case $e_K=e$, it suffices to apply the last statement of \cite[Thm. 4.6]{Zha1}, which has also been used to prove the uniqueness part of Theorem \ref{admissible1}. 
In translating the measure between the graph and the Berkovich space, it is helpful to have the following result, which is compatible with \cite[Thm. 2.3]{Zha1}.

\begin{prop} \label{measure translation}
For any $f \in F(\Gamma)$, denote by $\CO_C(f)$ the trivial line bundle $\CO_C$ endowed with the metric on $C^\an$ given by $\|1\|=e_K^{-r^*f}$. 
Then  
$$
c_1(\CO_C(f))=-i_*(\Delta f)
$$
as measures on $C^\an$.
\end{prop}
\begin{proof}
We only sketch the idea here.
The base case is when $\CO(f)$ is induced by an integral model $\CO_\CCC(V)$ of 
$\CO_C$ for an irreducible component $V$ of the special fiber of the minimal regular model $\CCC$. 
In this case, both the Chambert-Loir measure and Laplacian operator are easy to compute. 
By linear combination, this gives the case that $f$ is linear on edges of $\Gamma$. 
Note that the edges of $\Gamma$ have rational lengths, and the division points of the edges are called rational points of $\Gamma$.
By taking finite extensions of $K$, we can realize all rational points of $\Gamma$ as irreducible components of the special fibers of the minimal regular models. 
Then the equality holds for all piecewise linear function $f$ on $\Gamma$ whose  critical points are rational points.
Such functions are uniformly dense in $F(\Gamma)$. 
Then the result is extended by approximation.  
\end{proof}

Finally, the measure $d\mu_a$ on $C^\an$ defined in Proposition \ref{measure1} has the following explicit expression. 
Recall that the resistance function of a metrized graph is defined right before \cite[Prop. 3.3]{Zha1} or with more details in \cite[\S6, Def. 8]{BF}, and it is used in our proof of Lemma \ref{effective1}.

\begin{prop} \label{measure2}
The measure 
$$
d\mu_a=\frac{1}{g}\, i_*\left( \sum_{v\in V(\Gamma(C))} g_v\delta_v
+\sum_{e\in E(\Gamma(C))} \frac{1}{r_e+1}\delta_e \right),
$$
where
\begin{itemize}
\item
$V(\Gamma(C))$ is the vertex set of $\Gamma(C)$ (corresponding to irreducible components of the special fiber $\CCC_s$);
\item $g_v$ is the genus of the normalization of the irreducible component of $\CCC_s$ corresponding to $v$;
\item $\delta_v$ is the Dirac measure supported at $v$;
\item $E(\Gamma(C))$ is the edge set of $\Gamma(C)$ (corresponding to nodes of the special fiber $\CCC_s$);
 \item
$r_e$ is the resistance of the endpoints $p,q$ of $e$ in $\Gamma(C)\setminus e^0$ with $e^0=e\setminus \{p,q\}$;
\item $\delta_e$ is the Lebesgue measure on the edge $e$ of total integral $1$.
\end{itemize}

\end{prop}
\begin{proof}
It suffices to compute 
$$
c_1(\omega_{C/K},\|\cdot\|_a)
=c_1(\omega_{C/K},\|\cdot\|_\Ar)
+c_1(\CO_C(c+g_{\mu}(K_C,\cdot))).
$$
By the definition of the Chambert-Loir measure in \cite{CL}, 
$$c_1(\omega_{C/K},\|\cdot\|_\Ar)=i_*\delta_{K_C}.$$
By Proposition \ref{measure translation},
$$
c_1(\CO_C(c+g_{\mu}(K_C,\cdot)))
=-i_*\Delta g_{\mu}(K_C,\cdot)
=i_*((2g-2)\mu-\delta_{K_C}).
$$
It follows that
$$
c_1(\omega_{C/K},\|\cdot\|_a)
=i_*((2g-2)\mu).
$$
Apply the explicit formula for $\mu=\mu_{K_C}$ and $D=K_C$ in \cite[Lem. 3.7]{Zha1}.
\end{proof}

\begin{remark}
The result is compatible with the abstract result in \cite[Thm. 1.1]{Gub}. 
\end{remark}

\subsection{General valuation fields} \label{sec all val}

Let $C$ be a smooth projective curve of genus $g>0$ over a complete valuation field $K$ (instead of a non-archimedean field).
In this generality, we will still have a \emph{canonical admissible metric} $\|\cdot\|_{a}$ of $\omega_{C/K}$ on $C^\an$, and a \emph{canonical admissible metric}
$\|\cdot\|_{\Delta,a}$ of $\CO(\Delta)$ on $(C^2)^\an$
satisfying Theorem \ref{admissible1} in suitable senses.

In fact, if $K$ is a non-archimedean field, then $\|\cdot\|_{a}$ and $\|\cdot\|_{\Delta,a}$
are the Zhang metrics introduced in Theorem \ref{admissible1}.

If $K$ is an archimedean field, then
$\|\cdot\|_{a}=\|\cdot\|_{\Ar}$ and $\|\cdot\|_{\Delta,a}=\|\cdot\|_{\Delta,\Ar}$
are the Arakelov metrics in \S \ref{sec ara}.
Note that if $K$ is real, there should be a minor process to descend $\|\cdot\|_{\Ar}$ from $C(\ol K)$ to 
$$C^\an=C(\ol K)/\Gal(\ol K/K).$$ 
A similar process is also needed for $\|\cdot\|_{\Delta,\Ar}$. 
Our statement and proof of Theorem \ref{admissible1} also work in the archimedean case.

If $K$ is a trivially valued field, then we set the metrics $\|\cdot\|_{a}$ and
$\|\cdot\|_{\Delta,a}$ to be the \emph{canonical metrics} of the line bundles.
Recall that for a projective variety $X$ over $K$ and a line bundle $L$ on $X$, the line bundle $L$ induces a metric of itself on $X^\an$ by the usual way, called the \emph{canonical metric} or the \emph{admissible metric}; see \cite[\S3.4.1]{YZ2} for example. Note that in this case, any line bundle $L$ has a \emph{unique} admissible metric (instead of unique up to multiplicative constants).

Let $K$ be either a number field or the function field of one variable over a base field. 
With the above metrics over the complete fields, we have an adelic line bundle
$$
\overline\omega_{C/K,a}
=(\omega_{C/K}, \{\|\cdot\|_{a,v}\}_v)
$$
over $C$, 
and an adelic line bundle
$$
\overline\CO(\Delta)_a
=(\CO(\Delta), \{\|\cdot\|_{\Delta,a,v}\}_v)
$$
over $C^2$. 

To check the metrics are adelic, we need to check that in Theorem \ref{admissible1}, if $C$ has good reduction over $O_K$, then the metrics 
$\|\cdot\|_a$ and $\|\cdot\|_{\Delta,a}$ are induced by the natural integral models $\omega_{\CCC/O_K}$ and $\CO_{\CCC^2}(\Delta_\CCC)$ of $\omega_{C/K}$ and $\CO(\Delta)$. Here $\CCC$ is the smooth projective model of $C$ over $O_K$, and  $\Delta_\CCC$ denotes the diagonal morphism $\CCC\to \CCC^2=\CCC\times_{O_K}\CCC$.
By uniqueness, it suffices to check that the induced metrics satisfy the properties of the theorem. 
We omit it here.

\end{document}